\documentclass[11pt]{article}
\usepackage{amsmath, amssymb, amsthm, mathtools}
\usepackage{bm}
\usepackage{enumerate}
\usepackage[left=1in,right=1in,top=1in,bottom=1in]{geometry}
\usepackage[backend=biber,bibencoding=utf8,maxnames=100
,doi=false,isbn=false,url=true,eprint=false,giveninits=true,date=year,sortcites=true]{biblatex}
\renewbibmacro{in:}{}
\usepackage{subfiles}
\usepackage{subcaption}
\usepackage{caption}
\usepackage{psfrag}
\usepackage{graphicx}
\usepackage{xkeyval}
\usepackage{url}
\usepackage{xcolor}
\usepackage{authblk}

\date{}
\title{A Hybrid Semi-Lagrangian Cut Cell Method for Advection-Diffusion Problems with Robin Boundary Conditions in Moving Domains}
\author[1,*]{Aaron Barrett}
\author[2]{Aaron L. Fogelson}
\author[3,4,5,6]{Boyce E. Griffith}
\affil[1]{\footnotesize Department of Mathematics, University of Utah, Salt Lake City, UT, USA}
\affil[2]{\footnotesize Departments of Mathematics and Bioengineering, University of Utah, Salt Lake City, UT, USA}
\affil[3]{\footnotesize Departments of Mathematics, Applied Physical Sciences, and Biomedical Engineering, University of North Carolina, Chapel Hill, NC, USA}
\affil[4]{\footnotesize Carolina Center for Interdisciplinary Applied Mathematics, University of North Carolina, Chapel Hill, NC, USA}
\affil[5]{\footnotesize Computational Medicine Program, University of North Carolina, Chapel Hill, NC, USA}
\affil[6]{\footnotesize McAllister Heart Institute, University of North Carolina, Chapel Hill, NC, USA}
\affil[*]{\footnotesize barrett@math.utah.edu}
\providecommand{\keywords}[1]
{
  \small	
  \textbf{\textit{Keywords---}} #1
}

\newcommand{\size}[1]{\left\vert #1\right\vert}

\makeatletter
\def\gsh#1{%
  \vbox{\hbox{%
    \let\\\cr
    \offinterlineskip
    \valign{&\hb@xt@2\p@{\hss$##$\hss}\vskip.2ex\cr#1\crcr}%
  }\vskip-.36ex}%
}
\def\gshsym{\@ifstar\gsh@ssym\gsh@sym}
\def\gsh@sym#1#2{\mathrlap{\overset{#1}{\phantom{#2}}}#2}
\def\gsh@ssym#1#2{\overset{#1}{#2}{\vphantom{#2}}}
\makeatother
\DeclarePairedDelimiter\norm{\lVert}{\rVert}%

\newcommand{\uu}{\mathbf{u}}
\newcommand{\grad}{\nabla}
\newcommand{\Lap}{\Delta}
\newcommand{\ff}{\mathbf{f}}
\newcommand{\parens}[1]{\mathopen{}\left(#1\right)\mathclose{}}
\newcommand{\xx}{\mathbf{x}}

\newcommand{\XX}{\mathbf{X}}

\newcommand{\bchi}{\bm{\chi}}

\newcommand{\nn}{\mathbf{n}}

\newcommand{\BB}{\mathcal{B}}

\newcommand{\bigOh}[1]{\mathcal{O}\left(#1\right)}

\newcommand{\cc}{\mathbf{c}}
\newcommand{\rr}{\mathbf{r}}
\newcommand{\etal}{et al.}
\renewcommand{\AA}{\mathbf{A}}

\makeatletter
\newlength{\sfp@hseplen}\newlength{\sfp@vseplen}
\define@cmdkey{subfigpos}[sfp@]{pos}[ul]{}
\define@cmdkey{subfigpos}[sfp@]{font}[\small]{}
\define@cmdkey{subfigpos}[sfp@]{vsep}[2\baselineskip]{\setlength{\sfp@vseplen}{\sfp@vsep}}
\define@cmdkey{subfigpos}[sfp@]{hsep}[10pt]{\setlength{\sfp@hseplen}{\sfp@hsep}}
\newcommand{\subfigimg}[3][,]{%
  \setkeys{Gin,subfigpos}{pos,font,vsep,hsep,#1}
  \setbox1=\hbox{\includegraphics{#3}}
  \ifnum\pdfstrcmp{\sfp@pos}{ul}=0
    \leavevmode\rlap{\usebox1}
    \rlap{\hspace*{\sfp@hsep}\raisebox{\dimexpr\ht1-\sfp@vsep}{\sfp@font{#2}}}
    \phantom{\usebox1}
  \else\ifnum\pdfstrcmp{\sfp@pos}{ur}=0
    \leavevmode\usebox1
    \llap{\raisebox{\dimexpr\ht1-\sfp@vsep}{\sfp@font{#2}}\hspace*{\sfp@hsep}}
  \else\ifnum\pdfstrcmp{\sfp@pos}{lr}=0
    \leavevmode\usebox1
    \llap{\raisebox{\sfp@vsep}{\sfp@font{#2}}\hspace*{\sfp@hsep}}
  \else
    \leavevmode\rlap{\usebox1}
    \rlap{\hspace*{\sfp@hseplen}\raisebox{\sfp@vsep}{\sfp@font{#2}}}
    \phantom{\usebox1}
  \fi\fi\fi
}
\makeatother
\begin{document}
\maketitle
\begin{abstract}
We present a new discretization approach to advection-diffusion problems with Robin boundary conditions on complex, time-dependent domains. The method is based on second order cut cell finite volume methods introduced by Bochkov \etal ~\cite{Bochkov2019} to discretize the Laplace operator and Robin boundary condition. To overcome the small cell problem, we use a splitting scheme along with a semi-Lagrangian method to treat advection. We demonstrate second order accuracy in the $L^1$, $L^2$, and $L^\infty$ norms for both analytic test problems and numerical convergence studies. We also demonstrate the ability of the scheme to convert one chemical species to another across a moving boundary.
\end{abstract}
\keywords{Irregular domain, Level set method, Robin boundary condition, Cartesian grid method}
\section{Introduction}

The convective transport and diffusion of chemical species occurs in a broad range of systems. Many applications involve chemical concentrations that evolve within complex, time-dependent regions, including blood flow and clotting in the cardiovascular system \cite{Fogelson2014, Leiderman2011}, particulate and chemical vapor transport in the lungs \cite{Gloede2011, Tian2010, Schroeter2006}, and drug absorption in the digestive tract \cite{Billat2017, Pang2003, Yu1996}. In some cases, critical interactions occur between fluid-phase and structure-bound chemicals. These interactions can appear in the model equations as a Robin boundary condition for the fluid-bound chemical. For example, the interaction between circulating proteins and membrane-bound proteins play a pivotal role in thrombus formation \cite{Fogelson2014}. Modeling these interactions becomes even more challenging as one considers the motion of the flow domain itself.

The numerical simulation of PDEs in complex domains has garnered significant attention for decades. Embedded boundary methods are popular approaches to such problems in which a fixed rectangular Cartesian grid is overlayed on the complex structure. Embedded boundary approaches typically alter the PDE to include an additional source term that is non-zero only near the boundary. For instance, the immersed boundary method \cite{Griffith2020,Peskin2002} uses an integral transform with a regularized delta function kernel to enforce  boundary conditions along irregular interfaces immersed in a background Cartesian grid. These methods typically are designed for Dirichlet boundary conditions, and they have to be modified with special interpolation procedures to allow for other types of boundary conditions \cite{Zhang2008, Pacheco-Vega2007, Bouchon2010}. Volume penalization \cite{Kadoch2012} methods and diffuse domain \cite{Yu2012, Li2009, Towers2018} methods both introduce phase field models to track the interface. A smoothed version of the interface is then used to modify the original PDE to account for the boundary conditions, ultimately yielding an approach similar to the immersed boundary method. These methods all have the effect of smoothing the boundary of the complex domain over several grid cells. Regularization lowers solution accuracy near the boundary, and this can limit the effectiveness to applications in which boundary interactions are important. One alternative approach is the immersed interface method \cite{Kolahdouz2020, Liu2014, Xu2006}, which derives jump conditions across the interface and then builds these jump conditions into the discretized equations. Using this method to impose boundary conditions requires relating the boundary conditions to jump conditions along the interface. Many applications have succesfully used the immersed interface method, but to our knowledge, the jump conditions for arbitrary spatially dependent Robin boundary conditions have not yet been determined. A closely related approach is the ghost fluid method \cite{Chai2020}, which uses an interpolation procedure to build stencils near the interface that incorporate the boundary conditions. This approach can lead to non-symmetric and, in some cases, ill-conditioned systems that require specialized solvers \cite{Yao2012}. Another alternative to these methods is the immersed boundary smooth extension (IBSE) method \cite{Stein2017, Qadeer2021}, which can solve PDEs on complex geometries by embedding the geometry in a simpler region and solving the PDE on the extended domain that now includes both ``physical'' and ``non-physical'' subdomains. In the IBSE approach, a body force is incorporated in the ``non-physical" regime to extend the physical solution smoothly outside the physical domain. This allows for high order accuracy to the boundary, but at the cost of solving an additional multiharmonic problem, with the order depending on the number of interface conditions.

The approach that we use here to impose boundary conditions along a complex boundary while retaining a Cartesian grid discretization framework is based on a cut-cell finite volume formulation. In these flux-based methods, fluxes are carefully calculated to account for the portion of cells that are inside the physical domain. This approach allows for accurate solutions along the boundary, but comes at the expense of computing cell geometries at the interface. This expense can be alleviated through the use of a level set function to track the surface interface \cite{Min2007}. Recent work by Helgadottir \etal ~\cite{Helgadottir2015} demonstrated the ability to impose Dirichlet, Neumann, and Robin boundary conditions for Poisson problems. Cut-cell methods can suffer from the so-called ``small cell problem," however, in which cells with small volumes necessitate the use of extremely small time step sizes for conditionally stable time stepping schemes. This becomes a serious problem for hyperbolic equations for which accurate, efficient, and unconditionally stable implicit time-stepping schemes are difficult to create. Approaches to alleviate the small cell problem for advective PDEs involve using an implicit method only for the cut-cells \cite{May2017}, merging small cells with their neighbors to effectively create a larger cell \cite{Schneiders2013}, or partitioning fluxes into ``shielded" and ``unshielded" zones based on the geometry of the cut-cell \cite{Klein2009}. These methods have been succesfully deployed in two spatial dimensions, but extending them to three spatial dimensions remains a challenge. Recently, cut-cell methods that do not suffer from the small cell problem have been introduced to handle moving boundaries \cite{Schneiders2013,Strychalski2010}; however, a formulation involving Robin boundary conditions that achieves second order accuracy has yet to be developed. The major contribution of this study is the construction of such a method with second order accuracy.

Herein, we develop a cut-cell method for advection-diffusion equations on moving domains that allows for the imposition of general Robin boundary conditions. To avoid the small cell problem, we introduce a split cut-cell semi-Lagrangian scheme, in which the diffusive operator is handled using well established finite volume methods, and the advective operator is treated using a semi-Lagrangian method. The benefit of splitting the two operators is two fold. First, in the diffusive solve, the flux from the boundary does not need to account for the change in location of the boundary. This allows us to leverage recent work by Papac \etal ~\cite{Papac2010} and Bochkov \etal ~\cite{Bochkov2019} to solve a Poisson-like problem with stationary boundaries. Second, the semi-Lagrangian scheme for advection has no CFL stability constraint, so the small cell problem is no longer an issue. We demonstrate that this method can accurately resolve concentrations near the boundary for both Robin and Neumann boundary conditions, including inhomogeneous Robin boundary conditions. In addition, we show that this method is able to handle conversion of one concentration into another across the boundary, effectively modelling surface reactions.

\section{Continuous Equations}\label{sec:CEM}

We consider the advection and diffusion of a chemical species with concentration $q\parens{\xx,t}$ in an arbitrary domain $\Omega_t$, embedded in a larger rectangular region $\mathcal{B}$. Both $q\parens{\xx,t}$ and $\Omega_t$ are transported by the same velocity field $\uu\parens{\xx,t}$. We assume that $q\parens{\xx,t}$ diffuses with diffusion coefficient $D$. We consider the particular case in which $q\parens{\xx,t}$ satisfies Robin boundary conditions on the boundary $\Gamma_t$ of $\Omega_t$, so that
\begin{subequations}
\begin{align}
\frac{\partial q\parens{\xx,t}}{\partial t} + \grad\cdot\Big(\uu\parens{\xx,t} q\parens{\xx,t} - D\grad q\parens{\xx,t}\Big) &= f\parens{\xx,t}, \xx\in\Omega_t, \label{eq:main}\\
D\grad q\parens{\xx,t}\cdot\nn\parens{\xx,t} + a\parens{\xx,t} q\parens{\xx,t} &= g\parens{\xx,t}, \xx\in \Gamma_t, \label{eq:bc}
\end{align}
\end{subequations}
in which $\nn\parens{\xx,t}$ is the outward unit normal of $\Gamma_t$ and $f\parens{\xx,t}$ is a given volumetric source function. In the implementation, $q\parens{\xx,t}$ is defined only inside $\Omega_t$. The bounding region $\mathcal{B}$ is used only to define the cells contained within $\Omega_t$.

We describe the boundary $\Gamma_t$ using a signed distance function $\phi\parens{\xx,t}$ such that 
\begin{equation}
\Gamma_t = \{\xx\in\BB |~ \phi\parens{\xx,t} = 0\}.
\end{equation}
The signed distance function is passively advected by the prescribed velocity,
\begin{equation} \label{eq:level_set_adv}
\frac{\partial \phi\parens{\xx,t}}{\partial t} + \grad\cdot\parens{\uu\parens{\xx,t} \phi\parens{\xx,t}} = 0, \xx\in\mathcal{B}.
\end{equation}
Because $\phi\parens{\xx,t}$ will not in general remain a signed distance function under advection by \eqref{eq:level_set_adv}, a reinitialization procedure is used to maintain the signed distance property. This can be achieved by computing the steady-state solution to the Hamilton-Jacobi equation
\begin{align}\label{eq:reinitialization}
\frac{\partial \hat{\phi}\parens{\xx,\tau}}{\partial \tau} + \text{sgn}\parens{\hat{\phi}\parens{\xx,\tau}}\parens{\norm*{\grad\hat{\phi}\parens{\xx,\tau}} - 1} &= 0, \\
\hat{\phi}\parens{\xx,0} &= \phi\parens{\xx,t}, \nonumber
\end{align}
after which $\phi\parens{\xx,t}$ is set to the steady state solution.

\section{Numerical Methods}\label{sec:DISCRETIZATION}

To simplify notation, we describe the method in two spatial dimensions. Extensions to a third spatial dimension are straightforward, with the exception of calculating cut-cell geometries \cite{Min2007}.

We use a splitting scheme to split equation \eqref{eq:main} into a diffusion step,
\begin{subequations} \label{eq:diff}
\begin{align}
\frac{\partial q\parens{\xx,t}}{\partial t} - \grad\cdot\parens{D\grad q\parens{\xx,t}} &= f\parens{\xx,t}, \xx\in\Omega_t \label{eq:main_diff}\\
D\grad q\parens{\xx,t}\cdot\nn\parens{\xx,t} + a\parens{\xx,t} q\parens{\xx,t} &= g\parens{\xx,t}, \xx\in \Gamma_t, \label{eq:bc_diff}
\end{align}
\end{subequations}
where the domain $\Omega_t$ remains fixed, followed by an advection step,
\begin{subequations}\label{eq:main_adv}
\begin{align}
\frac{\partial \phi\parens{\xx,t}}{\partial t} &+ \grad\cdot\parens{\uu\phi\parens{\xx,t}} = 0 \mbox{ for } \xx\in\BB \\
\frac{\partial q\parens{\xx,t}}{\partial t} &+ \grad\cdot\parens{\uu q\parens{\xx,t}} = 0 \mbox{ for } \xx\in\Omega_t,
\end{align}
\end{subequations}
in which we evolve both the level set and the function $q\parens{\xx,t}$. We note that because $\Omega_t$ and $q\parens{\xx,t}$ evolve with the same velocity, no boundary condition is needed with equation \eqref{eq:main_adv}. The use of this splitting procedure inccurs an additional error cost. This cost can be reduced to second order temporal accuracy using Strang splitting \cite{Strang1968}, which will be described later.

We overlay a Cartesian grid on top of $\BB$ such that $\BB$ consists of rectangular grid cells $\cc_{i,j}$ and $\BB = \cup \cc_{i,j}$ with a grid cell spacing of $\Delta x = \Delta y = h$. The concentration field is approximated at the cell centroid of each full or partial cell contained within $\Omega_t$.

In the following sections, we describe the discretization of equations \eqref{eq:diff} and \eqref{eq:main_adv}. In what follows, unless otherwise noted, $q_{i,j}$ refers to the concentration at the cell centroid of the cell $\cc_{i,j}\cap\Omega_t$.

\subsection{Diffusion Step}\label{sec:DIFFUSION}
\begin{figure}
\begin{center}
\phantomsubcaption\label{fig:discretization:diff}\phantomsubcaption\label{fig:discretization:boundary}
\psfragscanon
\psfrag{Gamma}{$\Gamma_t$}
\psfrag{OmN}{$\Omega_t$}
\psfrag{cij}{$\cc_{i,j}$}
\psfrag{LiNj}{$L^\text{g}_{i-\frac{1}{2},j}$}
\psfrag{LijN}{$L^\text{g}_{i,j-\frac{1}{2}}$}
\psfrag{qij1hat}{$\hat{q}_{i,j+1}$}
\psfrag{qijhat}{$\hat{q}_{i,j}$}
\psfrag{dy}{$\Delta y$}
\psfrag{LG}{$L^\Gamma_{i,j}$}
\psfrag{rij}{$\rr_{i,j}$}
\psfrag{dij}{$d_{i,j}$}
\subfigimg[width=0.35\linewidth,hsep=-1.25em,vsep=1em,pos=ul]{(a)}{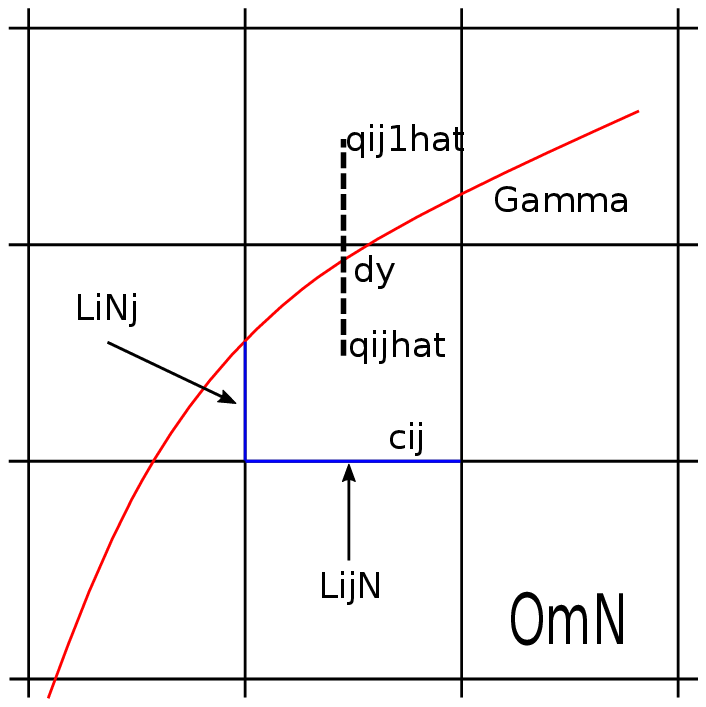}\quad\quad
\subfigimg[width=0.35\linewidth,hsep=-1.25em,vsep=1em,pos=ul]{(b)}{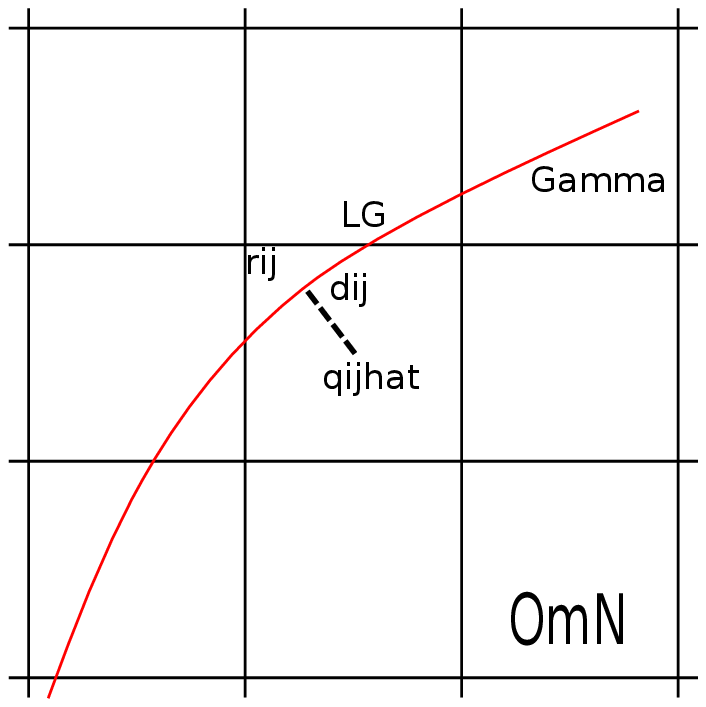}
\end{center}
\caption{Depiction of the nomenclature used in the diffusion discretization for the cell length fractions (\subref{fig:discretization:diff}) and for the boundary fluxes (\subref{fig:discretization:boundary}).}
\label{fig:discretization}
\end{figure}

To solve the diffusion step from equations \eqref{eq:diff}, we employ a cut-cell finite volume method based on the approaches of Papac \etal ~\cite{Papac2010} and Arias \etal ~\cite{Arias2018}. We summarize the derivation here, and refer the interested reader to a more detailed description in previous work. Integrating equation \eqref{eq:main_diff} over a cell $\cc_{i,j}$ that is entirely or partially interior to $\Omega_t$ and dividing by the volume of the cell, we get
\begin{equation}\label{eq:integrate_diff}
\frac{1}{\size{\cc_{i,j}\cap\Omega_t}}\int_{\cc_{i,j}\cap\Omega_t} \frac{\partial q\parens{\xx,t}}{\partial t} \mathrm{d}\xx = \frac{1}{\size{\cc_{i,j}\cap\Omega_t}}\int_{\cc_{i,j}\cap\Omega_t} D\Lap q\parens{\xx,t} \mathrm{d}\xx.
\end{equation}
We define $Q_{i,j}$ as the cell average of $q\parens{\xx,t}$ in the cell $\cc_{i,j} \cap\Omega_t$. Replacing the cell average in the left side of equation \eqref{eq:integrate_diff} and employing the divergence theorem on the right-hand side, we get
\begin{equation}\label{eq:diff_integral}
\frac{d Q_{i,j}}{d t} = \frac{1}{\size{\cc_{i,j}\cap\Omega_t}}\int_{\partial\parens{\cc_{i,j}\cap\Omega_t}} D\grad q\parens{\xx,t}\nn\cdot d\AA,
\end{equation}
in which $\nn$ is the outward unit normal of $\partial\parens{\cc_{i,j}\cap\Omega_t}$. We can further divide the integral in equation \eqref{eq:diff_integral} into an integral over the cell boundary $\partial \cc_{i,j}\cap\Omega_t$ and an integral over the physical boundary $\Gamma_t\cap\cc_{i,j}$
\begin{equation}
\int_{\partial\parens{\cc_{i,j}\cap\Omega_t}} D\grad q\parens{\xx,t}\cdot\nn dA = \parens{\int_{\partial \cc_{i,j}\cap\Omega_t} + \int_{\cc_{i,j}\cap\Gamma_t}} D\grad q\parens{\xx,t}\nn\cdot d\AA.
\end{equation}
We approximate the first integral by
\begin{align}\label{eq:first_integral}
\int_{\partial\cc_{i,j}\cap\Omega_t} D\grad q\parens{\xx,t}\cdot\nn dA \approx &L^\text{g}_{i+\frac{1}{2},j}\frac{\hat{q}_{i+1,j}-\hat{q}_{i,j}}{\Delta x} - L^\text{g}_{i-\frac{1}{2},j}\frac{\hat{q}_{i,j} - \hat{q}_{i-1,j}}{\Delta x} \\
 + &L^\text{g}_{i,j+\frac{1}{2}}\frac{\hat{q}_{i,j+1}-\hat{q}_{i,j}}{\Delta y} - L^\text{g}_{i,j-\frac{1}{2}}\frac{\hat{q}_{i,j}-\hat{q}_{i,j-1}}{\Delta y}, \nonumber
\end{align}
in which $\hat{q}_{i,j}$ is the point-wise concentration at the center of the cell $\cc_{i,j}$ and $L^\text{g}_{i+\frac{1}{2},j}$ is the length fraction of the face $\parens{i+\frac{1}{2}}\times \left[j-\frac{1}{2},j+\frac{1}{2}\right]$ covered by the irregular domain, see Figure \ref{fig:discretization:diff}. It is challenging to compute $L^\text{g}_{i+\frac{1}{2},j}$ exactly. The smoothness of the level set is determined by the zero contour. For shapes that have sharp features, the level set will inherit these features. However, to achieve second order accuracy for smooth level sets, it suffices to use a linear approximation
\begin{equation}\label{eq:L}
L^\text{g}_{i+\frac{1}{2},j} = \left\{ \begin{array}{cc}
\Delta y \size{\frac{\phi_{i+\frac{1}{2},j-\frac{1}{2}}}{\phi_{i+\frac{1}{2},j-\frac{1}{2}}-\phi_{i+\frac{1}{2},j+\frac{1}{2}}}} &\mbox{ if } \phi_{i+\frac{1}{2},j-\frac{1}{2}}\cdot\phi_{i+\frac{1}{2},j+\frac{1}{2}} < 0 \\
\Delta y &\mbox{ if } \phi_{i+\frac{1}{2},j-\frac{1}{2}} < 0 \mbox{ and } \phi_{i+\frac{1}{2},j+\frac{1}{2}} < 0 \\
0 &\mbox{ otherwise.}

\end{array}\right.
\end{equation}
While in the evolution of the level set $\phi\parens{\xx,t}$, the degrees of freedom live at cell centers, in equation \eqref{eq:L}, we require the value of $\phi\parens{\xx,t}$ at cell nodes. In our computations, we use a simple bi-linear interpolant to find these nodal values. In three spatial dimensions, evaluating the corresponding quantity $L^\text{g}_{i+\frac{1}{2},j,k}$ would involve computing the surface area.

We emphasize that $\hat{q}_{i,j}$ refers to the cell center of the (possibly cut) grid cell $\cc_{i,j}$ regardless of the location of the physical boundary $\Gamma_t$. In cases where the cell center does not correspond to the location of degrees of freedom, e.g.\ near cut-cells, we must reconstruct these values. Here, we perform this reconstruction using either a moving least squares (MLS) approximation, or a radial basis function (RBF) interpolant. Both procedures are described in more detail in Section \ref{sec:MLS_RBF}. We note that this reconstruction is not required in the methods of Arias \etal ~\cite{Arias2018} and Bochkov \etal ~\cite{Bochkov2019}, as those methods always define the degrees of freedom at cell centers.

As done in Bochkov \etal ~\cite{Bochkov2019}, the second integral is approximated using a linear approximation to $q\parens{\xx,t}$ on the boundary in the direction normal to $\Gamma_t$
\begin{equation}\label{eq:second_integral}
\int_{\cc_{i,j}\cap\Gamma_t} D\grad q\parens{\xx,t}\cdot\nn dA = \int_{\cc_{i,j}\cap\Gamma_t} \parens{g\parens{\xx,t} - a q\parens{\xx,t}}dA \approx L^\Gamma_{i,j}\parens{g\parens{\rr_{i,j}, t} - a q\parens{\rr_{i,j}, t}},
\end{equation}
in which $\rr_{i,j}$ is the closest point on $\Gamma_t$ to $\xx_{i,j}$ and $L^\Gamma_{i,j} = \size{\cc_{i,j}\cap\Gamma_t}$ (see Figure \ref{fig:discretization:boundary}). The value $q\parens{\rr_{i,j}, t}$ is found using a Taylor series expansion 
\begin{equation}\label{eq:bdry_1}
q\parens{\xx_{i,j}, t} = q\parens{\rr_{i,j}, t} + d_{i,j}\frac{\partial q\parens{\rr_{i,j}, t}}{\partial \nn} + \bigOh{h^2},
\end{equation}
in which $d_{i,j} = \frac{\phi\parens{\xx_{i,j}, t}}{\size{\grad\phi\parens{\xx_{i,j}, t}}}$ is the distance between $\rr_{i,j}$ and cell center $\xx_{i,j}$. On the boundary, we have that
\begin{equation}\label{eq:bdry_2}
D\frac{\partial q\parens{\rr_{i,j},t}}{\partial \nn} + a q\parens{\rr_{i,j}, t} = g\parens{\rr_{i,j},t}.
\end{equation}
We solve equations \eqref{eq:bdry_1} and \eqref{eq:bdry_2}, dropping the $\bigOh{h^2}$ terms, for $q\parens{\rr_{i,j},t}$ and use this value in the approximation to the integral in equation \eqref{eq:second_integral}. As shown previously \cite{Bochkov2019}, this system is well posed for a sufficiently refined grid.

The cut-cell volume $\size{\cc_{i,j}\cap\Omega_t}$ and physical boundary length $\size{\cc_{i,j}\cap\Gamma_t}$ are found by decomposing the cut-cell and boundary into simplices, for which analytic formulas for the volume exist \cite{Min2007}.

We approximate equation \eqref{eq:diff} using the implicit trapezoidal rule:

\begin{align}\label{eq:diff_timestepping}
\frac{1}{\Delta t} Q_{i,j}^{n+1} &- \frac{D}{2\size{\cc_{i,j}\cap\Omega^-}}\left(L^\text{g}_{i+\frac{1}{2},j}\frac{q^{n+1}_{i+1,j}-q^{n+1}_{i,j}}{\Delta x} - L^\text{g}_{i-\frac{1}{2},j}\frac{q^{n+1}_{i,j} - q^{n+1}_{i-1,j}}{\Delta x}\right.& \nonumber \\
&\left.+ L^\text{g}_{i,j+\frac{1}{2}}\frac{q^{n+1}_{i,j+1}-q^{n+1}_{i,j}}{\Delta y} - L^\text{g}_{i,j-\frac{1}{2}}\frac{q^{n+1}_{i,j}-q^{n+1}_{i,j-1}}{\Delta y}\right)  \nonumber \\
&= \frac{1}{\Delta t} Q_{i,j}^n+ \frac{D}{2\size{\cc_{i,j}\cap\Omega^-}}\left(L^\text{g}_{i+\frac{1}{2},j}\frac{q^n_{i+1,j}-q^n_{i,j}}{\Delta x} - L^\text{g}_{i-\frac{1}{2},j}\frac{q^n_{i,j} - q^n_{i-1,j}}{\Delta x}\right. \\
&\left.+ L^\text{g}_{i,j+\frac{1}{2}}\frac{q^n_{i,j+1}-q^n_{i,j}}{\Delta y} - L^\text{g}_{i,j-\frac{1}{2}}\frac{q^n_{i,j}-q^n_{i,j-1}}{\Delta y}\right)+ \frac{1}{\size{\cc_{i,j}\cap\Omega^-}}f_{i,j}^{n+\frac{1}{2}}. \nonumber
\end{align}
In our computational examples, we use GMRES without preconditioning to solve this system of equations. The solver typically converges to a solution with tolerance $10^{-12}$ after approximatley 20 to 50 iterations.

\subsection{Advection Step}\label{sec:ADVECTION}
For the advection step, we use a semi-Lagrangian method \cite{Durran1999,Russo2009} to advance both the level set $\phi\parens{\xx,t}$ and the concentration field $q\parens{\xx,t}$. We solve equation \eqref{eq:main_adv} by first finding the preimage $\XX_{i,j} = \bchi^{-1}\parens{\xx_{i,j},t^{n+1}}$ of the fluid parcel located at $\xx_{i,j}$ at time $t^{n+1}$. The mapping $\bchi\parens{\xx,t}$ satisfies the differential equation
\begin{equation}\label{eq:material_eq}
\frac{\partial\bchi\parens{\xx,t}}{\partial t} = \uu\parens{\xx,t} \mbox{ for } \xx\in\BB.
\end{equation}
The preimage $\XX_{i,j}$ can be found by integrating equation \eqref{eq:material_eq} backwards in time. We use an explicit two step Runge Kutta method:
\begin{subequations}\label{eq:material_approx}
\begin{align}
\XX^\star_{i,j} = \xx_{i,j} - \frac{\Delta t}{2}\uu\parens{\xx_{i,j}, t^{n+1}}, \\
\XX_{i,j}^n = \xx_{i,j} - \Delta t \uu\parens{\XX^\star_{i,j},t^{n+\frac{1}{2}}}.
\end{align}
\end{subequations}
Having found the preimage, we can interpolate the solution at time $t^n$ at the preimage locations $\XX_{i,j}^n$.

Our interpolation procedure consists of one of two possible methods depending on whether $\XX_{i,j}^n$ is near cut-cells. Away from cells pierced by the zero level set, we use a tensor product of special piecewise Hermite polynomials called Z-splines \cite{Becerra-Sagredo2016}. The Z-spline $Z_m\parens{x}$ interpolating the data $\parens{x_i,f_i = f\parens{x_i}}_{i = 1}^n$ is a piecewise polynomial function that satisfies
\begin{subequations}
\begin{align}
Z_m\parens{x}&\in C^m\parens{\left[x_1,x_n\right]}, \\
\left.\frac{d^p}{d x^p} Z_m\parens{x}\right\rvert_{x_j} &= f_{m,j}^p \mbox{ for } p=0,\ldots,m \mbox { and } j=1,\ldots,n, \\
Z_m\parens{x}&\in \pi_{2m+1}\parens{\left[x_i,x_{i+1}\right]} \mbox{ for } i=1,\ldots,n-1,
\end{align}
\end{subequations}
in which $f_{m,j}^p$ is the approximation of the $p^\text{th}$ order derivative of $f(x)$ computed from high-order finite differences of $f_j$ using $2m + 1$ points, and $\pi_{2m+1}$ is the space of polynomials of degree less than or equal to $2m+1$.

It is possible to define $Z_m$ for a general set of data points (or abscissae) via the cardinal Z-splines \cite{Becerra-Sagredo2016}
\begin{equation}\label{eq:interpolant}
Z_m\parens{x} = \sum_i f_i\tilde{Z}_m\parens{x - x_i}.
\end{equation}
The cardinal Z-splines have two key properties. First, the cardianl Z-splines have compact support, $\tilde{Z}_m\parens{x} = 0$ for $\size{x} > m+1$, which allows for efficient local evaluation of interpolants. Second, cardinal Z-splines are interpolatory, which makes the interpolant trivial to form. Because data values are defined on a regular Cartesian grid, the computation of the full interpolant uses a tensor product of cardinal Z-splines
\begin{equation}\label{eq:interpolant_full}
Z_m\parens{\xx} = \sum_i\sum_j q_{i,j}\tilde{Z}_m\parens{x - x_i}\tilde{Z}_m\parens{y - y_j}.
\end{equation}
In our computations, we use the quintic Z-splines, which are defined in terms of
\begin{equation}
Z_2\parens{x} = \left\{ \begin{array}{cc}
1 - \frac{15}{12} x^2 - \frac{35}{12} x^3 + \frac{63}{12} x^4 - \frac{25}{12} x^5 &\mbox{ if } \size{x} \leq 1 \\
-4 + \frac{75}{4} x - \frac{245}{8} x^2 + \frac{545}{24} x^3 - \frac{63}{8} x^4 + \frac{25}{24} x^5 &\mbox{ if } 1 < \size{x} \leq 2 \\
18 - \frac{153}{4} x + \frac{255}{8} x^2 - \frac{313}{24} x^3 + \frac{21}{8} x^4 - \frac{5}{24} x^5 &\mbox{ if } 2 < \size{x} \leq 3 \\
0 &\mbox{ otherwise.}
\end{array}\right.
\end{equation}

In cases where we do not have enough points to define the Z-spline, e.g., near cut-cells, we use the procedure described in Section \ref{sec:MLS_RBF}. We note that although the interpolation procedure described here has a discontinuous switch between operators, this does not appear to affect the overall convergence rates in our numerical tests of the methodology.

We note that this form of the semi-Lagrangian method is not conservative, but, in our experiments, the change in the amount of material was very small. Conservative verions of semi-Lagrangian methods have been developed \cite{Lentine2011, Iske2004}. The core change is to advect \emph{grid cells} instead of individual points, and then integrate the resulting interpolating polynomial over the grid cell. Because conservation is not critical to our ultimate problem of interest, we use the simpler non-conservative approach in this work.

The evolution of the level set and the concentration use the same procedure, except for the choice of reference grid. The reference grid for the level set consists of the entire domain $\BB$, whereas the reference grid for the concentration consists of the domain $\Omega_{t^{n+1}} = \left\{ \xx\in\BB | \phi\parens{\xx,t^{n+1}} < 0\right\}$. Therefore, we update the level set prior to updating the concentration.

\subsection{Moving Least Squares and Radial Basis Function Interpolation}\label{sec:MLS_RBF}
Near cut-cells, it is necessary to form interpolants on unstructured data. In this study, we compare the accuracy of a moving least squares (MLS) approximation and a local radial basis function (RBF) interpolant. This procedure is used both in the diffusion step, to extrapolate data from cut-cell centroids to full-cell centers, and also in the advection step, for cells where the Z-spline interpolant can not be formed. We describe the procedure in the context of reconstructing a function $f\parens{\xx}$ from the data points $\{\xx_i,f\parens{\xx_i}\}$ with $i = 1,2,\ldots,n$ in which $n$ is an arbitrary number of points.

The MLS approximation is formed by finding the best approximation $q\parens{\xx}$ at a point $\xx_\text{c}$ of the form
\begin{equation}
q\parens{\xx} = \sum_{j = 1}^{N} c_j p_j\parens{\xx}
\end{equation}
in which the $p_j\parens{\xx}$ form a basis for the space of polynomials up to a certain degree and $N$ is the number of polynomials in the basis. The approximation $q\parens{\xx}$ is chosen to be optimal with respect to the standard weighted $L^2$ inner product with weight function $w\parens{\xx}$
\begin{equation}
\langle f\parens{\xx},g\parens{\xx}\rangle = \int_\Omega f\parens{\xx}g\parens{\xx}w\parens{\xx}d\xx.
\end{equation}
We use a stencil width of approximately two grid cells. It has been shown that the weight function must be singular at the data locations $\xx_i$ for the MLS approximation to interpolate the data \cite{Bayona2019}. Here, we use
\begin{equation}
w\parens{\xx} = e^{-\norm*{\xx - \xx_\text{c}}^2}.
\end{equation}
Because this weighting function is not singular at $\xx_\text{c}$, the reconstructed polynomial will not be interpolatory. The MLS calculation reduces to solving a linear system of the form
\begin{equation}
\Lambda A \cc = \Lambda \ff
\end{equation}
in which $A$ is a matrix whose entries consist of $a_{i,j} = p_j\parens{\xx_i}$ and $\Lambda$ is a diagonal matrix with $\lambda_{i,i} = w\parens{\xx_i}$. In preliminary computational tests, we observed that using a lower order reconstruction for the diffusion step is more robust without affecting the order of accuracy while quadratic reconstructions can be used in the advection step. Consequently, we use a quadratic polynomial in the reconstructions for the advection step and a linear reconstruction with the diffusion step for the tests reported in section \ref{sec:RESULTS}.

The RBF approximation is constructed via a polyharmonic spline of the form
\begin{equation}
q\parens{\xx} = \sum_{j = 1}^n \lambda_j\phi\parens{\norm*{\xx-\xx_j}} + \sum_{j = 1}^s \beta_j p_j\parens{\xx}
\end{equation}
in which $\phi\parens{\xx}$ is a polyharmonic radial basis function of degree $m$ and $p_j\parens{\xx}$ are a set of $s$ polynomial basis functions \cite{Bayona2019, Flyer2016}. The degree of the polynomials $k$ is chosen such that $m = 2k + 1$. The coefficients are chosen so that
\begin{align}
q\parens{\xx_j} &= f_j \mbox{ for } j = 1,\ldots,n,\\
\sum_{k=1}^n\lambda_j p_k\parens{\xx_j} &= 0  \mbox{ for } j = 1,\ldots,n.
\end{align}
This leads to the linear system for the coefficients
\begin{equation}
\left(\begin{array}{cc}
A & P \\
P^T & 0 \\
\end{array}\right)
\left(\begin{array}{c}
\lambda \\
\beta
\end{array}\right) = 
\left(\begin{array}{c}
\ff \\
0
\end{array}\right)
\end{equation}
in which $A$ is a square matrix with elements $A_{i,j} = \phi\parens{\norm*{\xx_i - \xx_j}}$, $P$ is a rectangular matrix with elements $P_{i,j} = p_j\parens{\xx_i}$. In our computations, we use the polyharmonic spline $\phi\parens{\xx} = \norm*{\xx}^3$ combined with linear polynomials.

\subsection{Full procedure and Implementation}
We summarize the full procedure to advance the solution from time $t^{n}$ to time $t^{n+1}$:
\begin{enumerate}[1)]
\item If needed, reinitialize the level set $\phi\parens{\xx,t}$ into the signed distance function by iterating equation \eqref{eq:reinitialization} to steady state using the procedure as described by Nangia \etal ~\cite{Nangia2019}.
\item Update the diffusion equation to half time $t^{n+\frac{1}{2}}$ using the methods described in Section \ref{sec:DIFFUSION}.
\item Update the level set using the prescribed velocity $\uu$ to time $t^{n+1}$ by the methods described in Section \ref{sec:ADVECTION}.
\item Advect the concentration using the prescribed velocity $\uu$ by the methods described in Section \ref{sec:ADVECTION}.
\item Update the diffusion equation to full time $t^{n+1}$ using the concentration and level sets from the previous two steps.
\end{enumerate}

The above procedure is implemented using the SAMRAI \cite{Hornung2002} infrastructure, which provides an efficient, parallelized environment for structured adaptive mesh refinement. The diffusion solve is computed using matrix free solvers with operators provided by IBAMR \cite{Griffith2007, IBAMR} and Krylov methods provided by PETSc \cite{petsc-efficient, petsc-user-ref, petsc-web-page}.

\section{Results}\label{sec:RESULTS}

Here we demonstrate the capabilities of the method both using a prescribed level set and an level set advected with the fluid. We start with diffusion dominated examples before exploring inclusion of advection and spatially varying boundary conditions.
\subsection{Diffusion from a point source}

We consider the advection-diffusion of a point source within a disk of radius $R$. The disk is advected with velocity $\uu = \parens{\cos\parens{\frac{\pi}{4}},\sin\parens{\frac{\pi}{4}}}$. The exact concentration is 
\begin{equation}\label{eq:exact}
q\parens{\xx,t} = \frac{10}{4 D \parens{t+\frac{1}{2}}}e^{-\frac{\norm*{\xx - \xx_\text{c}\parens{t}}^2}{4D\parens{t+\frac{1}{2}}}} \mbox{ for } \xx \in \BB,
\end{equation}
in which $\xx_\text{c}\parens{t}$ is the center of the disk. We apply Robin boundary conditions of the form
\begin{equation}\label{eq:bdry_conds}
D\frac{\partial q}{\partial\nn} = g\parens{\xx,t} - a\parens{\xx,t} q\parens{\xx,t}.
\end{equation}
We set $a\parens{\xx,t} = 1$ and use the method of manufactured solutions to determine the value of $g\parens{\xx,t}$. For this example, we specify the location of the level set $\phi$ at each timestep. This procedure allows us test the convergence of the advection-diffusion step without having to account for numerical error in the evolution of the level set function. We note that whereas this example includes a trivial semi-Lagrangian backwards integration step, the interpolation procedure and diffusion step are non-trivial.

The extended domain $\BB$ is the box $[0,12]\times[0,12]$ and is discretized using $N$ points in each direction. We use a disk radius of $R = 1$ and initial center $\xx_\text{c}(0) = \parens{1.521,1.503}$. The choice of center slightly offsets the disk from the background grid so that different cut-cells are generated on opposite sides of the disk. We use a diffusion coefficient of $D = 0.01$. The simulations are run at an advective CFL number of $C_{\text{CFL}} = \frac{\Delta x}{\size{\uu}\Delta t} = 0.5$.

\begin{figure}
\begin{center}
\phantomsubcaption\label{fig:adv_diff_pt_source:rbf}\phantomsubcaption\label{fig:adv_diff_pt_source:mls}
\phantomsubcaption\label{fig:adv_diff_pt_source:bdry}
\subfigimg[width=0.45\linewidth,hsep=-1em,vsep=2em,pos=ul]{(a)}{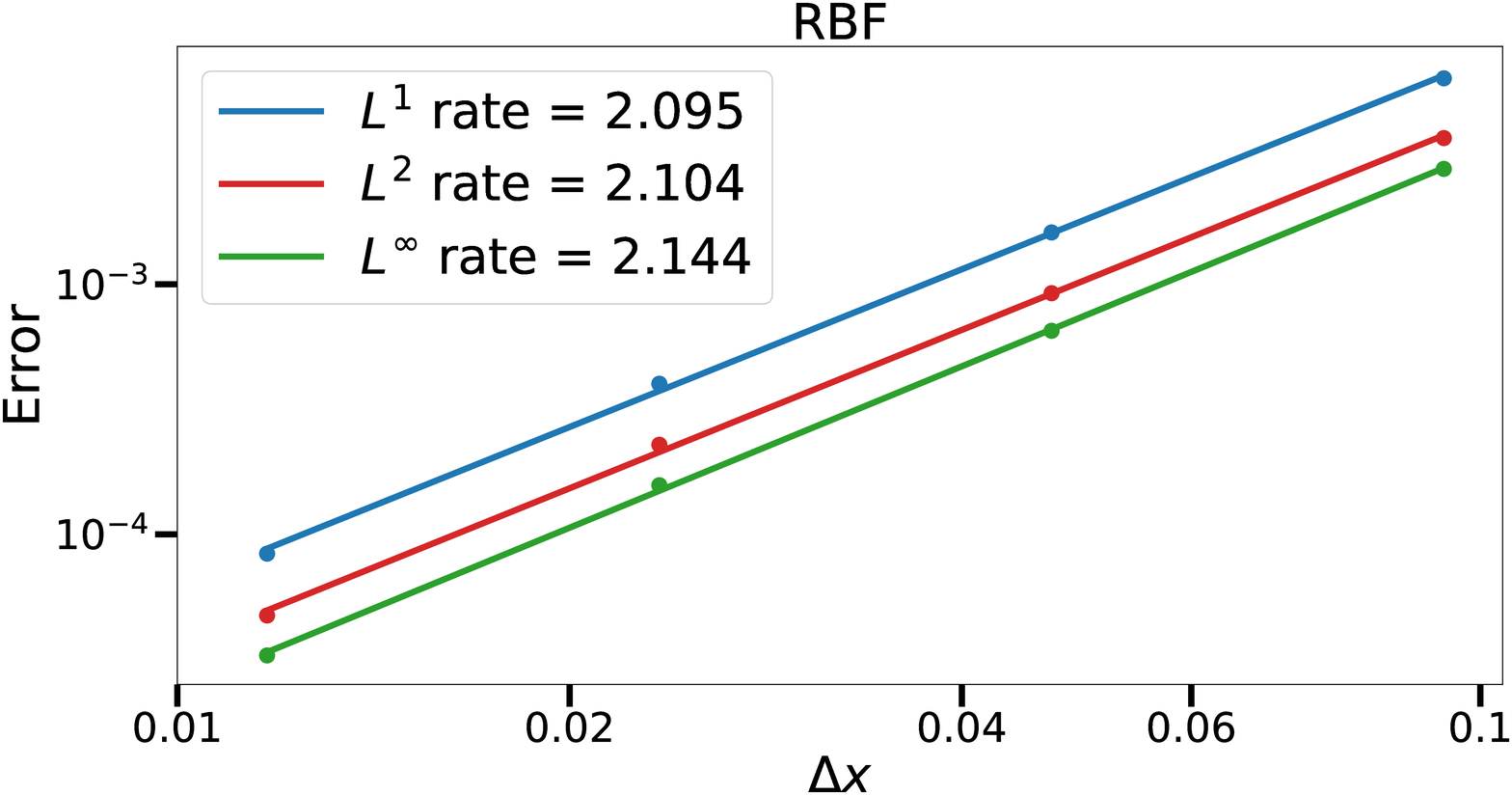}
\subfigimg[width=0.45\linewidth,hsep=-1em,vsep=2em,pos=ul]{(b)}{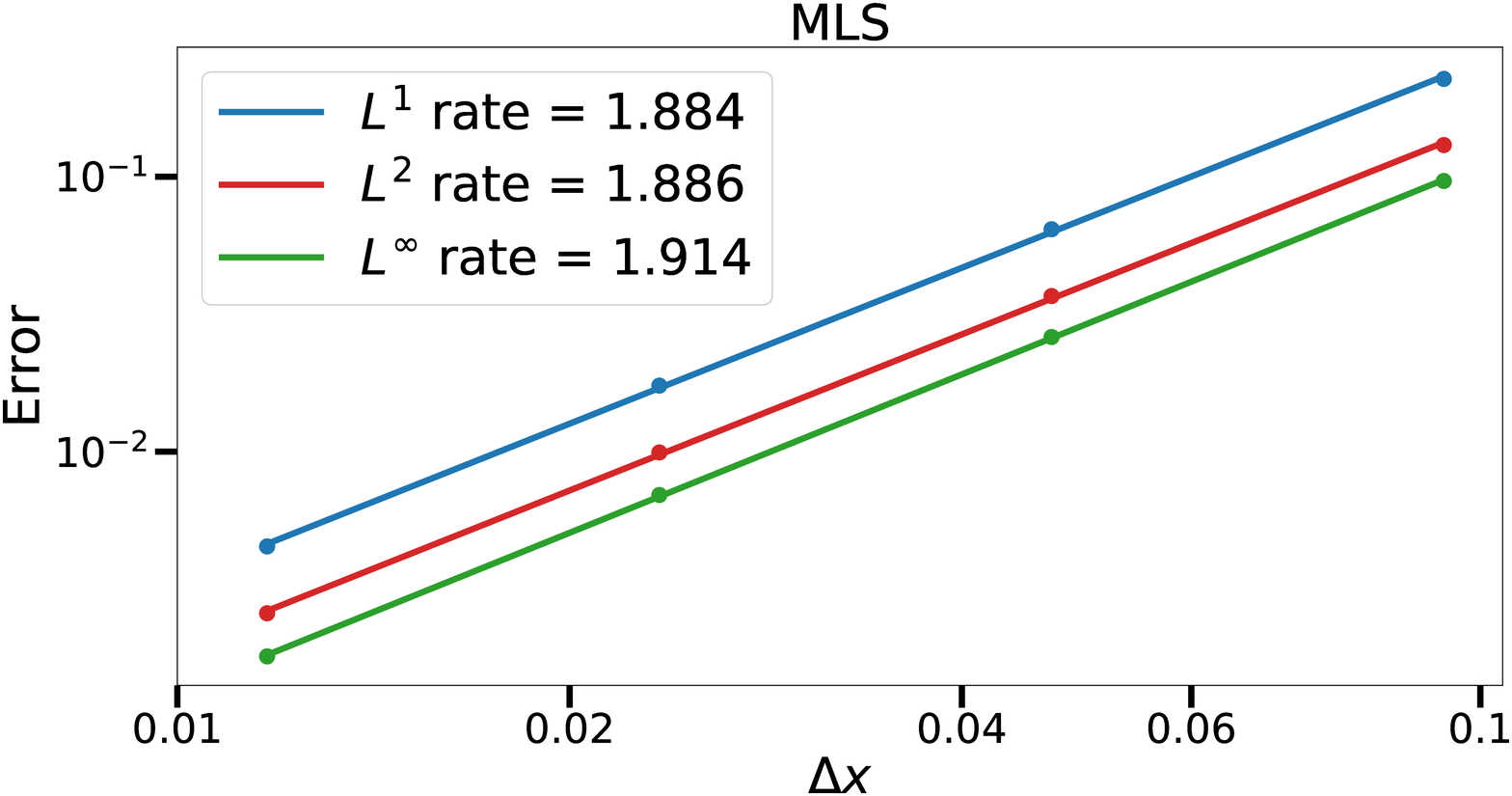}\\
\subfigimg[width=0.45\linewidth,hsep=-1em,vsep=2em,pos=ul]{(c)}{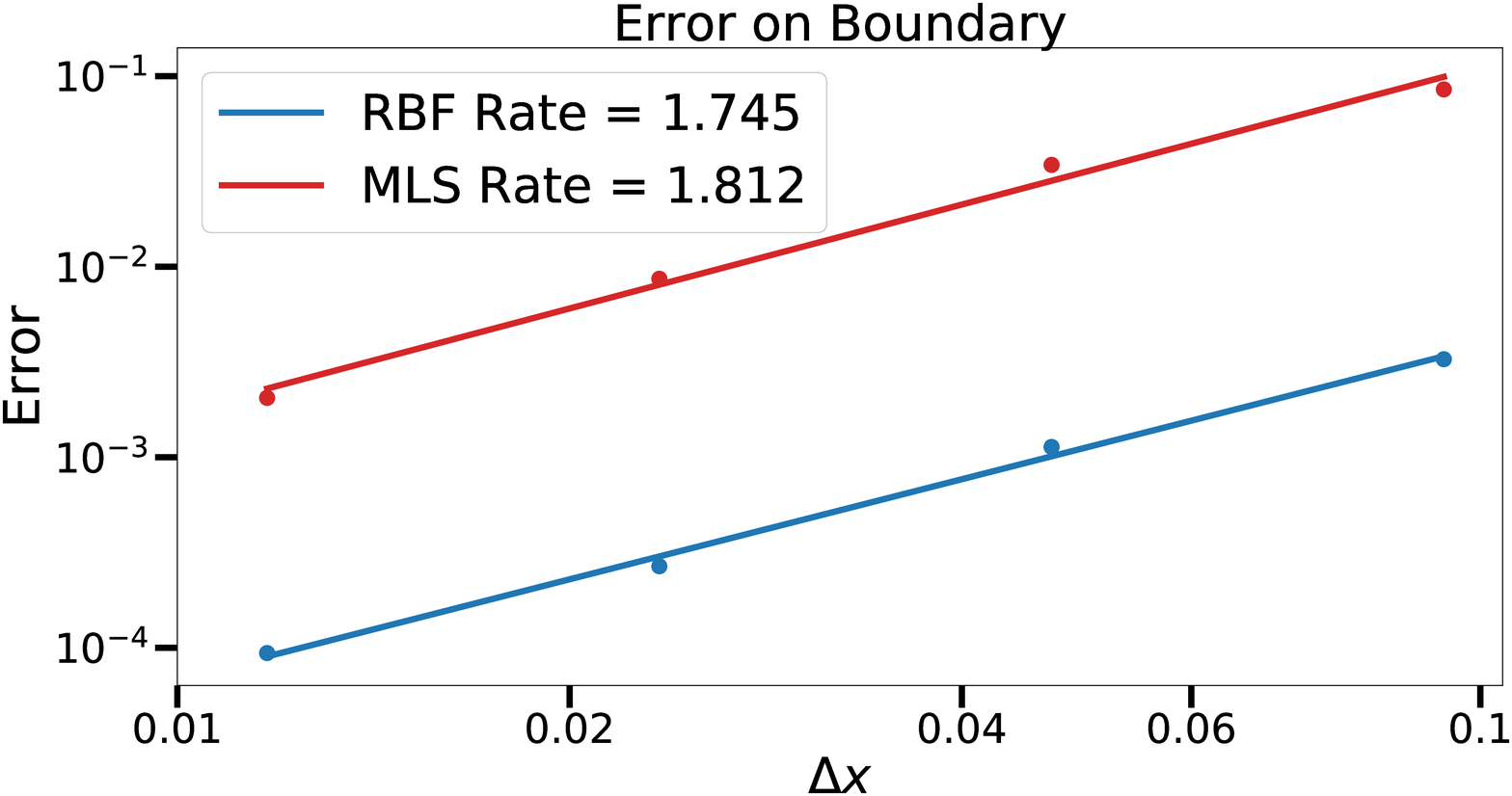}
\end{center}
\caption{Convergence rates for an advecting and diffusion concentration field initialized from a point source in a constant flow using either a RBF (\subref{fig:adv_diff_pt_source:rbf}) or MLS (\subref{fig:adv_diff_pt_source:mls}) reconstruction. Also shown is the convergence rate for the integral in equation \eqref{eq:Error_moving} using either an RBF or MLS reconstruction (\subref{fig:adv_diff_pt_source:bdry}). In all cases, the error using RBF reconstructions are more than an order of magnitude smaller than MLS reconstructions. Simulations are run to a final time of $T = 10$ and at a CFL number of $C_\text{CFL} = 0.5$. We use a diffusion coefficient of $D = 0.01$.}
\label{fig:adv_diff_pt_source}
\end{figure}

Figure \ref{fig:adv_diff_pt_source} shows convergence plots, which demonstrate second order accuracy in the $L^1$, $L^2$, and $L^\infty$ norms. The coarsest simulation uses $N = 128$ grid points in each direction of the background Cartesian mesh, which corresponds to approximately 26 grid cells covering the diameter of the disk. Despite the relatively coarse description of the disk, we still see errors that are on the order of one percent.

We also consider the case where the level set $\phi$ is advected with the same velocity. In this case, numerical issues make it difficult to assess the pointwise error on the Cartesian grid because the piecewise linear reconstruction used to determine the cut-cells gives different cut-cell geometries than that of a prescribed level set. Instead, we compare the solution at the boundary to the imposed boundary condition. We compute the solution on the boundary using a moving least squares linear extrapolation from cut-cells to the cell boundary. Further, to assess the convergence rate, we compute the error of the surface integral
\begin{equation}\label{eq:Error_moving}
E = \size{\int_0^{2\pi} q\parens{\xx,t} d\theta - \int_0^{2\pi} \tilde{q}\parens{\xx,t} d\theta},
\end{equation}
in which $\theta$ is the angle along the boundary from the center of the disk, $q\parens{\xx,t}$ is the exact solution from \eqref{eq:exact}, and $\tilde{q}\parens{\xx,t}$ is the approximate solution computed using the moving least squares extrapolation as described in Section \ref{sec:MLS_RBF}. A convergence plot is shown in Figure \ref{fig:adv_diff_pt_source} and indicates that the method achieves between first and second order convergence rates. In both the evolved and prescribed level set cases, the RBF reconstruction yields errors that are more than an order of magnitude smaller than the MLS reconstructions.

\subsection{Solid body rotation}
In this section, we again consider diffusion from a point source, but instead apply a solid body rotation to the disk. In this case, the semi-Lagrangian sub-step contains a non-trivial integration in time, and therefore tests all components of the solver. We again prescribe the initial condition as in equation \eqref{eq:exact}, and apply boundary conditions as in equation \eqref{eq:bdry_conds}.

\begin{figure}
\begin{center}
\phantomsubcaption\label{fig:rotational_init_final:initial}\phantomsubcaption\label{fig:rotational_init_final:final}
\subfigimg[width=0.45\linewidth,hsep=-1em,vsep=2em,pos=ul]{(a)}{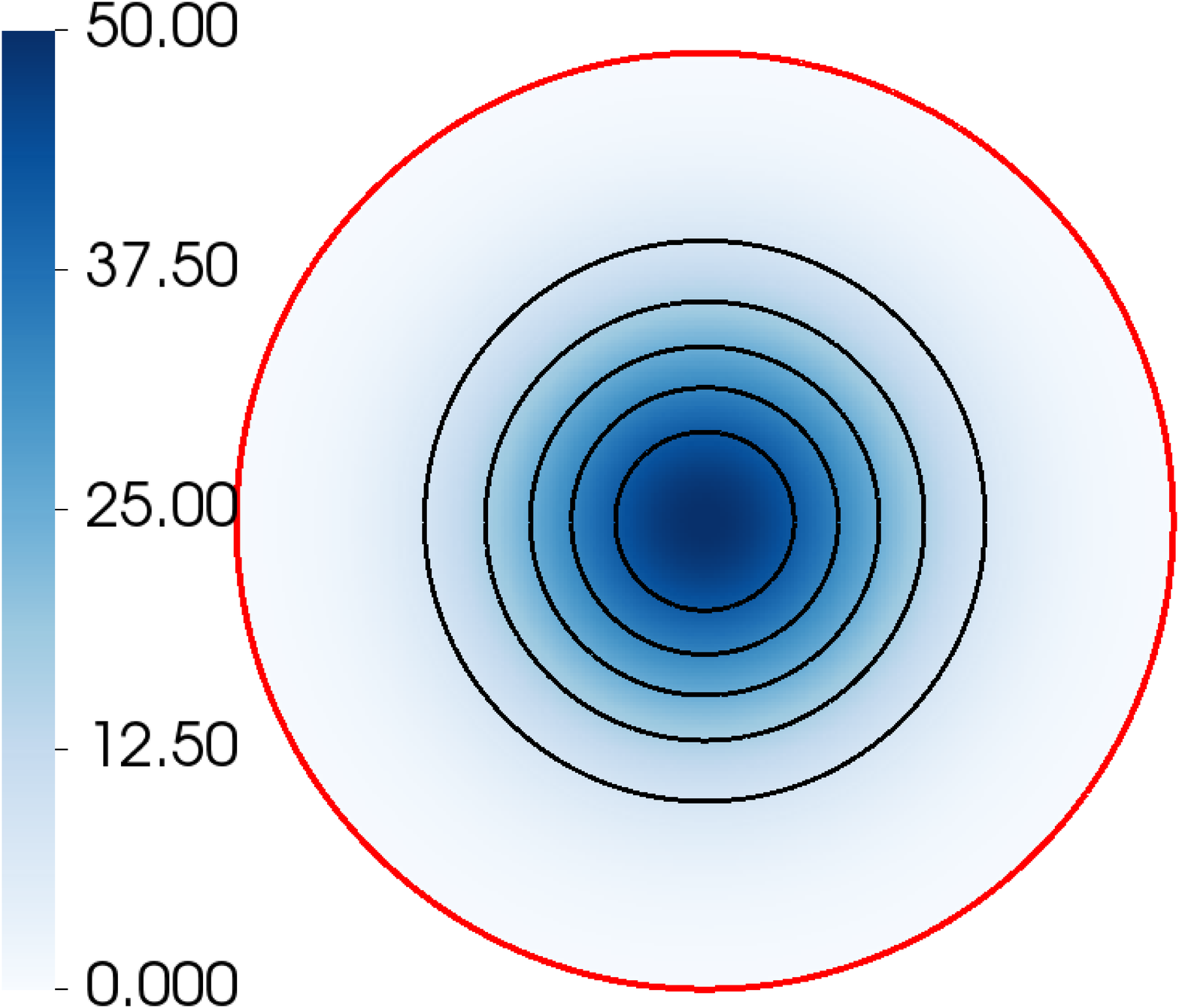}
\subfigimg[width=0.45\linewidth,hsep=-1em,vsep=2em,pos=ul]{(b)}{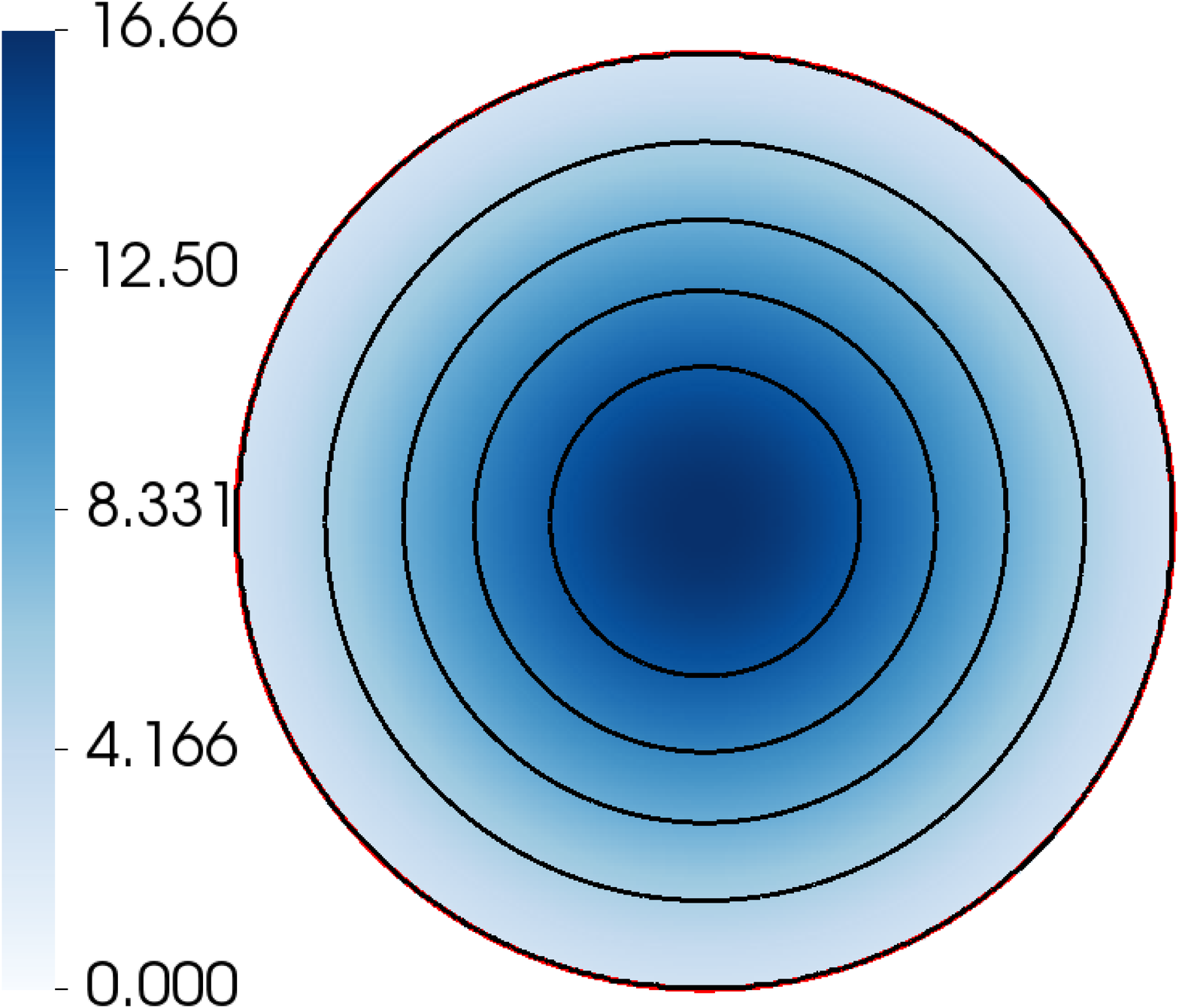}
\end{center}
\caption{Initial (\subref{fig:rotational_init_final:initial}) and final (\subref{fig:rotational_init_final:final}) configuration for a diffusing concentration field initialized from a point source in a rigid body rotational flow field. The zero contour of $\phi$ is shown in red. Simulations are run at an advective CFL number of $C_\text{CFL} = 0.5$ to a final time of $T = 1$.}
\label{fig:rotational_init_final}
\end{figure}

\begin{figure}
\begin{center}
\phantomsubcaption\label{fig:adv_diff_rotational:rbf}\phantomsubcaption\label{fig:adv_diff_rotational:mls}
\phantomsubcaption\label{fig:adv_diff_rotational:bdry}
\subfigimg[width=0.45\linewidth,hsep=-0.5em,vsep=2em,pos=ul]{(a)}{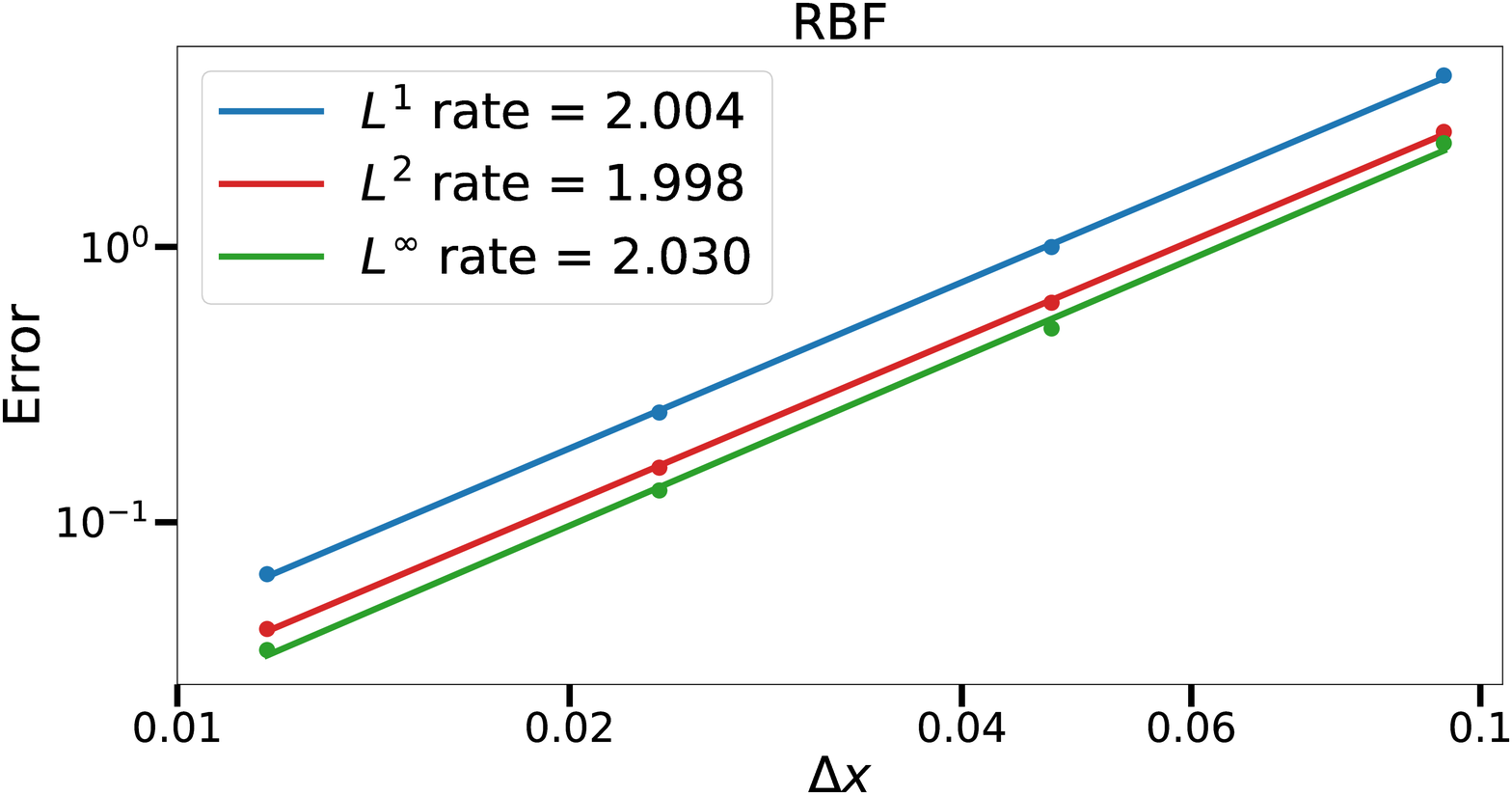}
\subfigimg[width=0.45\linewidth,hsep=-0.5em,vsep=2em,pos=ul]{(b)}{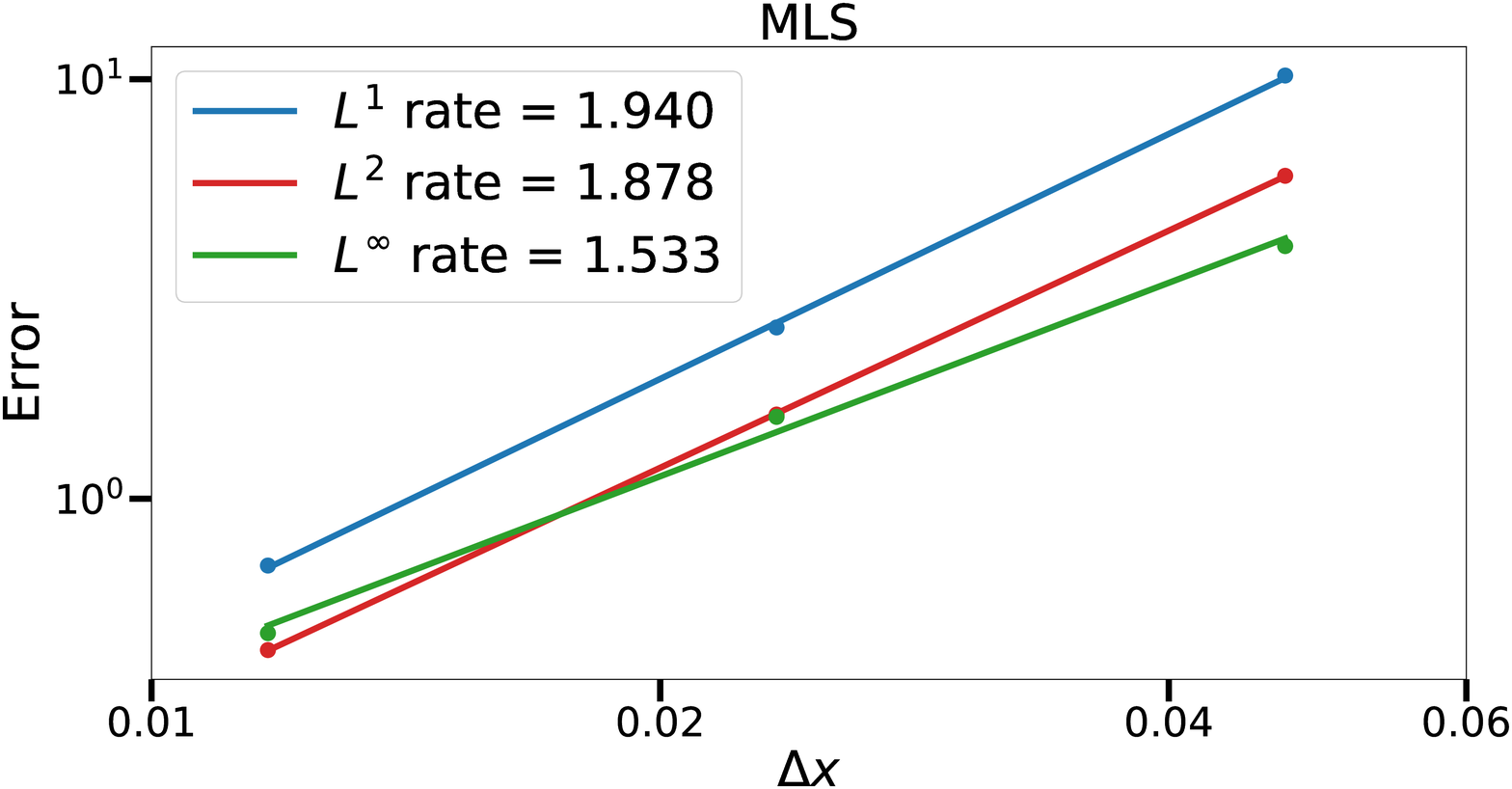}\\
\subfigimg[width=0.45\linewidth,hsep=-0.5em,vsep=2em,pos=ul]{(c)}{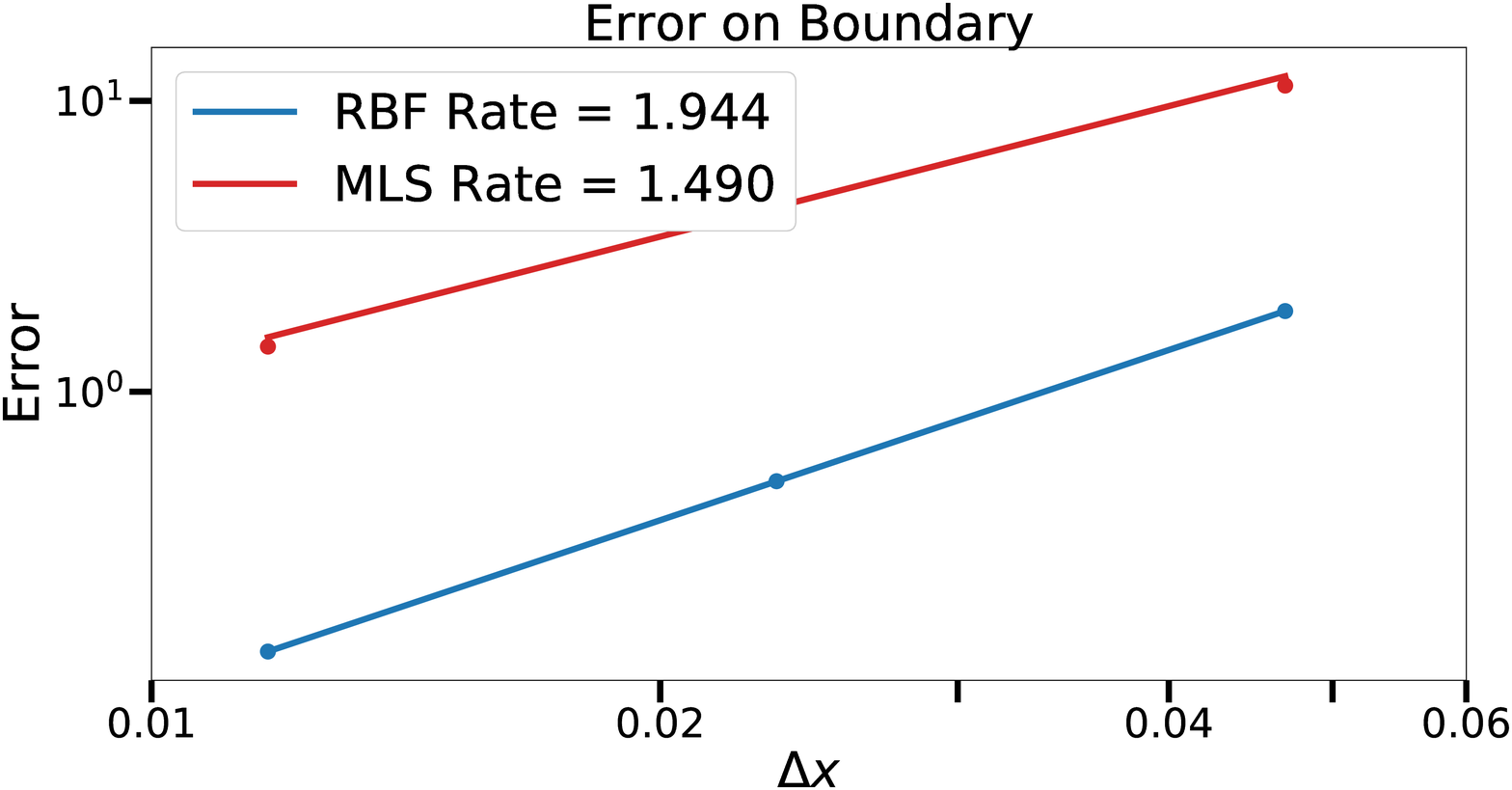}
\end{center}
\caption{Convergence rates for advecting and diffusion concentration field initialized from a point source in a solid body rotational flow using either a RBF (\subref{fig:adv_diff_rotational:rbf}) or MLS (\subref{fig:adv_diff_rotational:mls}) reconstruction. Also shown is the convergence rate for the integral in equation \eqref{eq:Error_moving} using both an RBF or MLS reconstruction (\subref{fig:adv_diff_rotational:bdry}). Simulations are run at an advective CFL number of $C_\text{CFL} = 0.5$ to a final time of $T = 1$.}
\label{fig:adv_diff_rotational}
\end{figure}

The computational domain $\BB$ is $[-4,4]\times[-4,4]$ and is again discretized using $N$ points in each direction. We use a disk radius of $R = 1$ and initial center $\xx_c(0) = \parens{1.521, 1.503}$. We use a diffusion coefficient of $D = 0.1$. The prescribed velocity is $\uu = 2\pi\parens{-y,x}$. Simulations are run at an advective CFL number of $0.5$ to a final time of $T = 1$. Figure \ref{fig:rotational_init_final} shows the initial and final solution. We can see that the solution maintains symmetry throughout the simulation. Figure \ref{fig:adv_diff_rotational} shows convergence rates for both RBF and MLS reconstructions. Figures \ref{fig:adv_diff_rotational:rbf} and \ref{fig:adv_diff_rotational:mls} show convergence rates in which the level set is prescribed. We see second order convergence rates in all norms for the prescribed level set case. Figure \ref{fig:adv_diff_rotational:bdry} shows the convergence rate for the solution extrapolated to the boundary for the case in which the level set is evolved with the fluid velocity. In this case, we see approximately second order convergence rates. In both of these cases, the error for the RBF reconstruction is almost an order of magnitude lower than that of the MLS reconstruction.

\subsection{Advection diffusion in oscillatory Couette flow}
We now consider the advection and diffusion of a concentration inside a vesicle under oscillating Couette flow with no flux boundary conditions. Specifically, we specify the rotational component of the velocity field as
\begin{subequations}
\begin{align}
u_\theta\parens{r} &= \parens{a r + \frac{b}{r}}\sin\parens{\pi t}, \\
a &= \frac{\Omega_2 R_2^2 - \Omega_1 R_1^2}{R_2^2 - R_1^2}, \\
b &= \frac{\parens{\Omega_1 - \Omega_2}R_1^2R_2^2}{R_2^2 - R_1^2},
\end{align}
\end{subequations}
where $R_1$ and $R_2$ are the radii of the coaxial cylinders, and $\Omega_1$ and $\Omega_2$ are the angular velocities of the cylinders. Here, we set the radii to be $R_1 = 0.5$ and $R_2 = 3.75$ with respective angular velocities $\Omega_1 = 10.0$ and $\Omega = 1.0$.

The vesicle is described by the level set $\phi_\text{v}\parens{\xx,t}$ so that the initial condition is
\begin{equation}
\phi_\text{v}\parens{\xx,0} = \norm*{\xx - \xx_{\text{v},\text{c}}} - R
\end{equation}
where $R = 1$ is the radius of the disk, and $\xx_{\text{v},\text{c}} = \parens{1.521,1.503}$ is the center of the disk. We initialize the concentration field $q_\text{v}$ via
\begin{equation}\label{eq:no_flux_inside}
q_\text{v}\parens{\xx,0} = \parens{\cos\parens{\pi \norm*{\xx-\xx_{\text{v},\text{c}}}} + 1}^2
\end{equation}

 For comparison, we define another concentration field $q_\text{o}\parens{\xx,t}$ outside the vesicle, but inside the cylinders. The second level set is defined by
\begin{equation}
\phi_\text{o}\parens{\xx,t} = \max\parens{-\phi_\text{v}\parens{\xx,t},\norm*{\xx} - R_1, R_2 - \norm*{\xx}}.
\end{equation}
The initial condition is given by
\begin{equation}
q_\text{o}\parens{\xx,0} = \left\{\begin{array}{cc}
\parens{\cos\parens{\pi \norm*{\xx-\xx_{\text{o},\text{c}}}} + 1}^2 &\mbox{ if } \norm*{\xx-\xx_{\text{o},\text{c}}} \leq 1 \\
0 &\mbox{ otherwise,}
\end{array}\right.
\end{equation}
in which $\xx_{\text{o},\text{c}} = - \xx_{\text{v},\text{c}}$.

\begin{figure}
\begin{center}
\phantomsubcaption\label{fig:two_couette:0}\phantomsubcaption\label{fig:two_couette:1}
\phantomsubcaption\label{fig:two_couette:2}\phantomsubcaption\label{fig:two_couette:3}
\phantomsubcaption\label{fig:two_couette:4}
\psfragscanon
\psfrag{4}{$4.0$}
\psfrag{2277}{$2.28$}
\psfrag{127}{$1.27$}
\psfrag{5546}{$0.55$}
\psfrag{0}{$0.0$}
\psfrag{2}{$2.0$}
\psfrag{1193}{$1.19$}
\psfrag{6781}{$0.68$}
\psfrag{2996}{$0.30$}
\subfigimg[width=0.45\linewidth,hsep=0em,vsep=2em,pos=ul]{(a)}{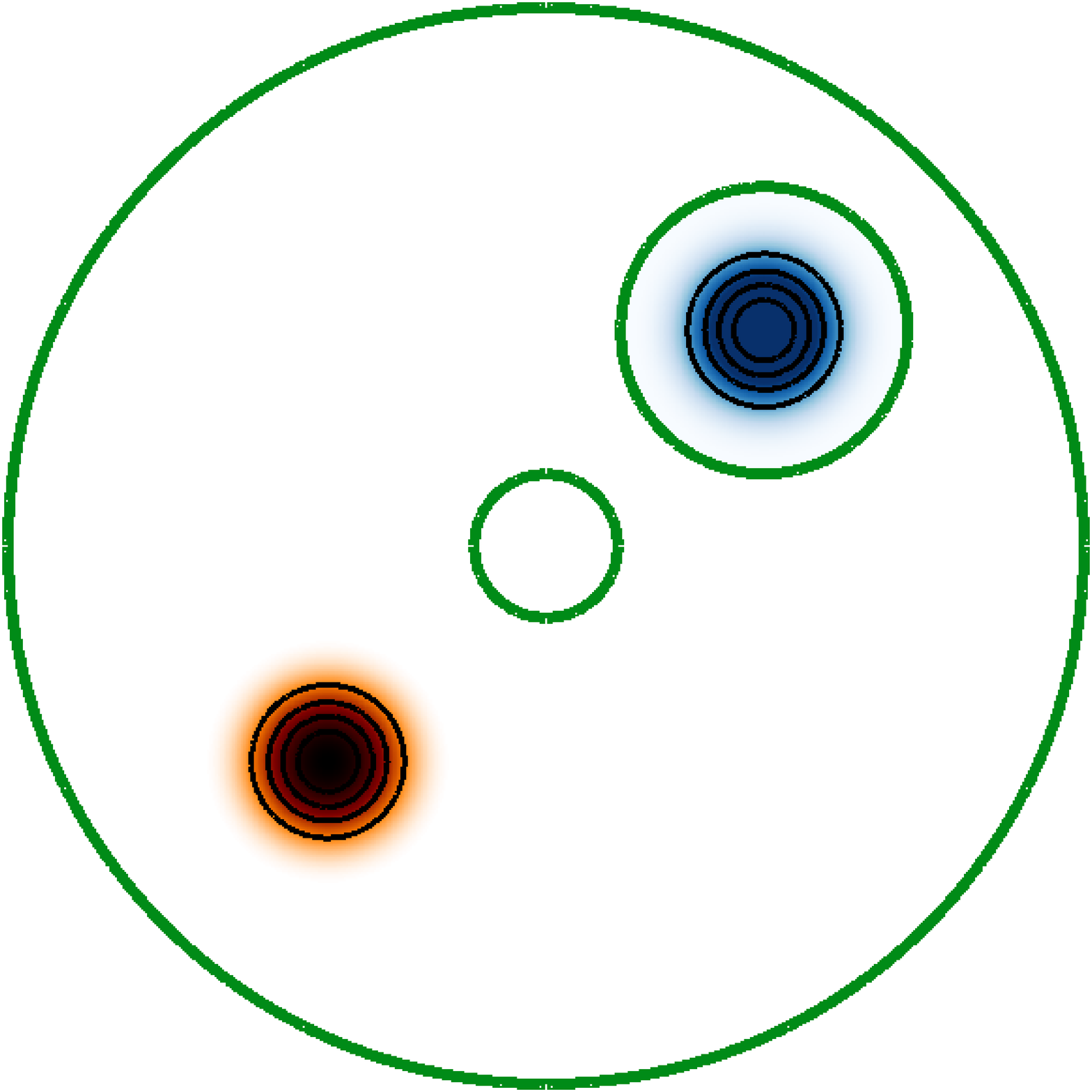}\quad
\subfigimg[width=0.45\linewidth,hsep=0em,vsep=2em,pos=ul]{(b)}{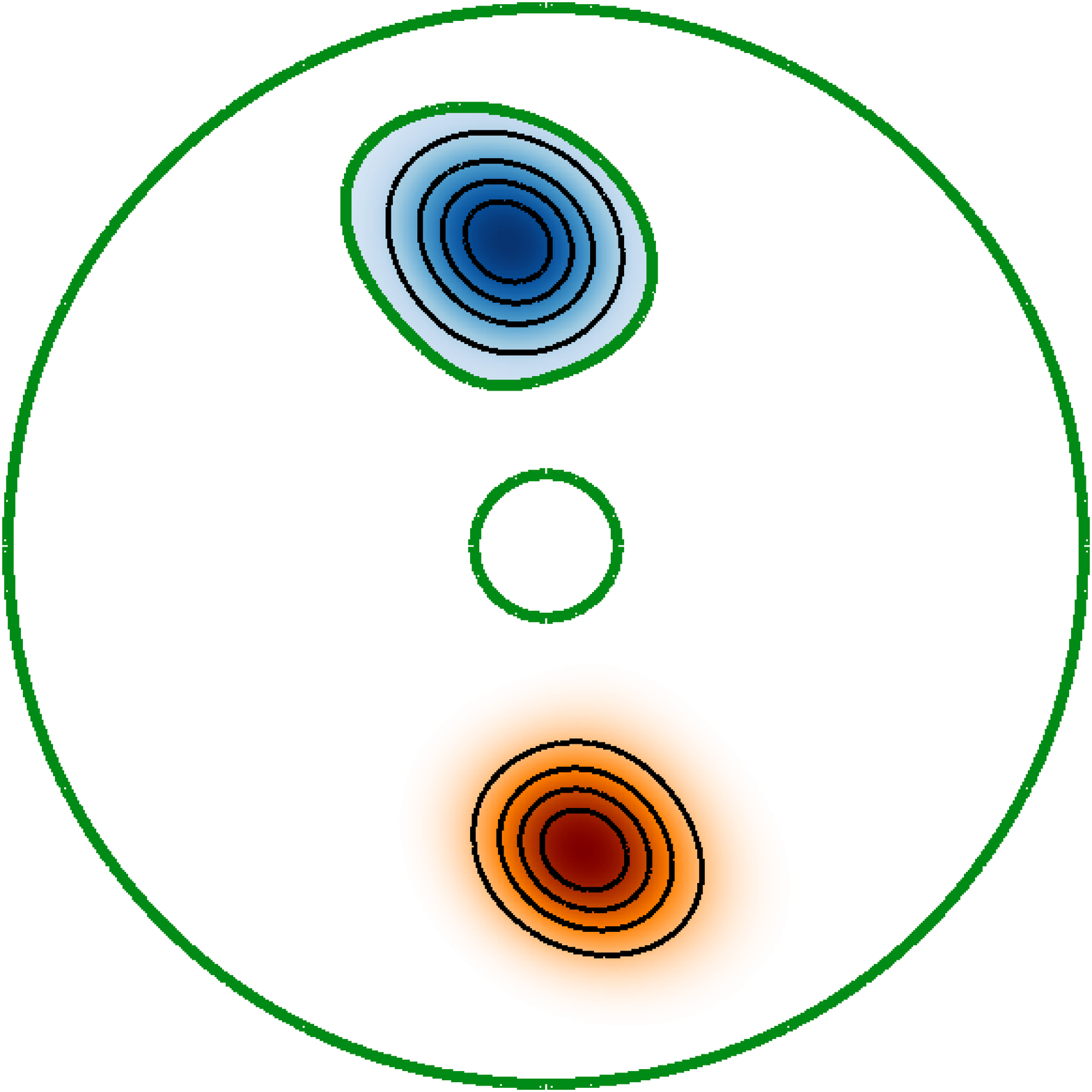}\\
\subfigimg[width=0.45\linewidth,hsep=0em,vsep=2em,pos=ul]{(c)}{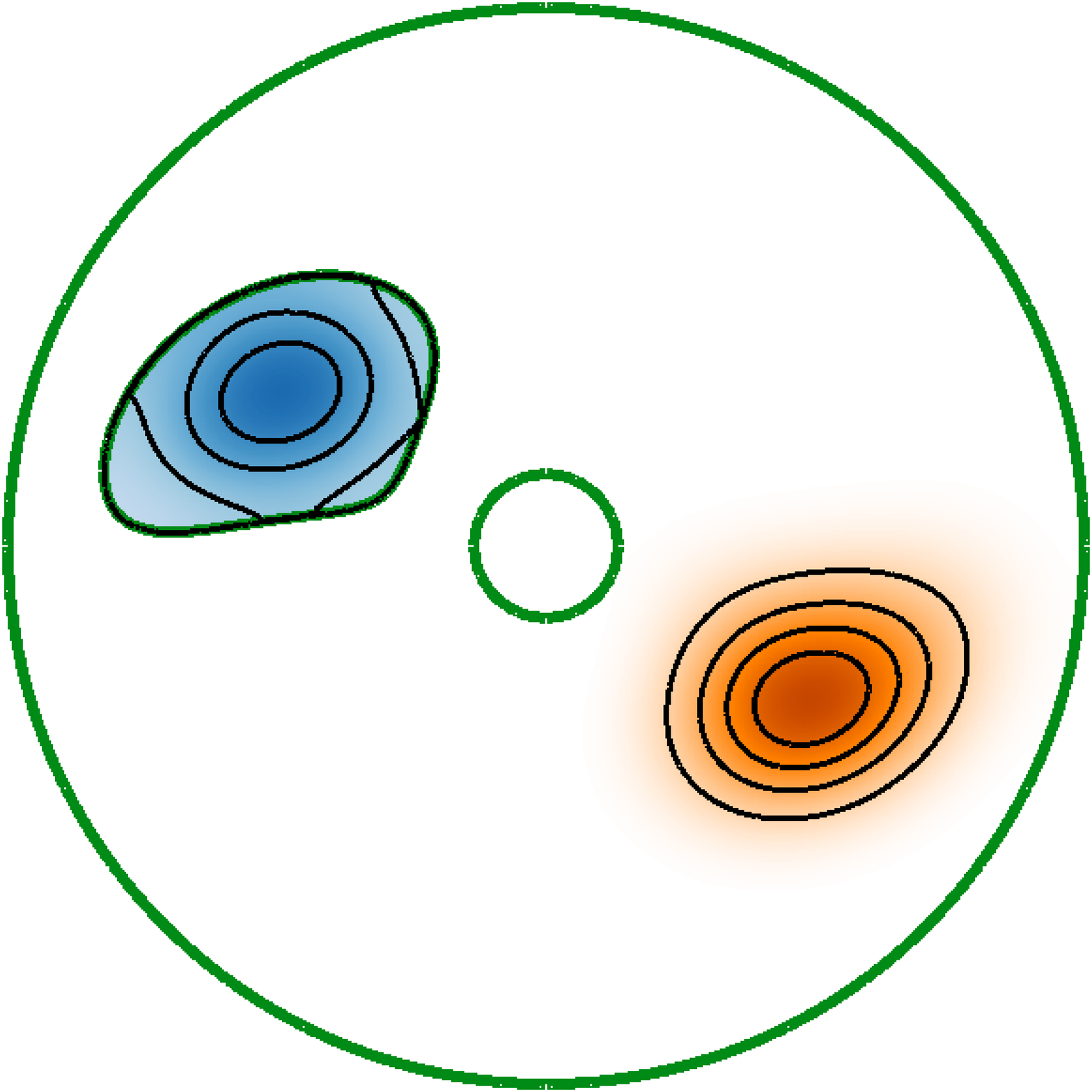}\quad
\subfigimg[width=0.45\linewidth,hsep=0em,vsep=2em,pos=ul]{(d)}{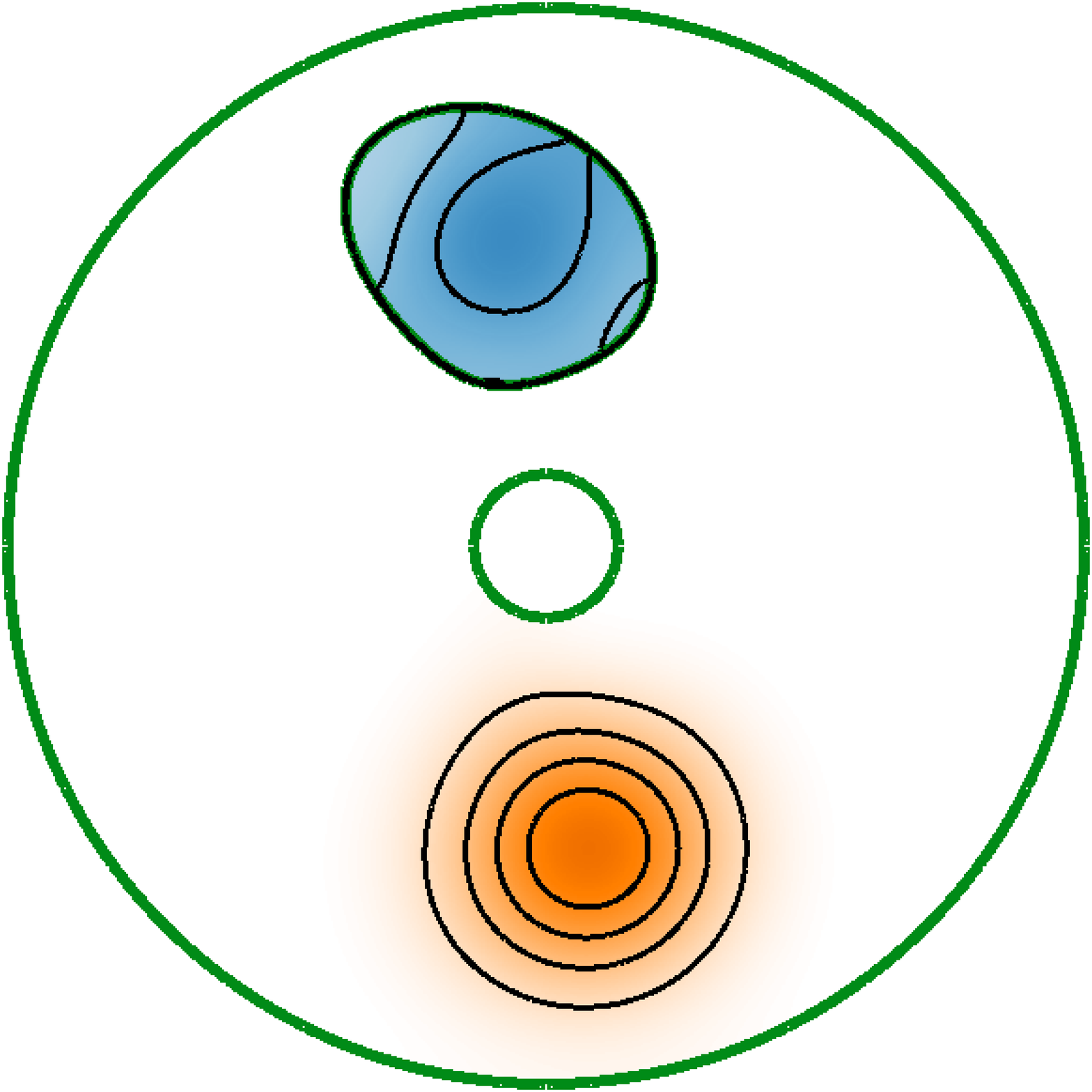}\\
\subfigimg[width=0.45\linewidth,hsep=0em,vsep=2em,pos=ul]{(e)}{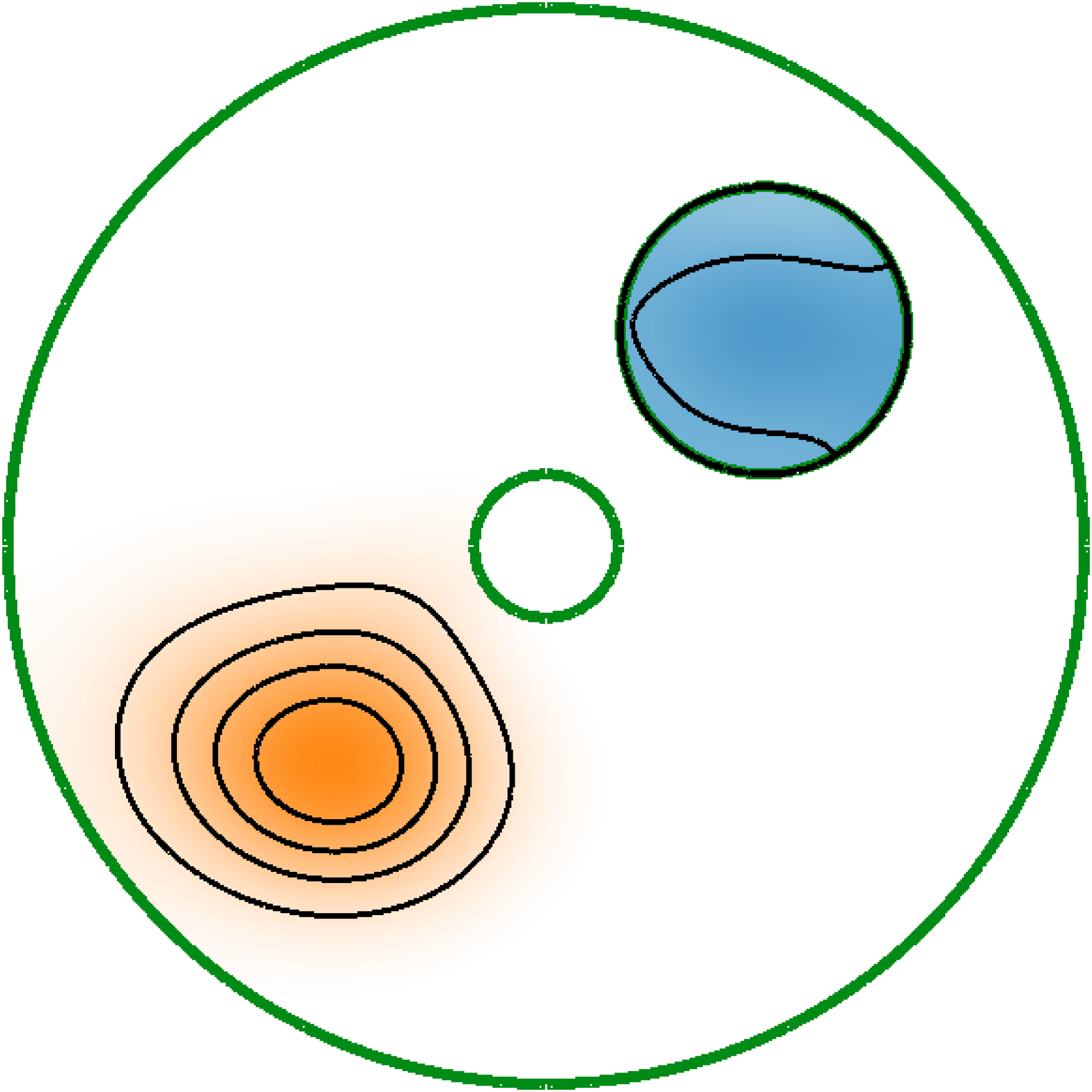}
\includegraphics[width=0.2\linewidth]{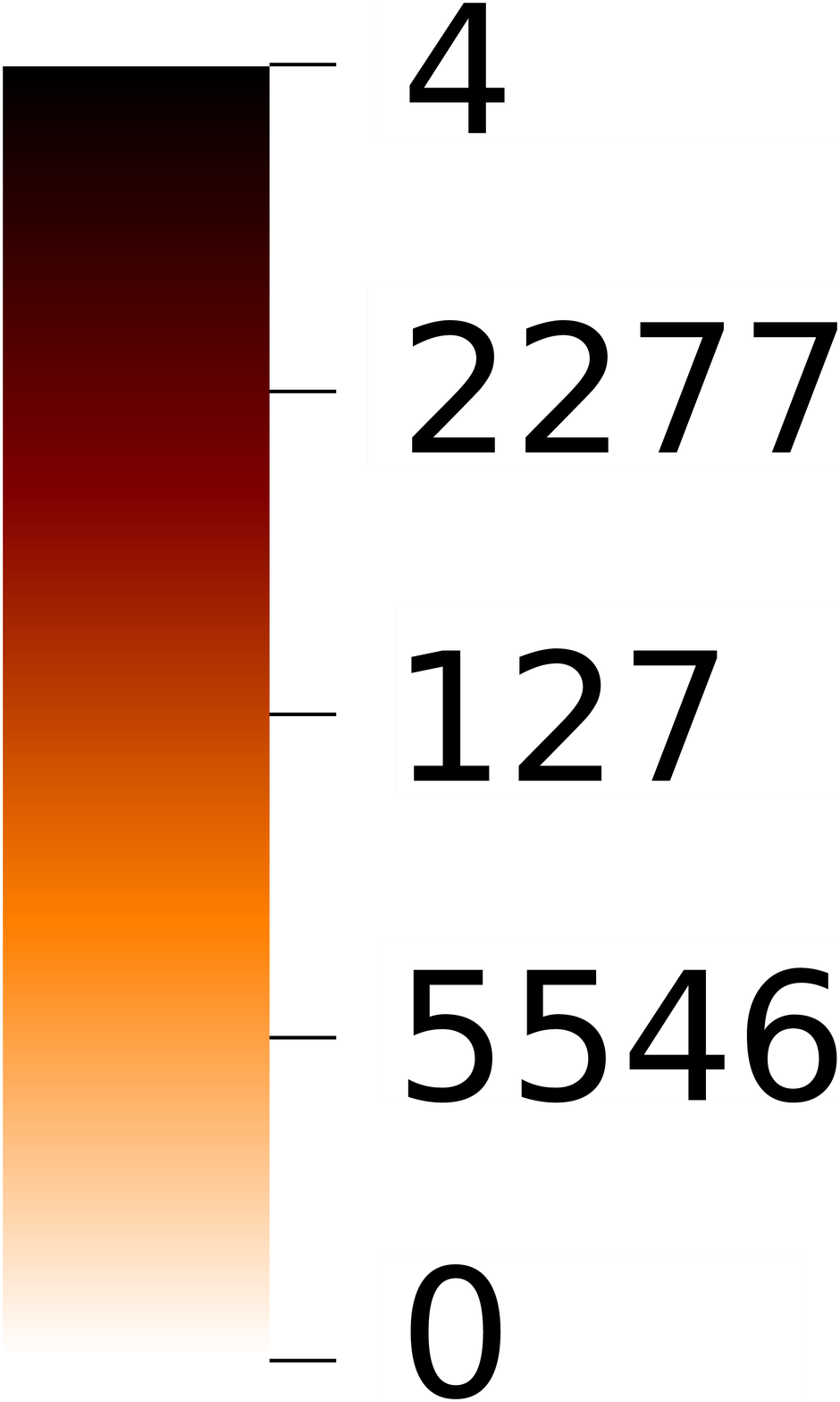}
\includegraphics[width=0.2\linewidth]{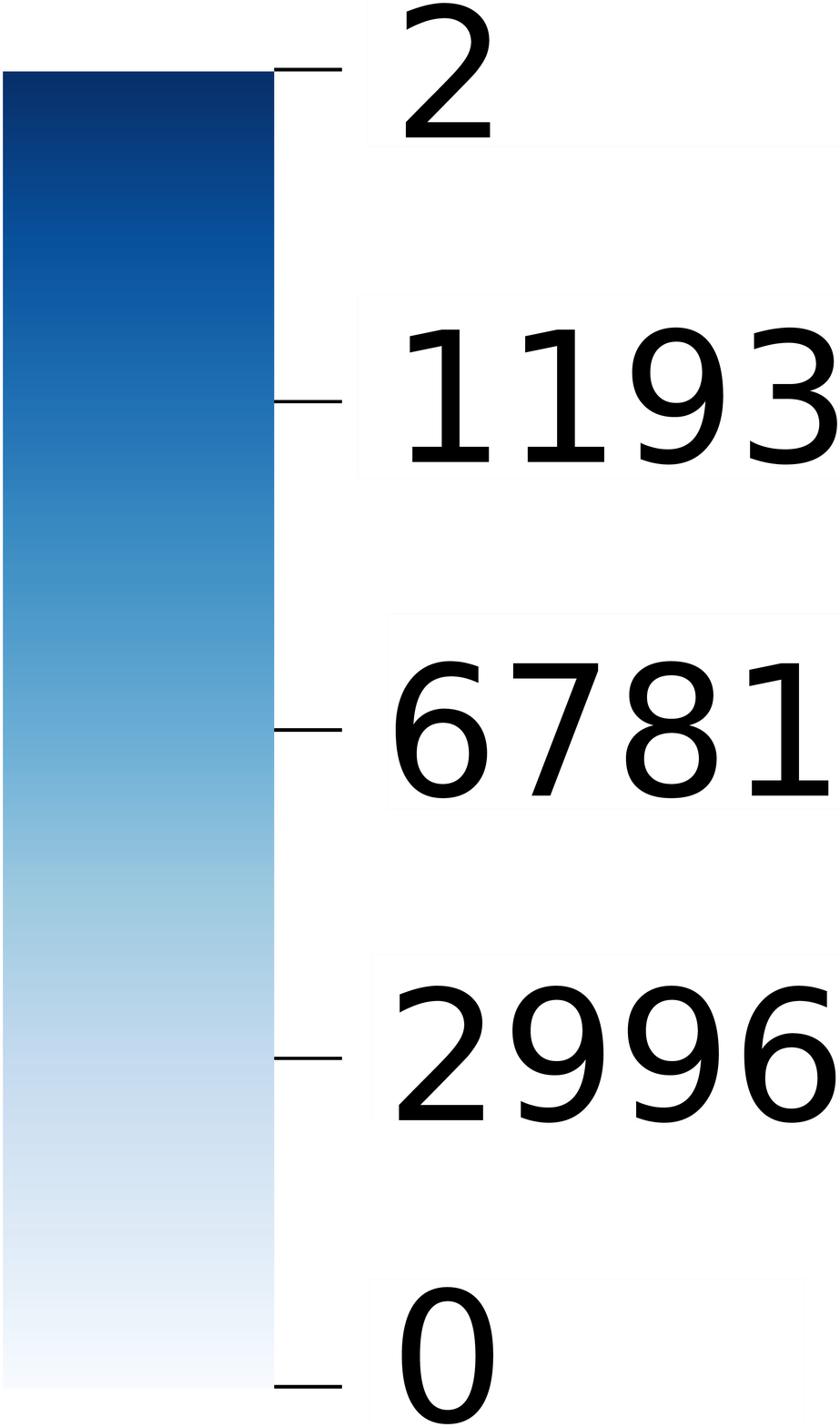}
\end{center}
\caption{Concentrations $q_\text{o}$ (orange) and $q_\text{v}$ (blue) at time points $0.0$ (\subref{fig:two_couette:0}), $0.5$ (\subref{fig:two_couette:1}), $1.0$ (\subref{fig:two_couette:2}), $1.5$ (\subref{fig:two_couette:3}), and $2.0$ (\subref{fig:two_couette:4}). The zero contour of the level set $\phi_\text{o}$ is shown in green. The diffusion coefficient is fixed at $D = 0.05$ with $1024$ points in each direction. We use a fixed time-step corresponding to a maximum CFL number of $C_\text{CFL} = 0.5$ to a final time of $T = 2$.}
\label{fig:two_couette}
\end{figure}

The computational domain $\BB$ is the box $[-4, 4]\times [-4,4]$ and is discretized using $N$ points in each direction, and the diffusion coefficient for both the interior and exterior concentration fields is $D = 0.05$. Figure \ref{fig:two_couette} shows the solution at five different time values during the simulation. Over the course of the simulations, the concentration field inside the disk stays completely contained within the disk, while the outside concentration field is free to diffuse between the coaxial cylinders.

\begin{figure}
\begin{center}
\phantomsubcaption\label{fig:couette_in_converge:l1_rate}\phantomsubcaption\label{fig:couette_in_converge:l2_rate}
\phantomsubcaption\label{fig:couette_in_converge:max_rate}\phantomsubcaption\label{fig:couette_in_converge:bdry_l1_rate}
\subfigimg[width=0.45\linewidth,hsep=-1em,vsep=2em,pos=ul]{(a)}{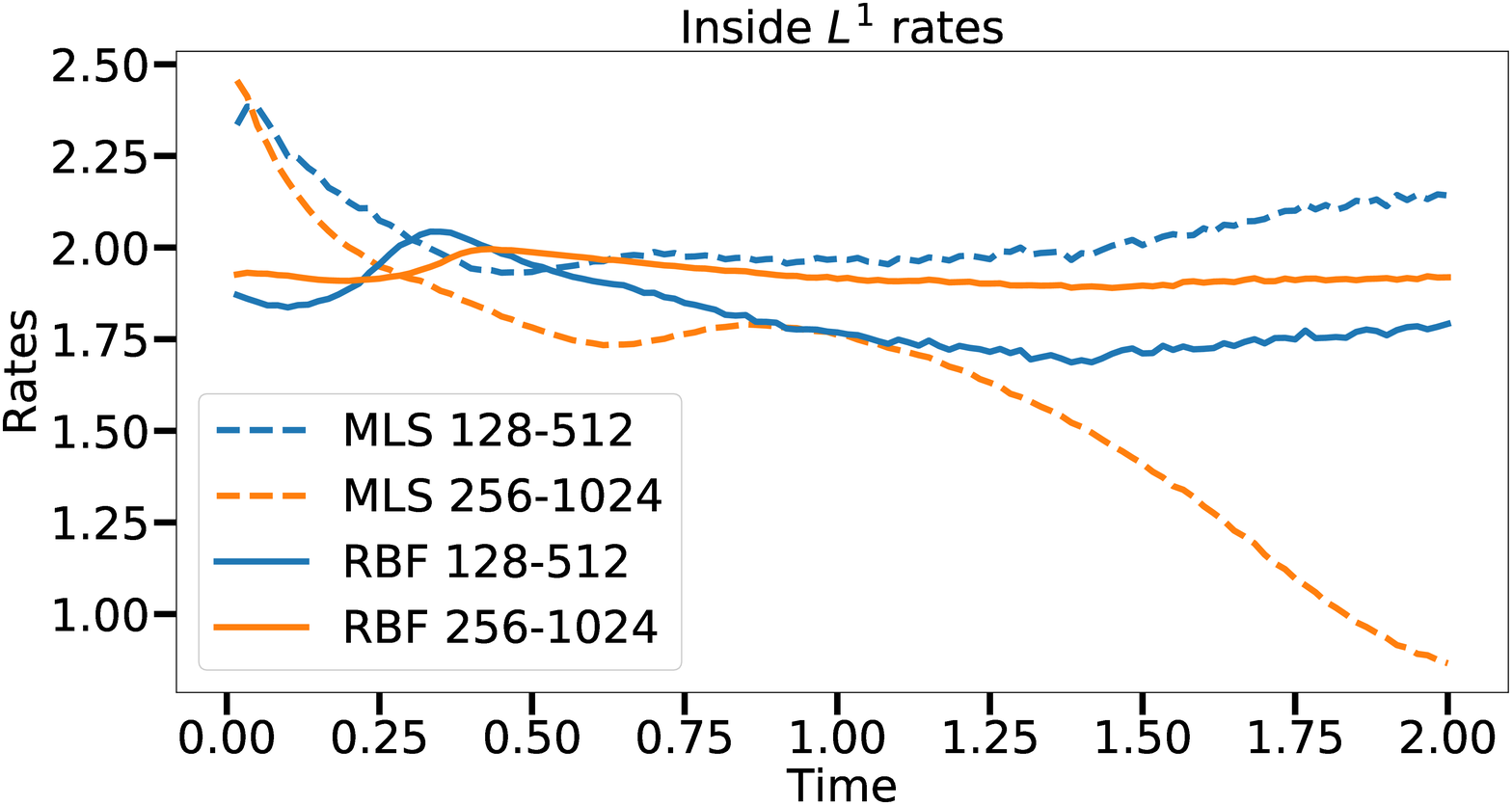}\quad
\subfigimg[width=0.45\linewidth,hsep=-1em,vsep=2em,pos=ul]{(b)}{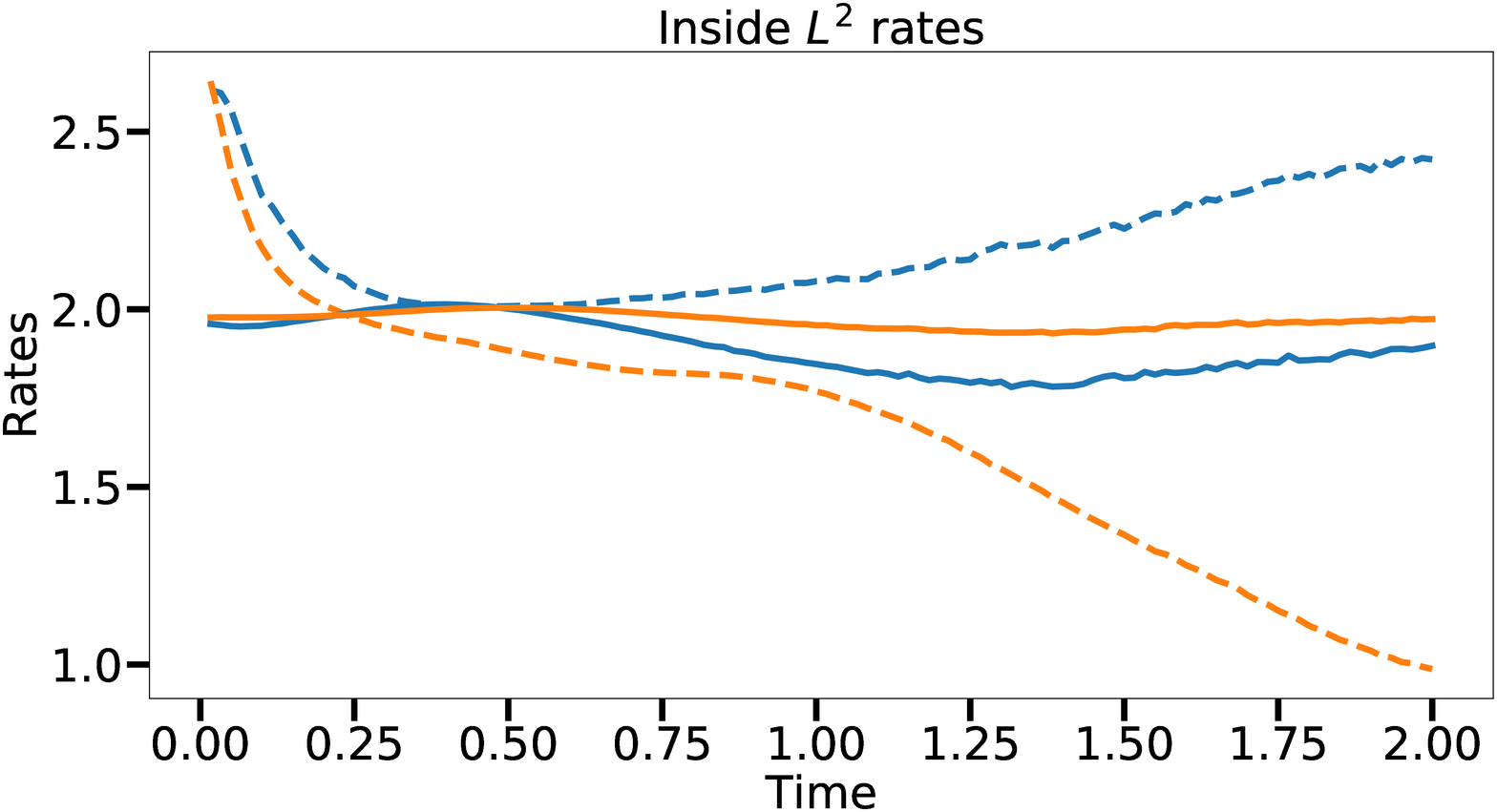}\\
\subfigimg[width=0.45\linewidth,hsep=-1em,vsep=2em,pos=ul]{(c)}{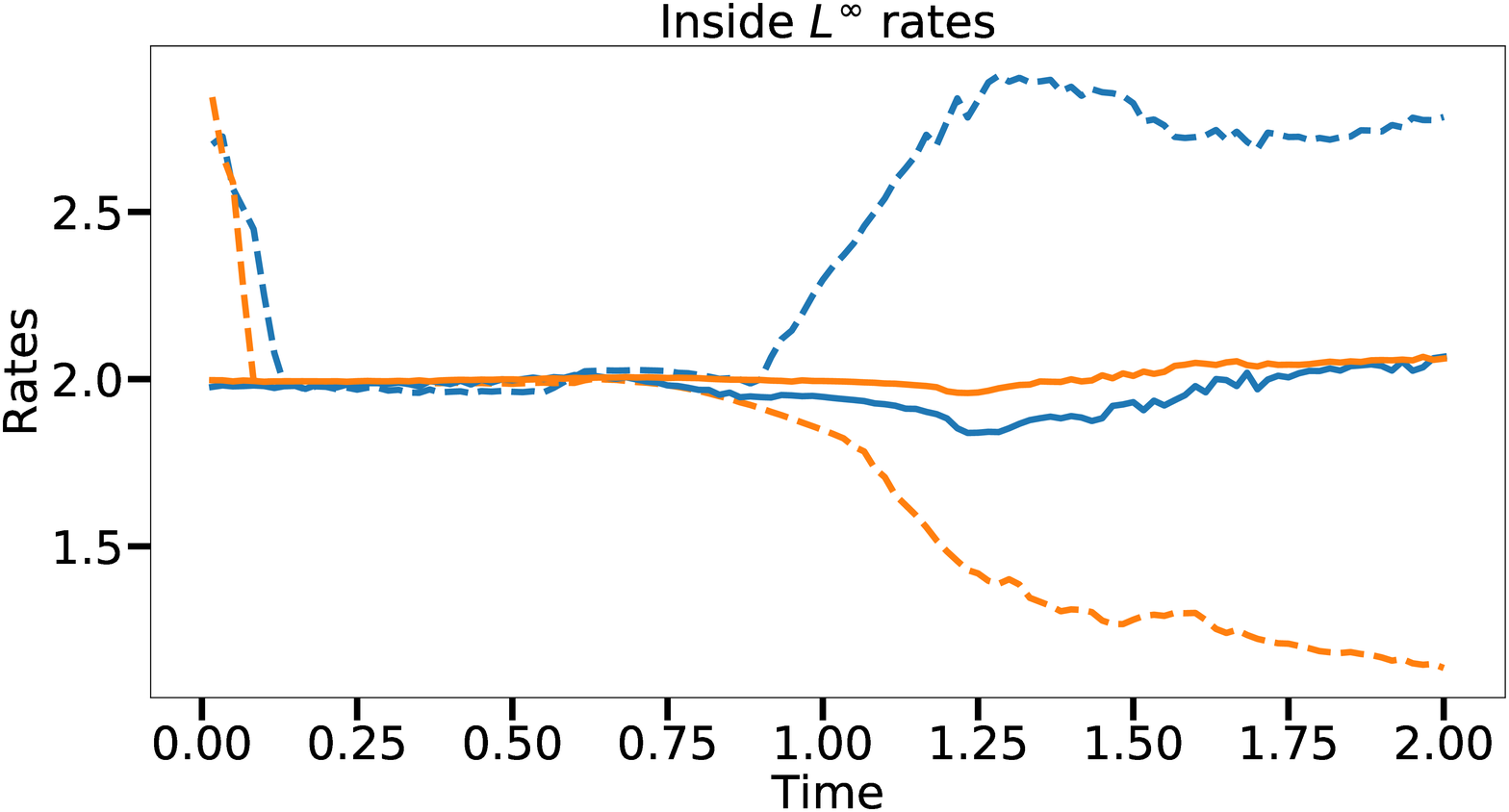}\quad
\subfigimg[width=0.45\linewidth,hsep=-1em,vsep=2em,pos=ul]{(d)}{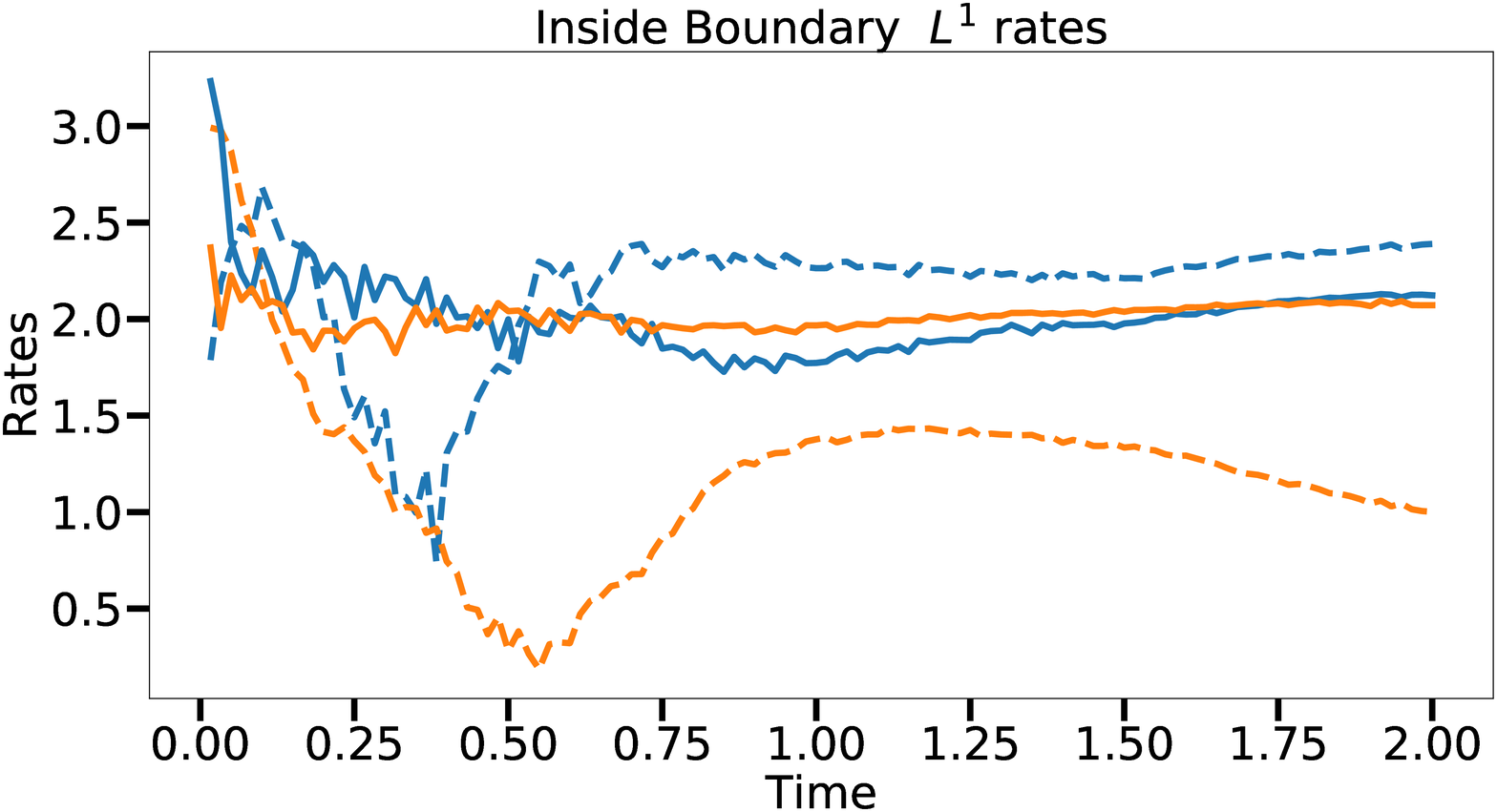}
\end{center}
\caption{Numerical convergence study for the interior of the disk in the $L^1$ (\subref{fig:couette_in_converge:l1_rate}), $L^2$ (\subref{fig:couette_in_converge:l2_rate}), and $L^\infty$ (\subref{fig:couette_in_converge:max_rate}) norms. The convergence rates are computed from simulations with $N = 128,256,512$ points (blue) and $N = 256,512,1024$ points (orange). Also shown is the convergence rate for the extrapolation of the solution to the boundary (\subref{fig:couette_in_converge:bdry_l1_rate}).}
\label{fig:couette_in_converge}
\end{figure}

\begin{figure}
\begin{center}
\phantomsubcaption\label{fig:couette_out_converge:l1_rate}\phantomsubcaption\label{fig:couette_out_converge:l2_rate}
\phantomsubcaption\label{fig:couette_out_converge:max_rate}
\subfigimg[width=0.45\textwidth,hsep=-1em,vsep=2em,pos=ul]{(a)}{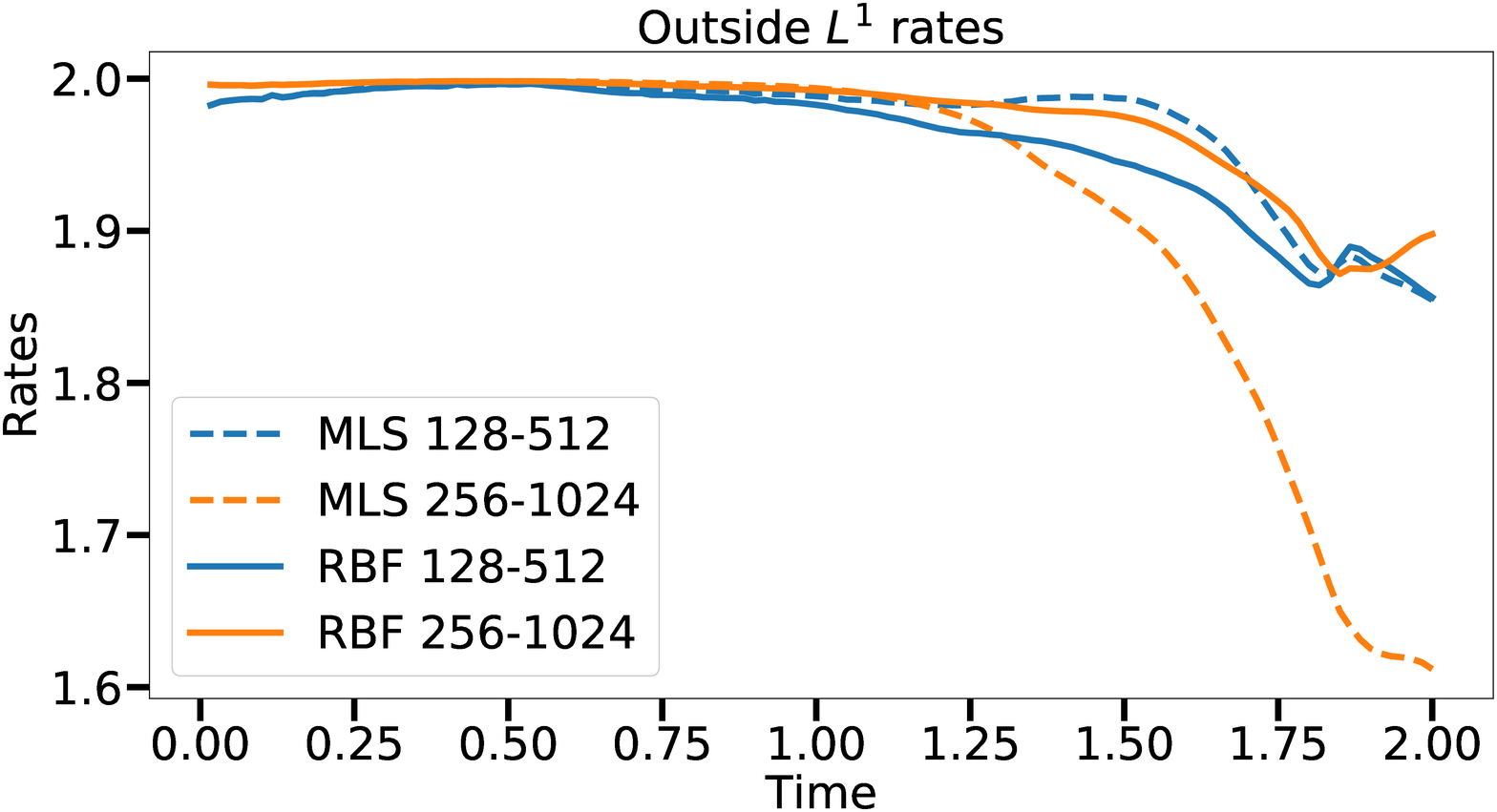}
\subfigimg[width=0.45\textwidth,hsep=-1em,vsep=2em,pos=ul]{(b)}{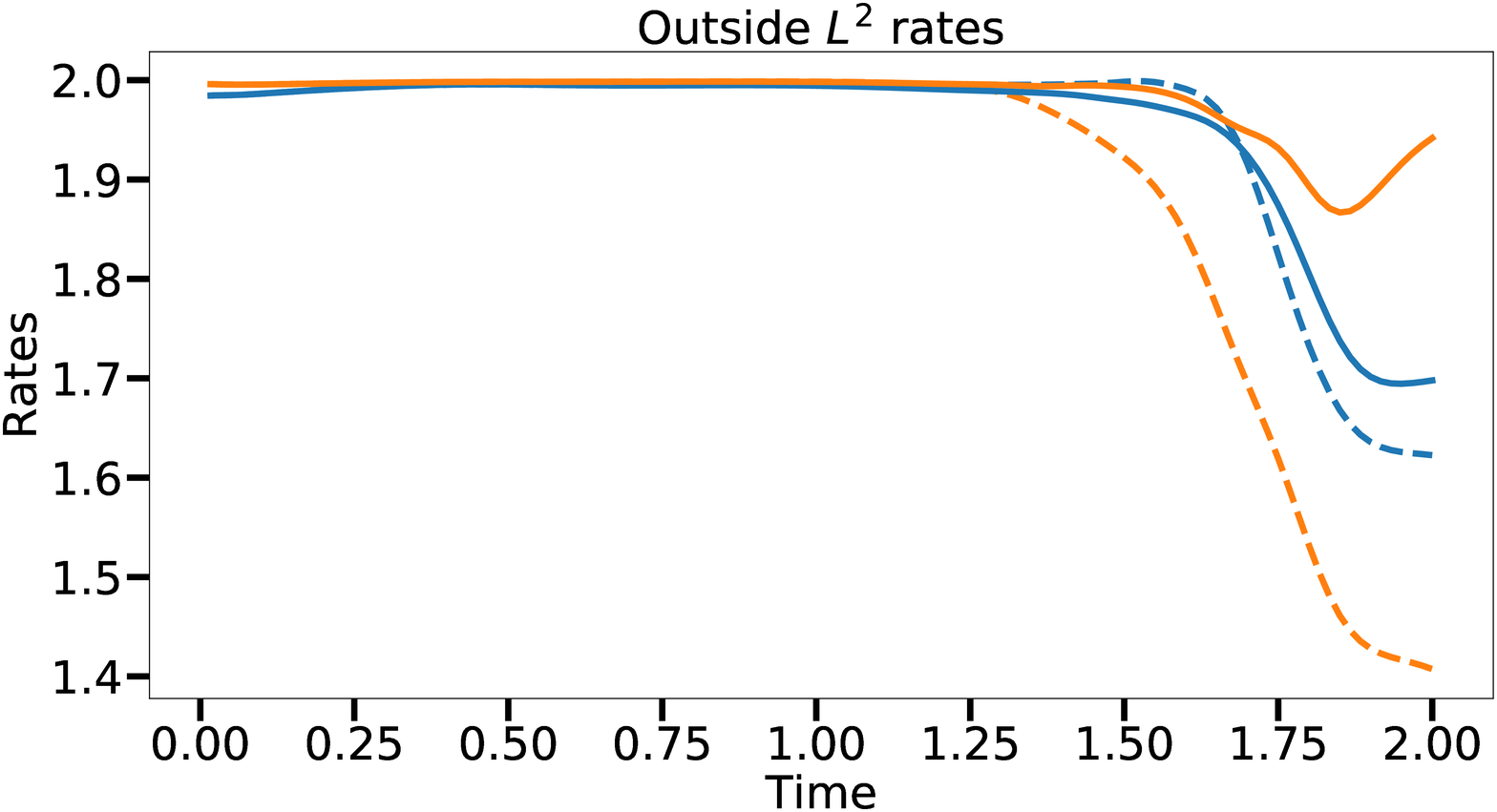}
\subfigimg[width=0.45\textwidth,hsep=-1em,vsep=2em,pos=ul]{(c)}{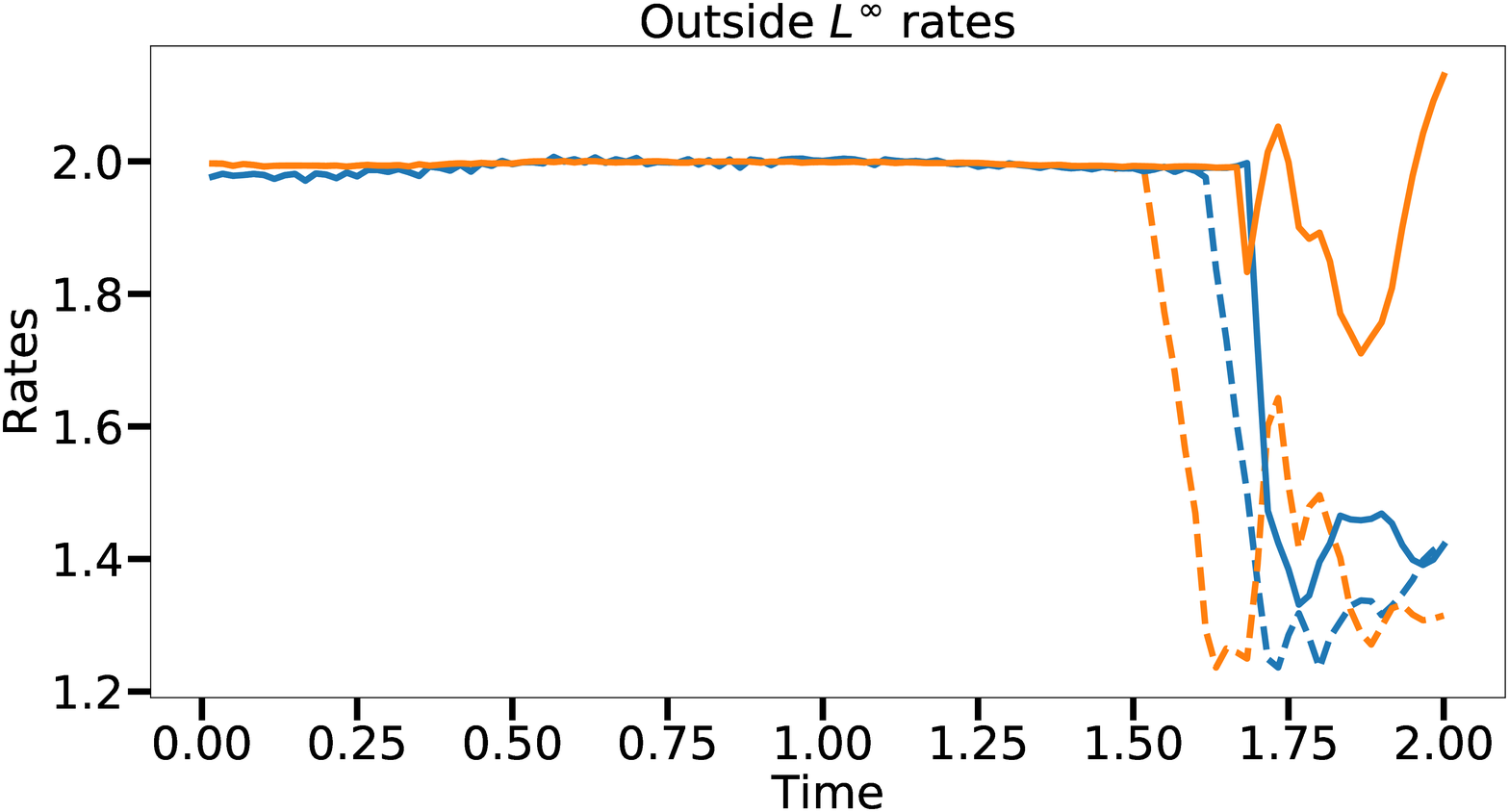}
\end{center}
\caption{Numerical convergence study for the exterior of the disk but inside the two cylinders in the $L^1$ (\subref{fig:couette_out_converge:l1_rate}), $L^2$ (\subref{fig:couette_out_converge:l2_rate}), and $L^\infty$ (\subref{fig:couette_out_converge:max_rate}) norms. The convergence study is computed from simulations with $N = 128,256,512$ points and $N = 256,512,1024$ points. The dip in convergence rates near the end of the simulation occur when there is a significant value of $q_\text{o}$ near the boundaries.}
\label{fig:couette_out_converge}
\end{figure}

\begin{figure}
\begin{center}
\phantomsubcaption\label{fig:couette_in_norms:l1_norms}\phantomsubcaption\label{fig:couette_in_norms:l2_norms}
\phantomsubcaption\label{fig:couette_in_norms:max_norms}\phantomsubcaption\label{fig:couette_in_norms:bdry_norms}
\subfigimg[width=0.45\textwidth,hsep=-1em,vsep=2em,pos=ul]{(a)}{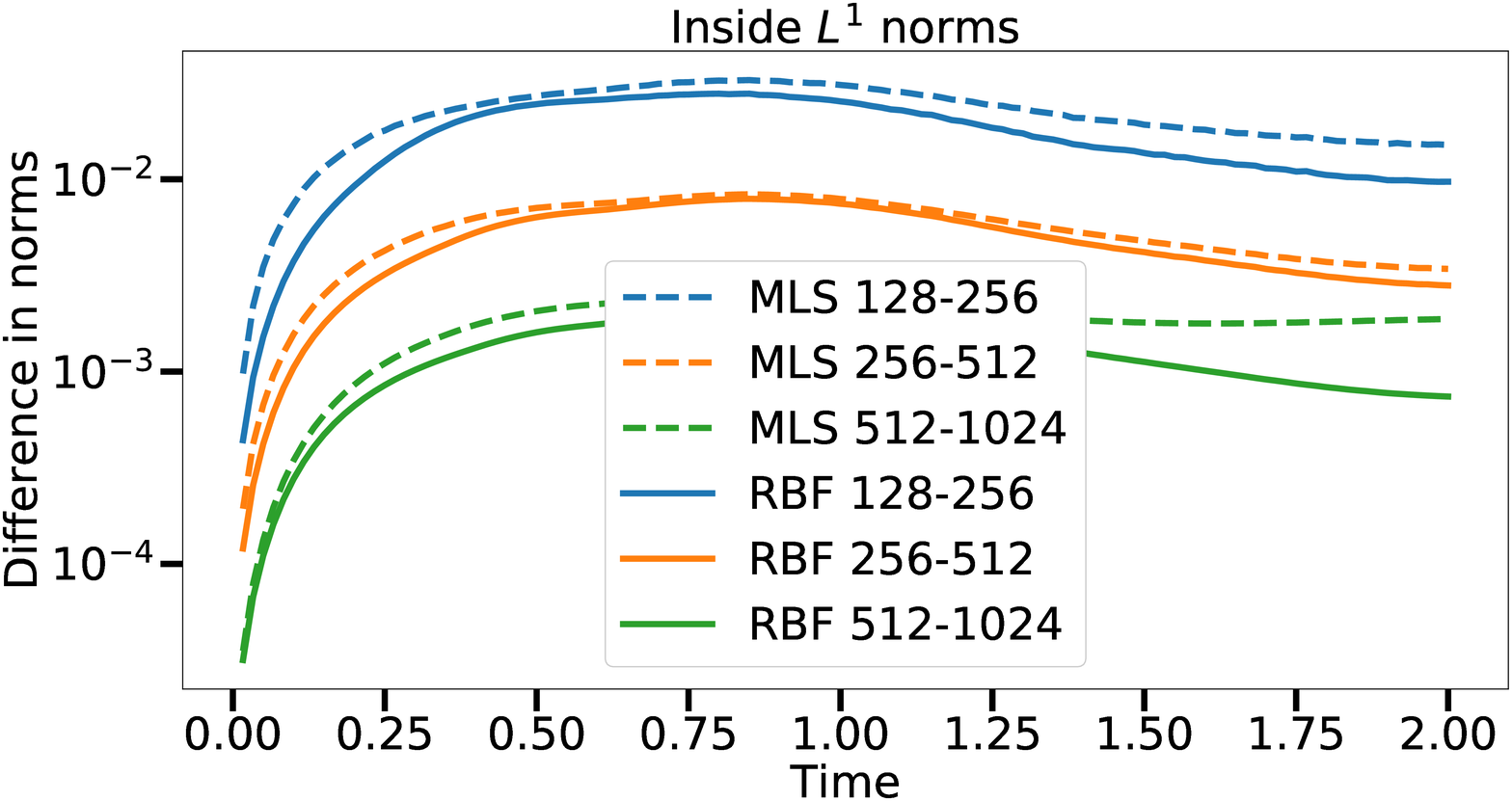}
\subfigimg[width=0.45\textwidth,hsep=-1em,vsep=2em,pos=ul]{(b)}{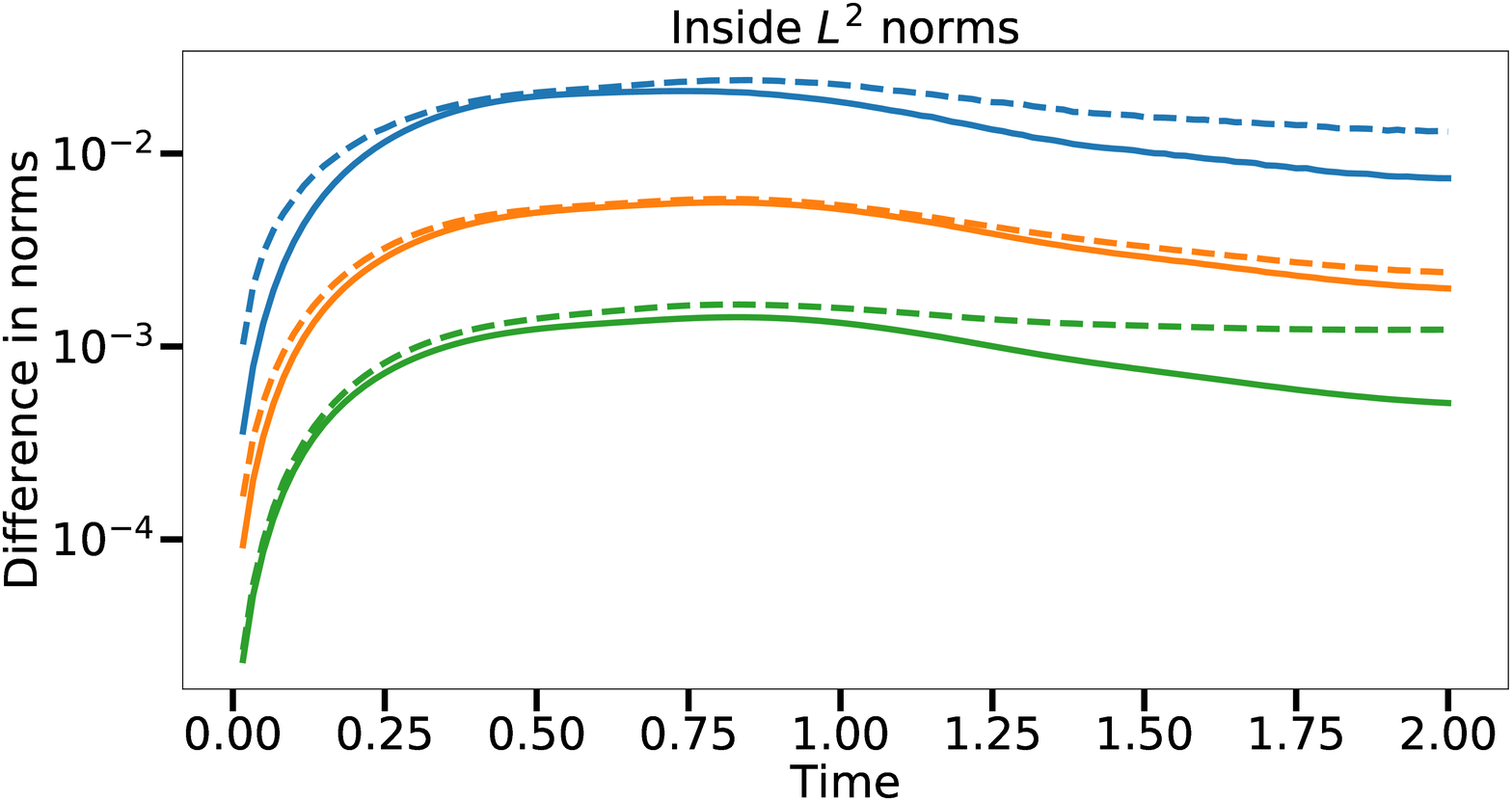}\\
\subfigimg[width=0.45\textwidth,hsep=-1em,vsep=2em,pos=ul]{(c)}{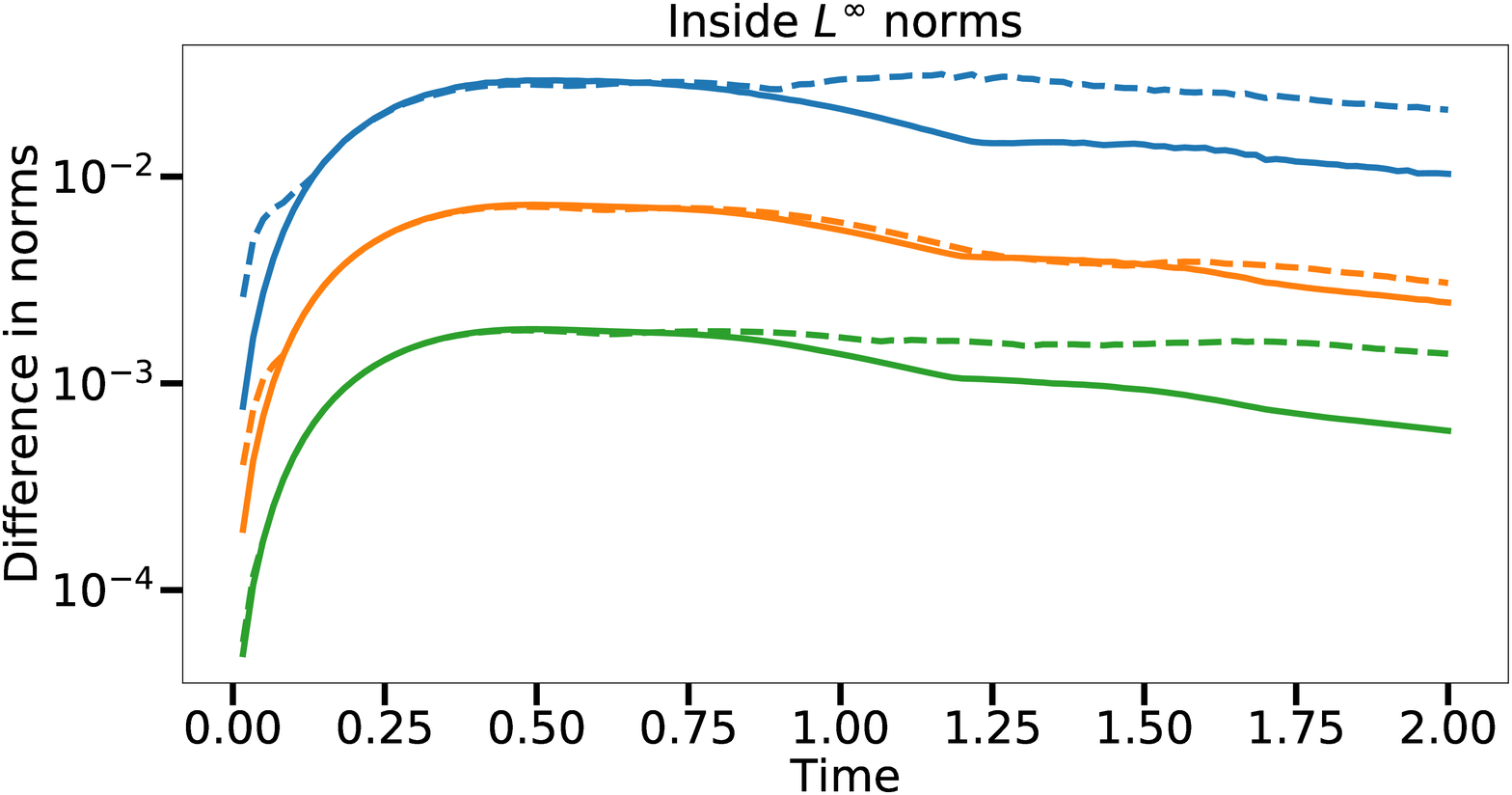}
\subfigimg[width=0.45\textwidth,hsep=-1em,vsep=2em,pos=ul]{(d)}{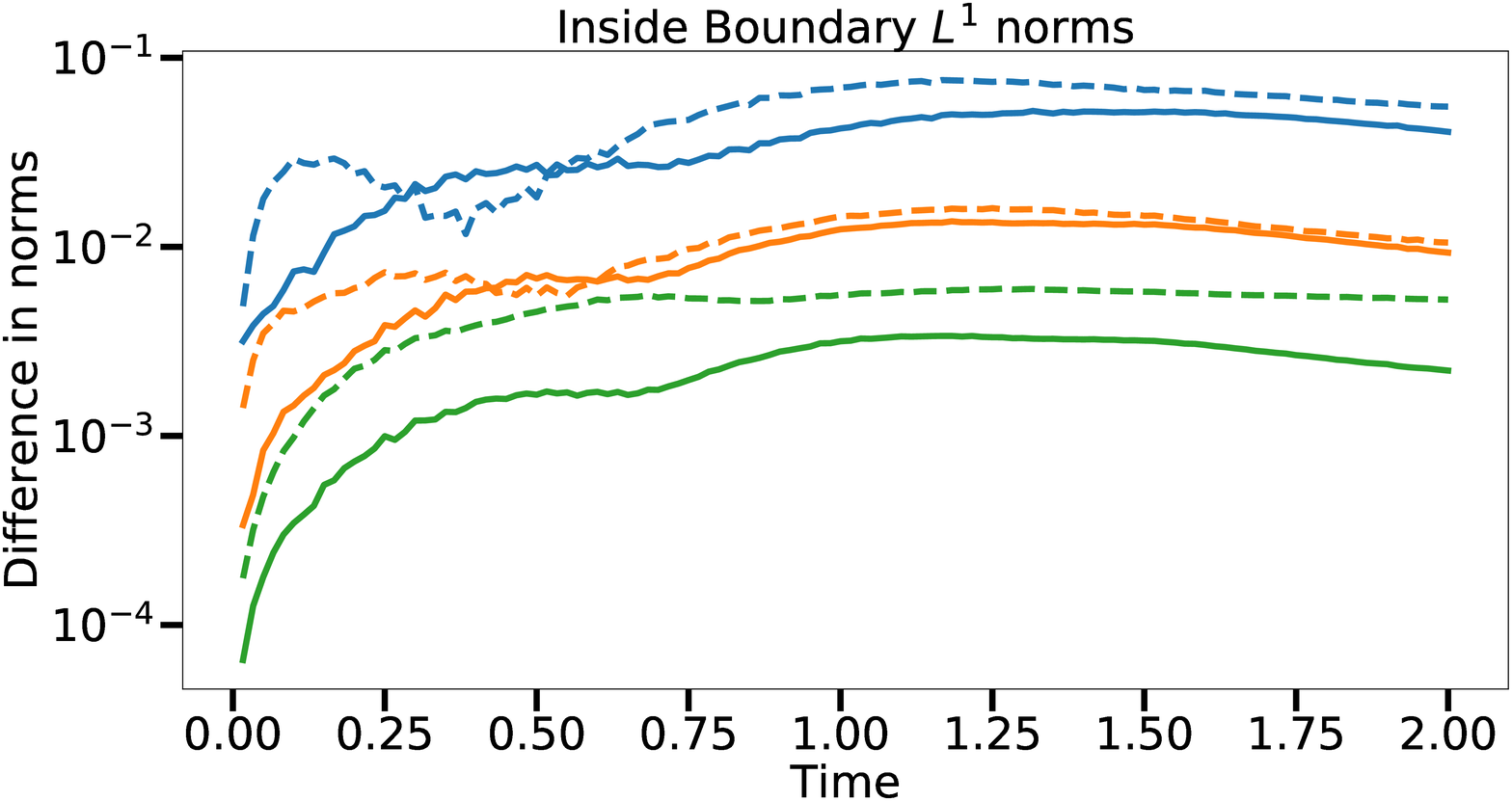}
\end{center}
\caption{Norms of solution differences for the concentration field inside the domain for the $L^1$ (\subref{fig:couette_in_norms:l1_norms}), $L^2$ (\subref{fig:couette_in_norms:l2_norms}), and $L^\infty$ (\subref{fig:couette_in_norms:max_norms}) norms. The dashed lines correspond to MLS reconstruction while the solid lines correspond to a RBF reconstruction. Also shown are the norms of solution differences for the extrapolation to the boundary (\subref{fig:couette_in_norms:bdry_norms}).}
\label{fig:couette_in_norms}
\end{figure}

\begin{figure}
\begin{center}
\phantomsubcaption\label{fig:couette_out_norms:l1_norms}\phantomsubcaption\label{fig:couette_out_norms:l2_norms}
\phantomsubcaption\label{fig:couette_out_norms:max_norms}
\subfigimg[width=0.45\textwidth,hsep=-1em,vsep=2em,pos=ul]{(a)}{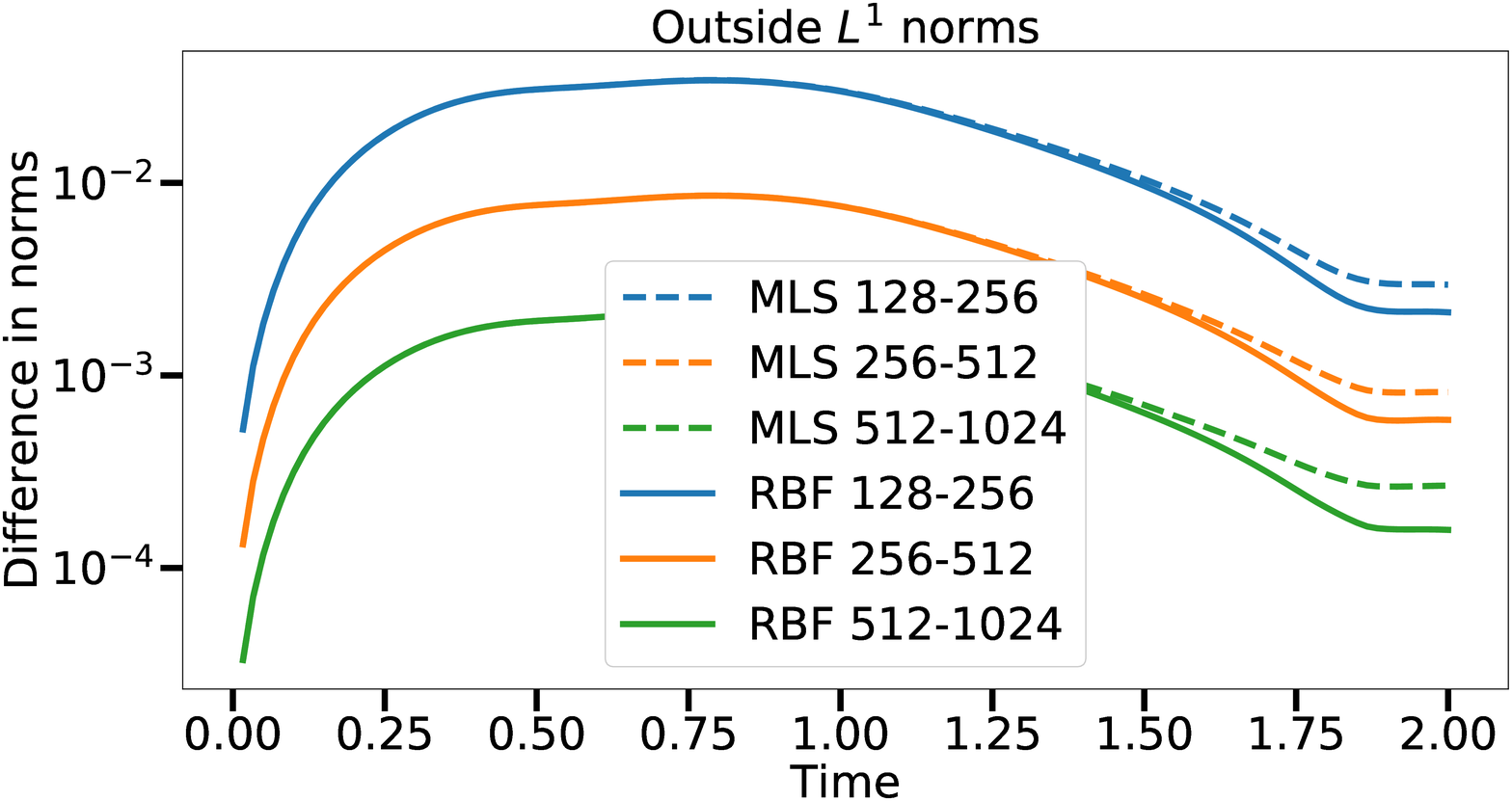}
\subfigimg[width=0.45\textwidth,hsep=-1em,vsep=2em,pos=ul]{(b)}{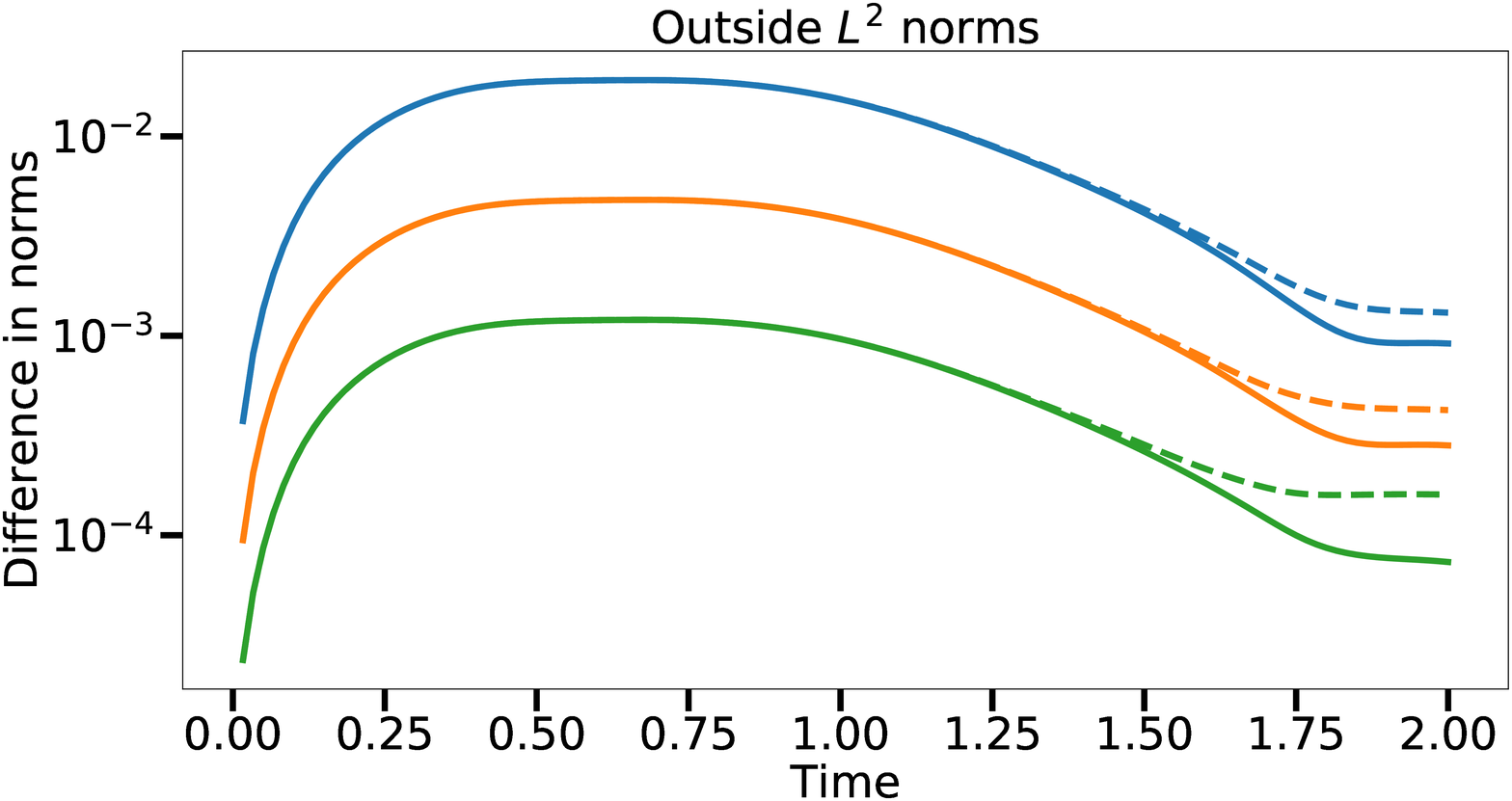}\\
\subfigimg[width=0.45\textwidth,hsep=-1em,vsep=2em,pos=ul]{(c)}{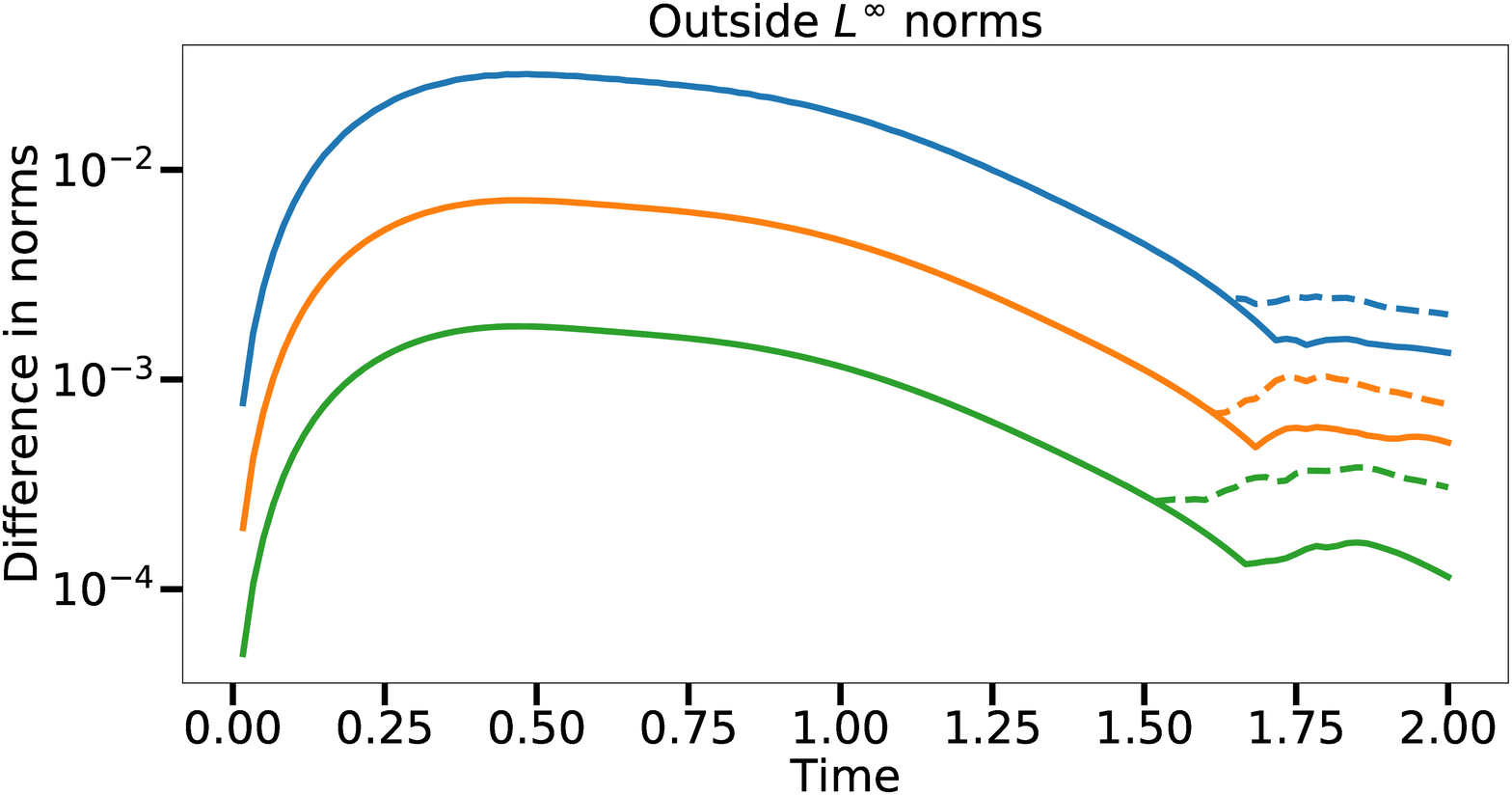}
\end{center}
\caption{Norms of solution differences for the concentration field outside the domain for the $L^1$ (\subref{fig:couette_out_norms:l1_norms}), $L^2$ (\subref{fig:couette_out_norms:l2_norms}), and $L^\infty$ (\subref{fig:couette_out_norms:max_norms}) norms. The dashed lines correspond to MLS reconstruction while the solid lines correspond to a RBF reconstruction.}
\label{fig:couette_out_norms}
\end{figure}

Figures \ref{fig:couette_in_converge} and \ref{fig:couette_out_converge} show numerical convergence results for the concentrations inside the vesicle and outside the vesicle but inside the two cylinders respectively. We estimate the convergence rate $r$ using Richardson extrapolation via
\begin{equation}
r = \log_2\frac{\norm*{q_{4h} - q_{2h}}}{\norm*{q_{2h} - q_{h}}},
\end{equation}
in which $q_h$ is the solution with a grid spacing of $h$. To compute the difference $\norm*{q_{2h} - q_{h}}$, we first interpolate the fine solution $q_h$ onto the coarser grid $q_{2h}$. Because the interpolation from a fine grid to a coarse grid is nontrivial near cut-cells, we exclude all cut-cells and any cells that neighbor cut-cells in the computation of the norm. To assess the error on the cut-cells, we perform a reconstruction of the the solution on the surface of the disk by extrapolating the solution to the disk in each cut-cell, then computing a cubic spline of these data points. Figure \ref{fig:couette_in_converge:bdry_l1_rate} shows convergence of this reconstruction for the solution inside the vessicle. We obtain second order convergence rates using RBFs at all points in time for the finest grids. In contrast, the MLS reconstruction does not show optimal convergence rates for the range of grid spacings considered. We expect that under additional grid refinement the MLS approximation will settle to second order accuracy, but this is not achieved over the range of grid spacings considered herein. The dips in convergence rates in Figure \ref{fig:couette_out_converge} occur when significant values of $q_\text{o}$ begin to appear in cut cells. Additionally, Figures \ref{fig:couette_in_norms} and \ref{fig:couette_out_norms} show the norms of the differences. To quantify the difference in accuracy between MLS and RBF reconstructions, we estimate the error coefficient $C$ in the expression
\begin{equation}
\norm*{q_h - q} \approx C h^p,
\end{equation}
in which $q$ is the exact solution and $p$ is the order of accuracy. We approximate $C$ by comparing differences between two grids of refinement $h$ and $2h$ through the following equation
\begin{equation}\label{eq:C}
C \approx \frac{\norm*{q_h - q_{2h}}}{h^p\parens{1 - 2^p}}.
\end{equation}
For both MLS and RBF reconstructions, we have $p = 2$. Table \ref{tab:c_coeff} shows our estimates of the $C$ coefficient for both reconstructions. In all cases, the RBF reconstruction is more accurate.

\begin{table}
\begin{center}
\begin{tabular}{cccc|ccc}
\hline
 & $q_\text{o}$ & & & $q_\text{v}$ & & \\
\hline
 & $L^1$ norm & $L^2$ norm & $L^\infty$ norm & $L^1$ norm & $L^2$ norm & $L^\infty$ norm \\
\hline
$C_\text{MLS}$ & 41.01 & 26.71 & 30.37 & 5.85 & 3.49 & 6.68\\
$C_\text{RBF}$  & 16.20  & 11.12 & 12.88 & 3.44 & 1.60 & 2.50\\
\hline
\end{tabular}
\end{center}
\caption{The leading order error coefficient $C$ for both MLS and RBF reconstructions is estimated using equation \eqref{eq:C}. In all cases, the RBF reconstruction yields superior accuracy.}\label{tab:c_coeff}
\end{table}

\subsection{Choice of time step size}
Previous sections use time step sizes based on a CFL number of $0.5$. We note that each substep consists of an unconditionally stable method. While the affect on stability of the splitting error is unclear, we briefly explore the choice of time step size in this section. We again use the oscillatory Couette flow example with radial basis function reconstructions as in the previous section; however, we now vary the CFL number. Recall we defined the CFL number as $C_{\text{CFL}} = \frac{\Delta x}{\size{\uu}\Delta t}$ in which $\size{\uu}$ is the maximum velocity over both space and time so that we used a fixed time step size throughout the entire simulation. As before, we estimate the $C$ coefficient (see equation \eqref{eq:C}), and the values are reported in table \ref{tab:c_coeff_CFL}. We observe a decrease in the error as the CFL number increases. This is not surprising as the error term for the semi-Lagrangian method contains a term that scales like $\bigOh{\frac{1}{\Delta t}}$, so that the error can decrease for larger timesteps \cite{Fletcher2019}, although this reduction in error will eventually break down.

\begin{table}
\begin{center}
\begin{tabular}{c|cccccc}
\hline
& $q_\text{o}$ & & & & &\\
\hline
$C_\text{CFL}$ & 0.5 & 1.0 & 1.5 & 2.0 & 2.5 & 3.0 \\
\hline
$L^1$ norm      & 16.20 & 12.44 & 10.08 & 8.18 & 7.04 & 6.98 \\
$L^2$ norm      & 11.12 & 8.95  & 7.47  & 6.38 & 5.69 & 5.25 \\
$L^\infty$ norm & 12.88 & 11.81 & 10.55 & 9.43 & 8.46 & 7.61 \\
\hline
& $q_\text{v}$ & & & & &\\
\hline
$C_\text{CFL}$ & 0.5 & 1.0 & 1.5 & 2.0 & 2.5 & 3.0 \\
\hline
$L^1$ norm      & 3.44 & 3.01 & 2.86 & 2.64 & 2.40 & 2.18 \\
$L^2$ norm      & 1.60 & 1.47 & 1.38 & 1.24 & 1.08 & 0.95 \\
$L^\infty$ norm & 2.50 & 2.65 & 2.49 & 2.17 & 1.81 & 1.60 \\
\hline
\end{tabular}
\end{center}
\caption{The leading order error coefficient $C$ as estimated by equation \eqref{eq:C} as a function of the CFL number $C_\text{CFL}$. We see a decrease in error as the CFL number is increased.}\label{tab:c_coeff_CFL}
\end{table}

\subsection{Interaction of two concentrations in oscillatory Couette flow}

We now consider the case of two fluid phase chemicals advecting and diffusing on two different domains, but interacting through a common boundary. Specifically, we consider oscillatory Couette flow of two chemicals, one of which is contained within a vesicle that passively advects with the flow and is initially in the shape of a disk. The other chemical exists outside the vesicle, but between the disks defining the domain of interest. Specifically, we solve the equations
\begin{subequations} \label{eq:reaction_couette}
\begin{alignat}{2}
\frac{\partial q_\text{v}}{\partial t} + \uu\cdot \grad q_\text{v} &= D \Lap q_\text{v}, && \xx\in\Omega_\text{v} \label{eq:reaction_couette_v}\\
-D\frac{\partial q_\text{v}}{\partial \nn} &= \kappa\parens{q_\text{v} - q_\text{o}}, && \xx \in\Gamma_\text{v} \label{eq:reaction_couette_v_bc}\\
\frac{\partial q_\text{o}}{\partial t} + \uu\cdot \grad q_\text{o} &= D \Lap q_\text{o}, && \xx\in\Omega_\text{o} \label{eq:reaction_couette_o}\\
-D\frac{\partial q_\text{o}}{\partial \nn} &= \kappa\parens{q_\text{o} - q_\text{v}}, && \xx \in\Gamma_\text{v} \label{eq:reaction_couette_o_bc1}\\
-D\frac{\partial q_\text{o}}{\partial \nn} &= 0, && \xx \in\Gamma_\text{o}\setminus\Gamma_\text{v}, \label{eq:reaction_couette_o_bc2}
\end{alignat}
\end{subequations}
in which $\uu$ is defined in the previous section and $\Omega_\text{v}$ and $\Omega_\text{o}$ are the domains for the concentration inside and outside the vesicle, respectively. The initial concentration for $q_\text{v}$ is given by equation \eqref{eq:no_flux_inside}, and for $q_\text{o}$ is initially uniformly zero.

We modify the time discretization of the diffusion step to use a modified trapezoidal rule for the boundary conditions. We first solve equation \eqref{eq:diff_timestepping} for initial approximations $\tilde{q}_\text{o}$ and $\tilde{q}_\text{v}$ using explicit approximations for the boundary conditions
\begin{align}
-D\frac{\partial \tilde{q}_\text{o}}{\partial \nn} = \kappa\parens{\tilde{q}_\text{o} - q^n_\text{v}}, \text{ and } \\
-D\frac{\partial \tilde{q}_\text{v}}{\partial \nn} = \kappa\parens{\tilde{q}_\text{v} - q^n_\text{o}}.
\end{align}
We then solve equation \eqref{eq:diff_timestepping} for $q_\text{o}^{n+1}$ and $q_\text{v}^{n+1}$ again using the intermediate result when evaluating the boundary conditions
\begin{align}
-D\frac{\partial q^{n+1}_\text{o}}{\partial \nn} = \kappa\parens{q_\text{o}^{n+1} - \tilde{q}_\text{v}}, \text{ and } \\
-D\frac{\partial q^{n+1}_\text{v}}{\partial \nn} = \kappa\parens{q_\text{v}^{n+1} - \tilde{q}_\text{o}},
\end{align}
in which $\tilde{q}_\text{v}$ and $\tilde{q}_\text{o}$ are the results from the intermediate result.

\begin{figure}
\begin{center}
\phantomsubcaption\label{fig:reacting_couette:0}\phantomsubcaption\label{fig:reacting_couette:1}
\phantomsubcaption\label{fig:reacting_couette:2}\phantomsubcaption\label{fig:reacting_couette:3}
\phantomsubcaption\label{fig:reacting_couette:4}
\psfragscanon
\psfrag{3}{$0.30$}
\psfrag{225}{$0.225$}
\psfrag{15}{$0.15$}
\psfrag{075}{$0.075$}
\psfrag{0}{$0.0$}
\psfrag{8}{$0.8$}
\psfrag{6}{$0.6$}
\psfrag{4}{$0.4$}
\psfrag{2}{$0.2$}
\subfigimg[width=0.45\linewidth,hsep=0em,vsep=0em,pos=ul]{(a)}{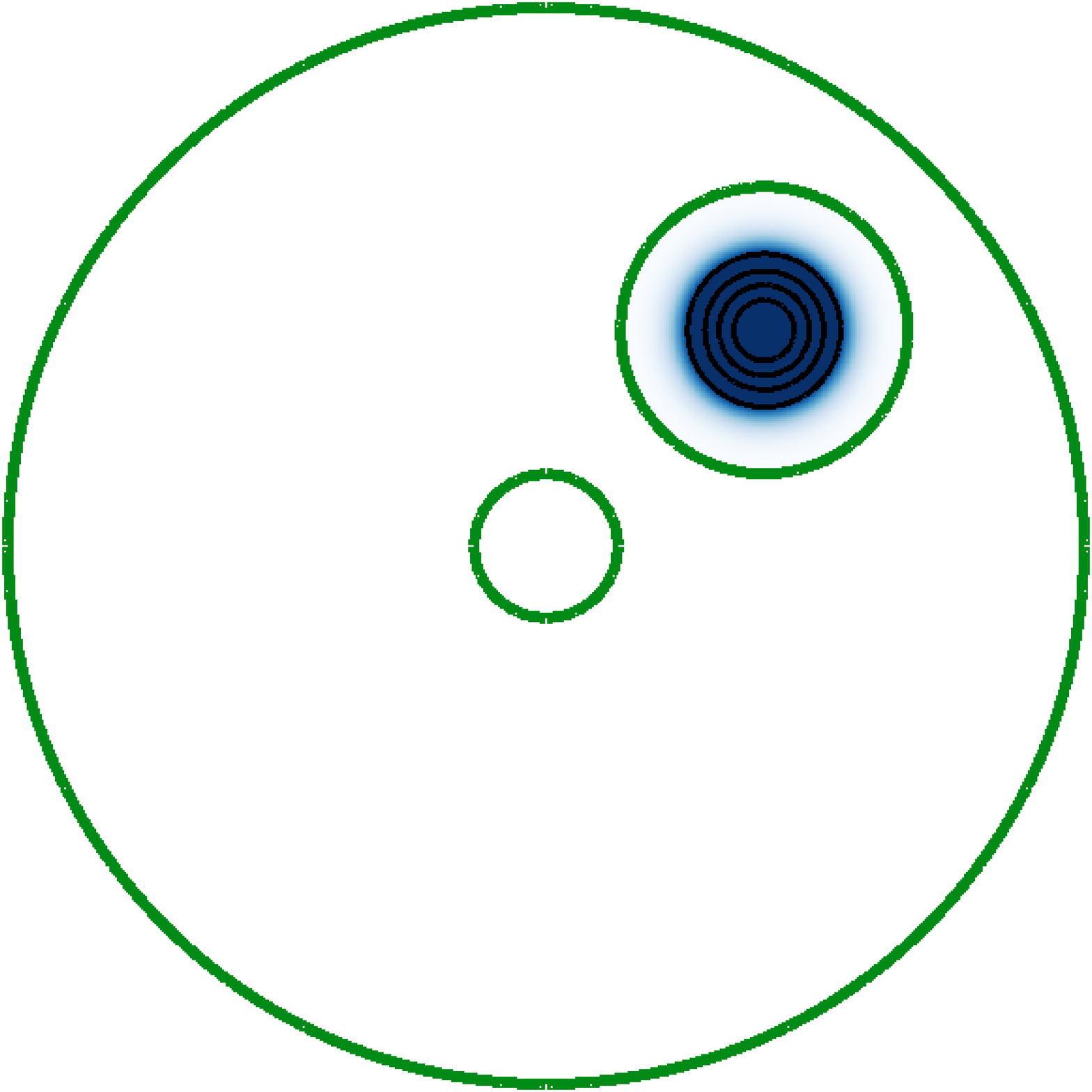}
\subfigimg[width=0.45\linewidth,hsep=0em,vsep=0em,pos=ul]{(b)}{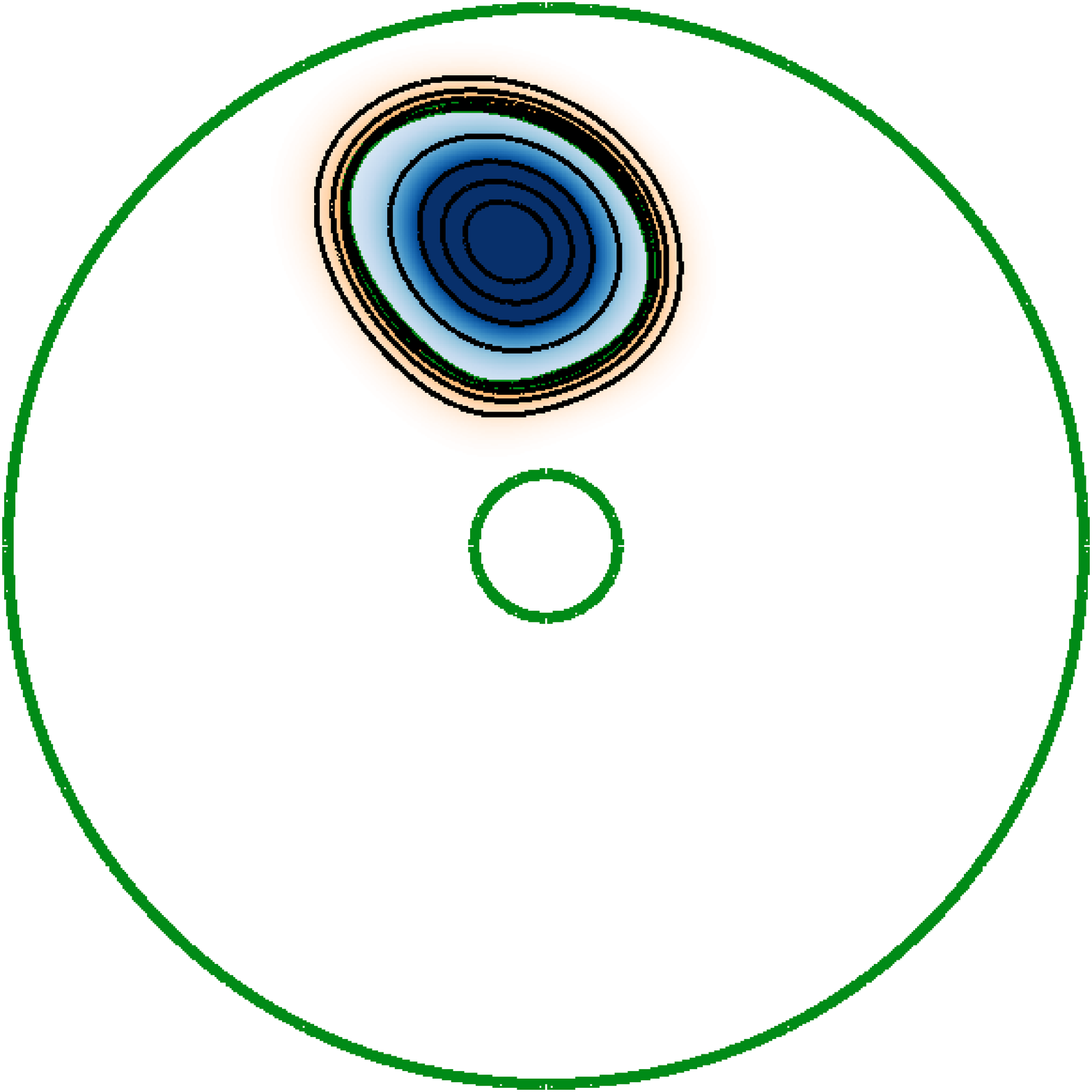}\\
\subfigimg[width=0.45\linewidth,hsep=0em,vsep=0em,pos=ul]{(c)}{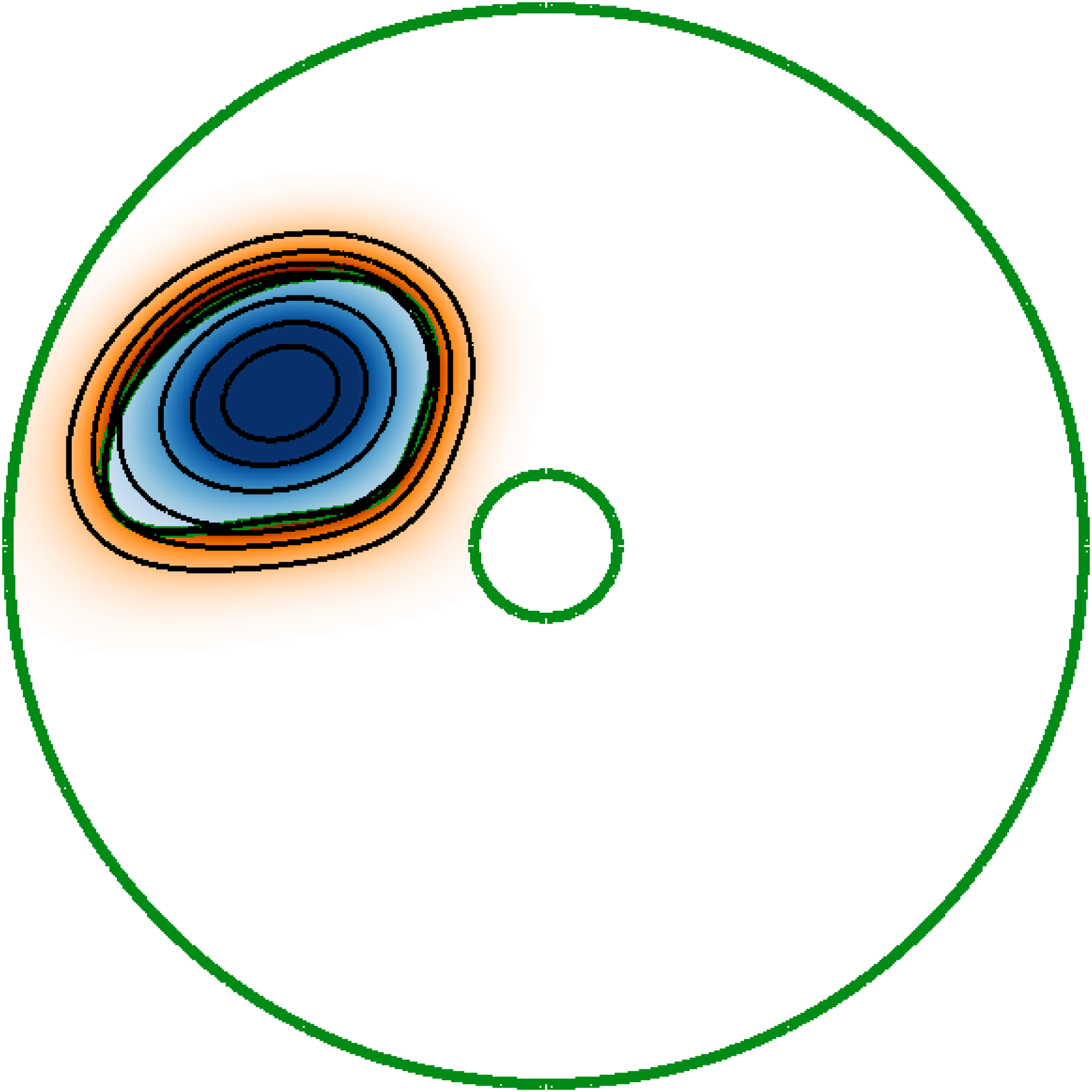}
\subfigimg[width=0.45\linewidth,hsep=0em,vsep=0em,pos=ul]{(d)}{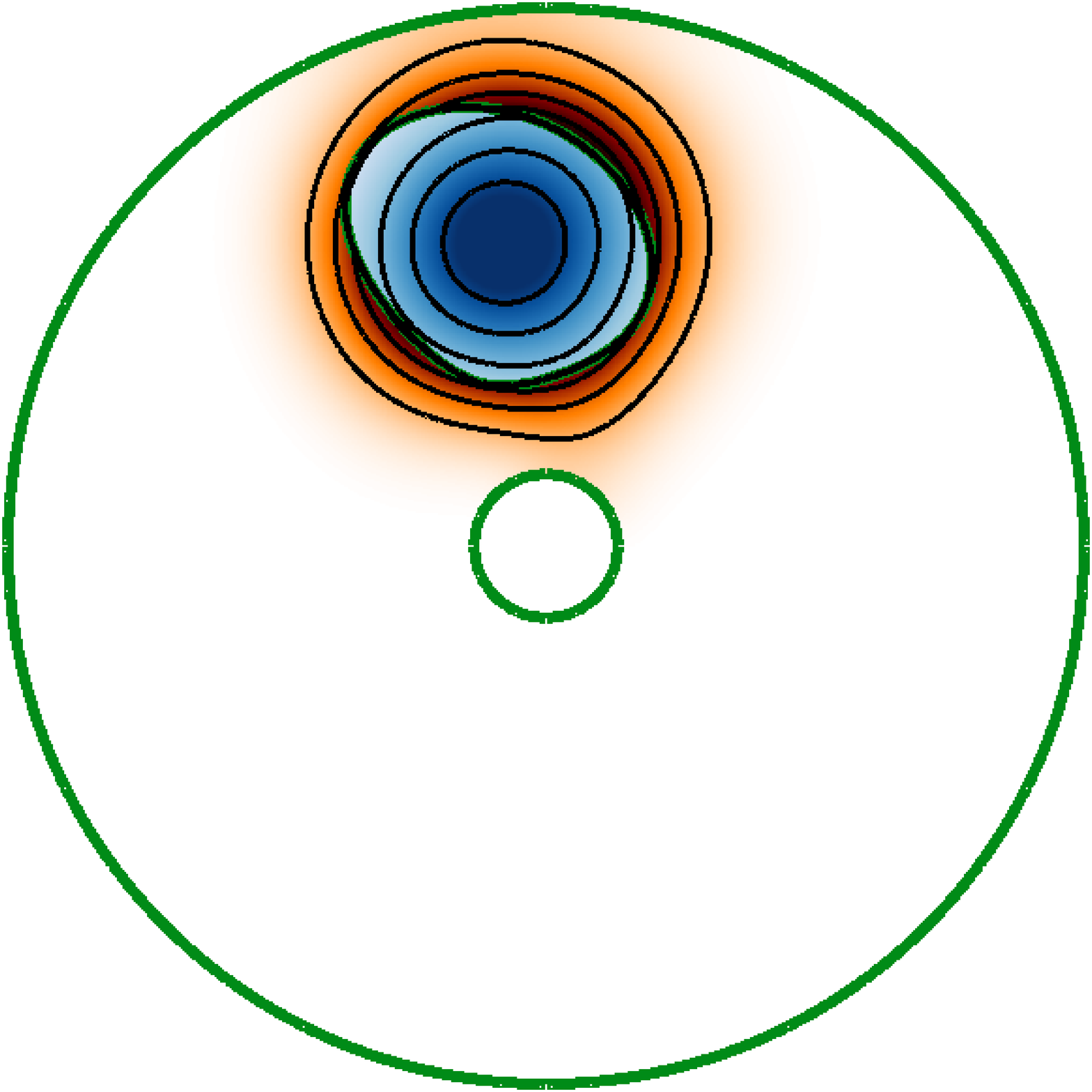}\\
\subfigimg[width=0.45\linewidth,hsep=0em,vsep=0em,pos=ul]{(e)}{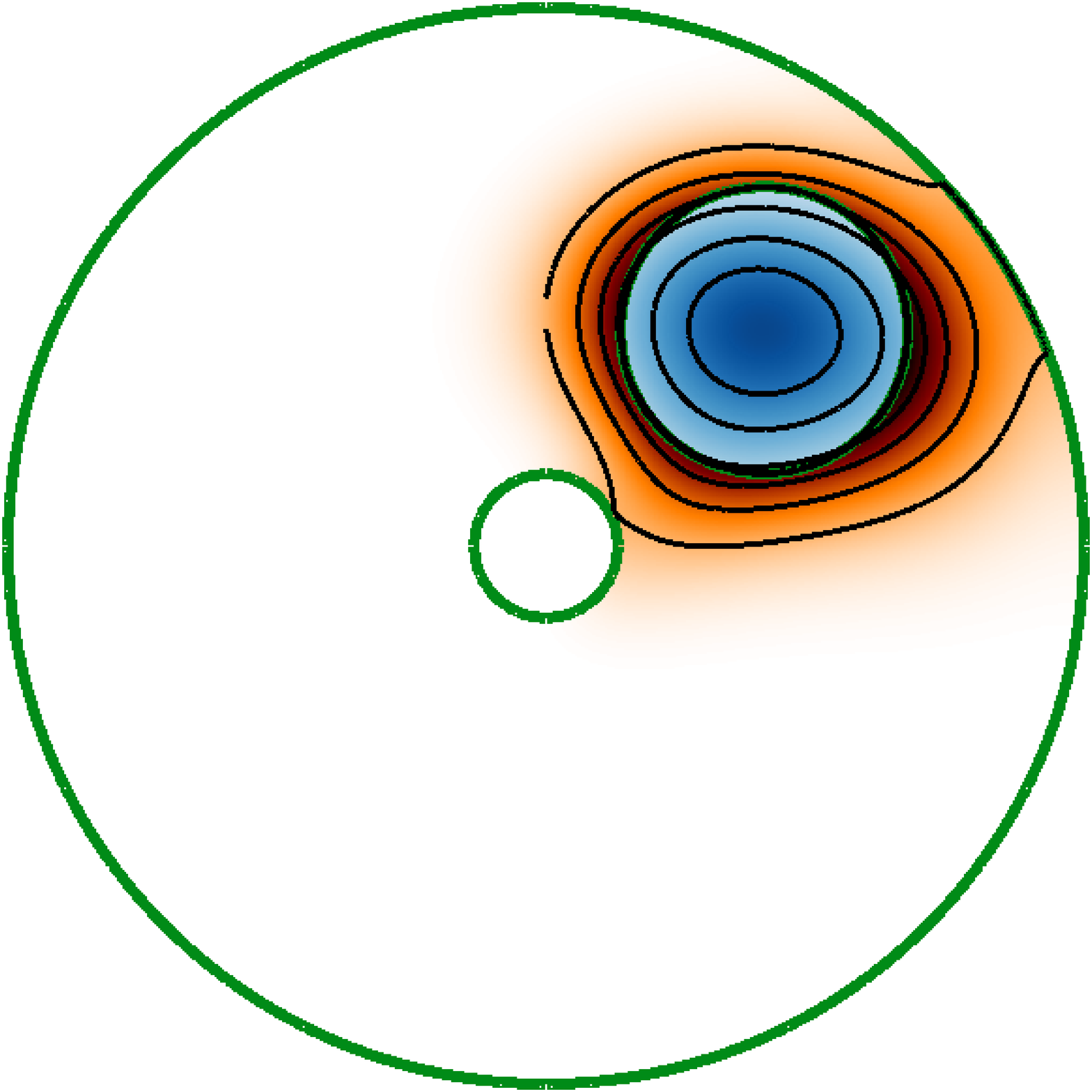}
\includegraphics[width=0.2\linewidth]{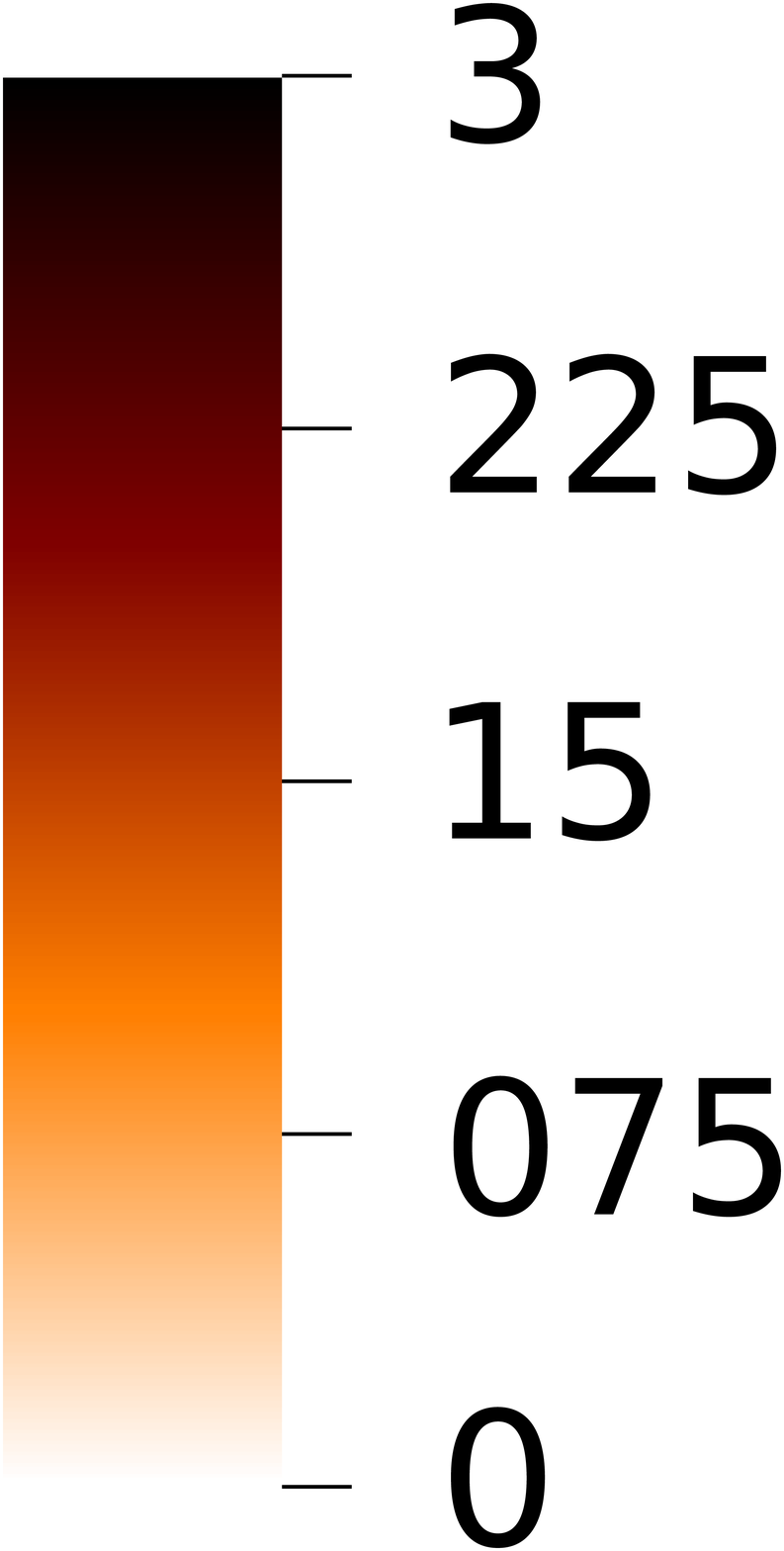}
\includegraphics[width=0.2\linewidth]{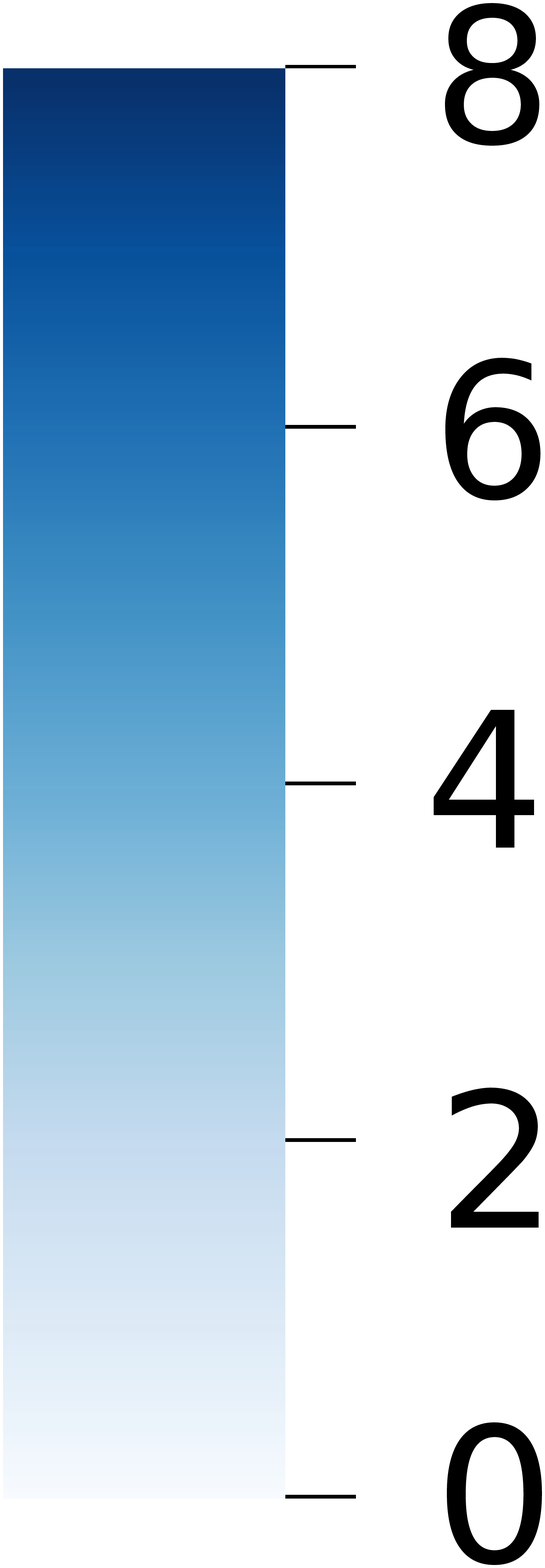}
\end{center}
\caption{Concentrations $q_\text{o}$ (orange) and $q_\text{v}$ (blue) at times $0.0$ (\subref{fig:reacting_couette:0}), $0.5$ (\subref{fig:reacting_couette:1}), $1.0$ (\subref{fig:reacting_couette:2}), $1.5$ (\subref{fig:reacting_couette:3}), and $2.0$ (\subref{fig:reacting_couette:4}). The zero contour of the level set $\phi_\text{o}$ is shown in green. The diffusion coefficient is fixed at $D = 0.05$ with $1024$ points in each direction. We use a fixed time-step corresponding to a maximum CFL number of $C_\text{CFL} = 0.5$ to a final time of $T = 2$.}
\label{fig:reacting_couette}
\end{figure}

For this example, we perform only RBF reconstructions at the boundaries. We set $\kappa = 1$. The remaining parameters are the same as in the previous section. Concentration values for various points in time are shown in Figure \ref{fig:reacting_couette}. As $q_\text{v}$ reaches the boundary $\Gamma_\text{v}$, it is converted to $q_\text{o}$ at a rate proportional to the difference in concentrations.

\begin{figure}
\begin{center}
\phantomsubcaption\label{fig:reaction_in_converge:l1_rates}\phantomsubcaption\label{fig:reaction_in_converge:l2_rates}
\phantomsubcaption\label{fig:reaction_in_converge:max_rates}\phantomsubcaption\label{fig:reaction_in_converge:bdry_rates}
\subfigimg[width=0.45\textwidth,hsep=-1em,vsep=2.5em,pos=ul]{(a)}{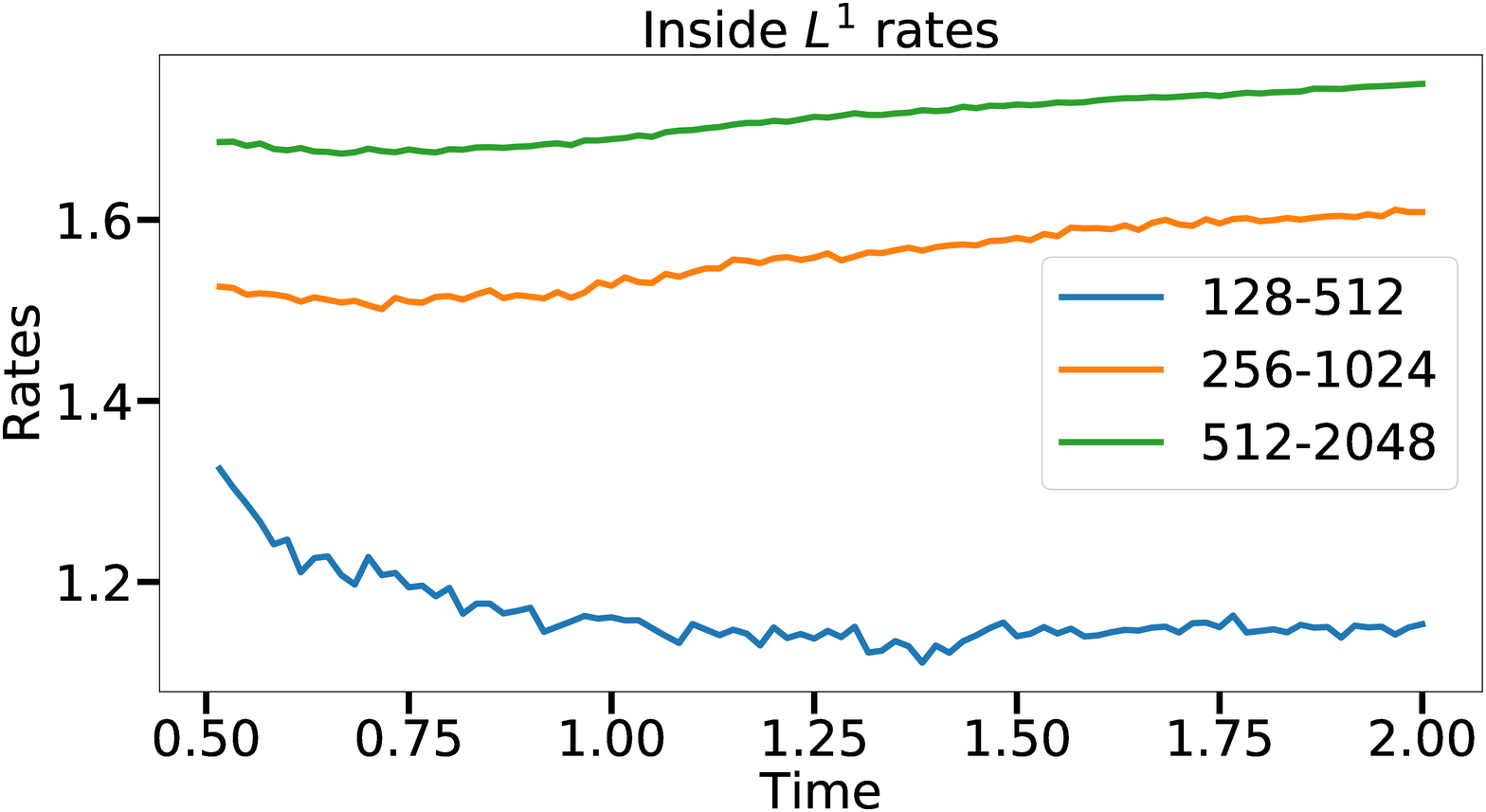}
\subfigimg[width=0.45\textwidth,hsep=-1em,vsep=2.5em,pos=ul]{(b)}{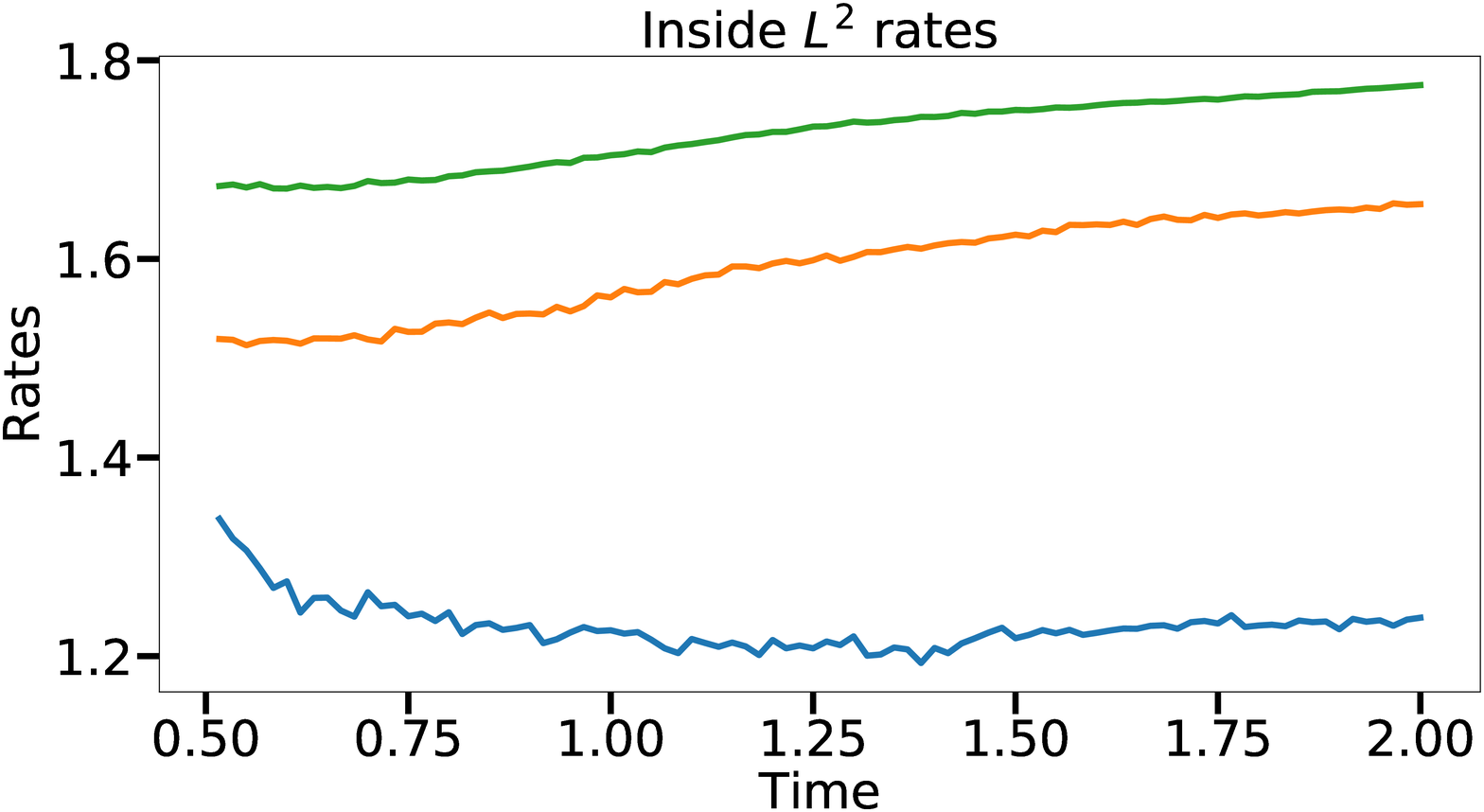}\\
\subfigimg[width=0.45\textwidth,hsep=-1em,vsep=2.5em,pos=ul]{(c)}{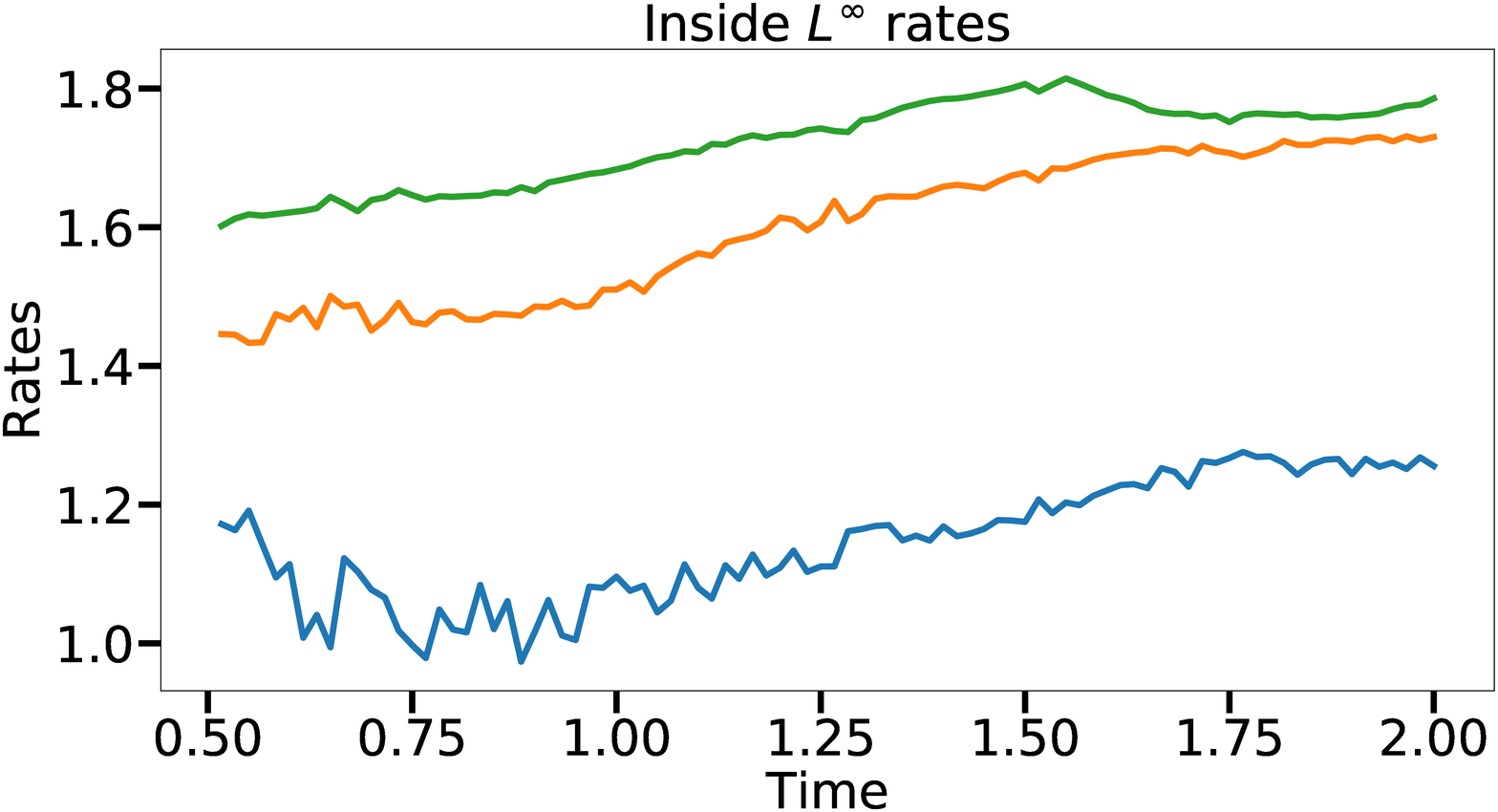}
\subfigimg[width=0.45\textwidth,hsep=-1em,vsep=2.5em,pos=ul]{(d)}{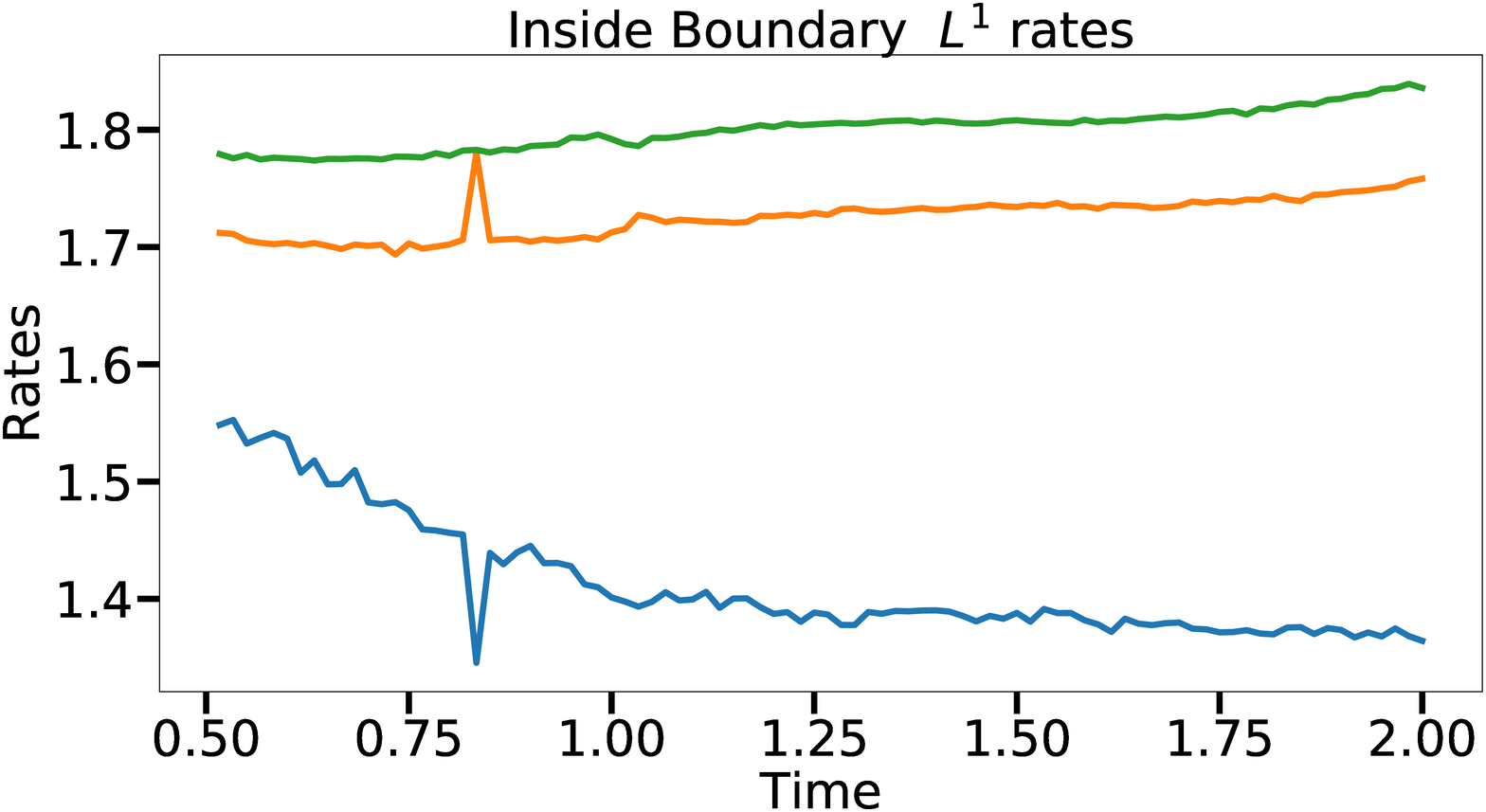}
\end{center}
\caption{Numerical convergence study for the interior of the disk in the $L^1$ (\subref{fig:reaction_in_converge:l1_rates}), $L^2$ (\subref{fig:reaction_in_converge:l2_rates}), and $L^\infty$ (\subref{fig:reaction_in_converge:max_rates}) norms. The convergence rates are computed from simulations with $N = 128,256, \text{ and } 512$ points (blue), $N = 256,512,\text{ and } 1024$ points (orange), and $N = 512,1024, \text{ and } 2048$ points (green). Also shown are the convergence rates for the extrapolation of the solution to the boundary (\subref{fig:reaction_in_converge:bdry_rates}).}
\label{fig:reaction_in_converge}
\end{figure}

\begin{figure}
\begin{center}
\phantomsubcaption\label{fig:reaction_out_converge:l1_rates}\phantomsubcaption\label{fig:reaction_out_converge:l2_rates}
\phantomsubcaption\label{fig:reaction_out_converge:max_rates}
\subfigimg[width=0.45\textwidth,hsep=-1em,vsep=2em,pos=ul]{(a)}{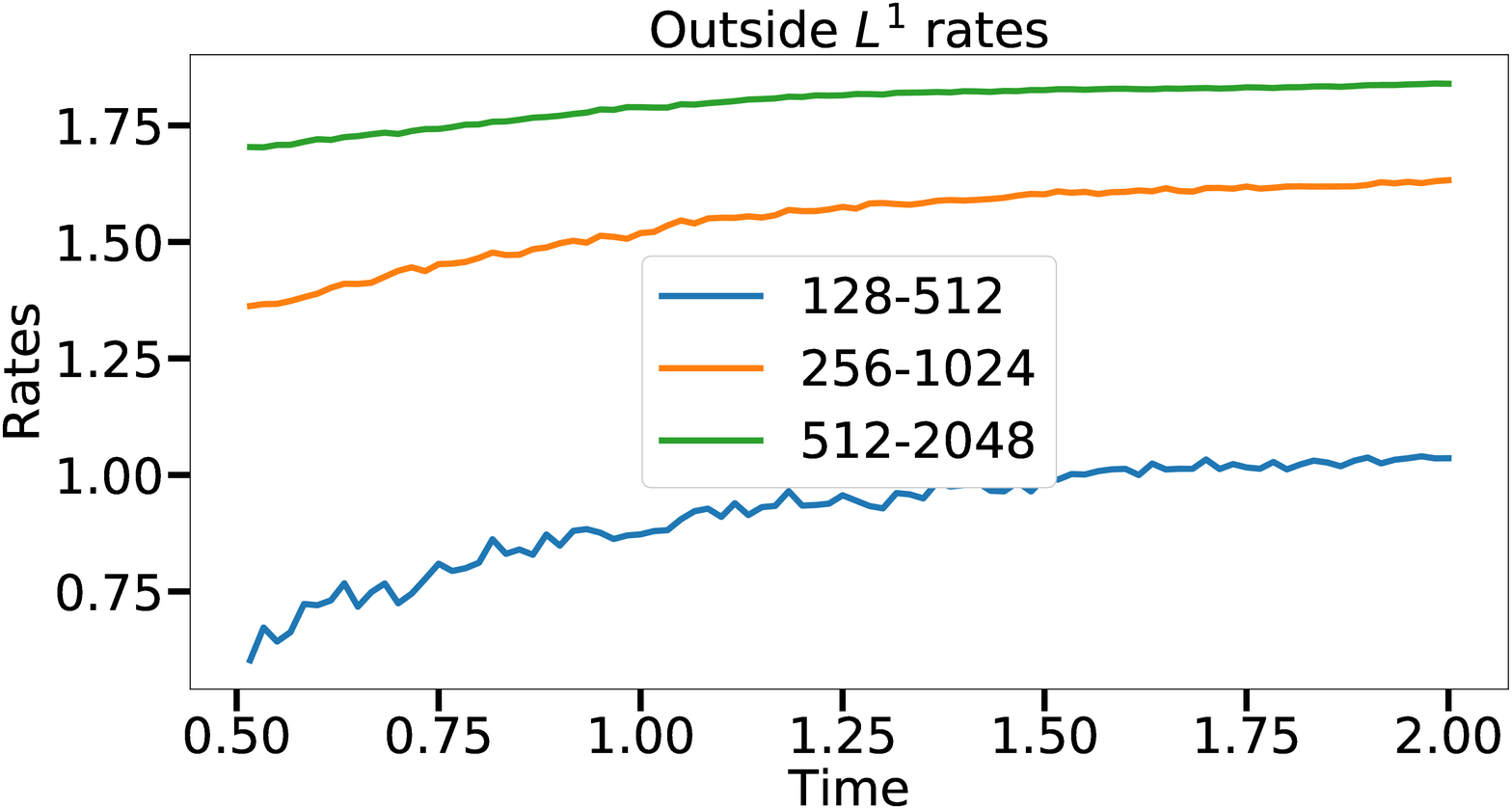}
\subfigimg[width=0.45\textwidth,hsep=-1em,vsep=2em,pos=ul]{(b)}{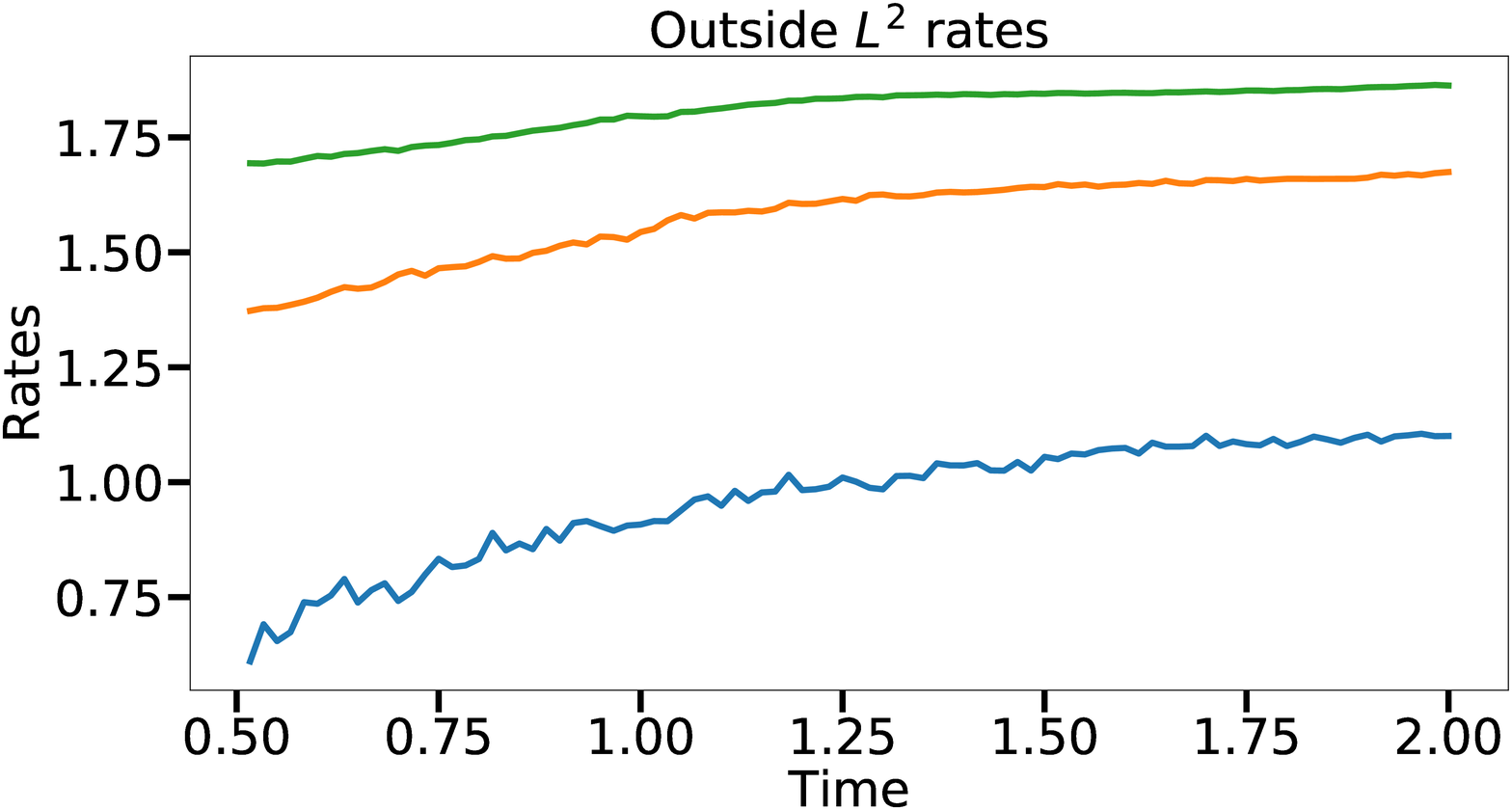}
\subfigimg[width=0.45\textwidth,hsep=-1em,vsep=2em,pos=ul]{(c)}{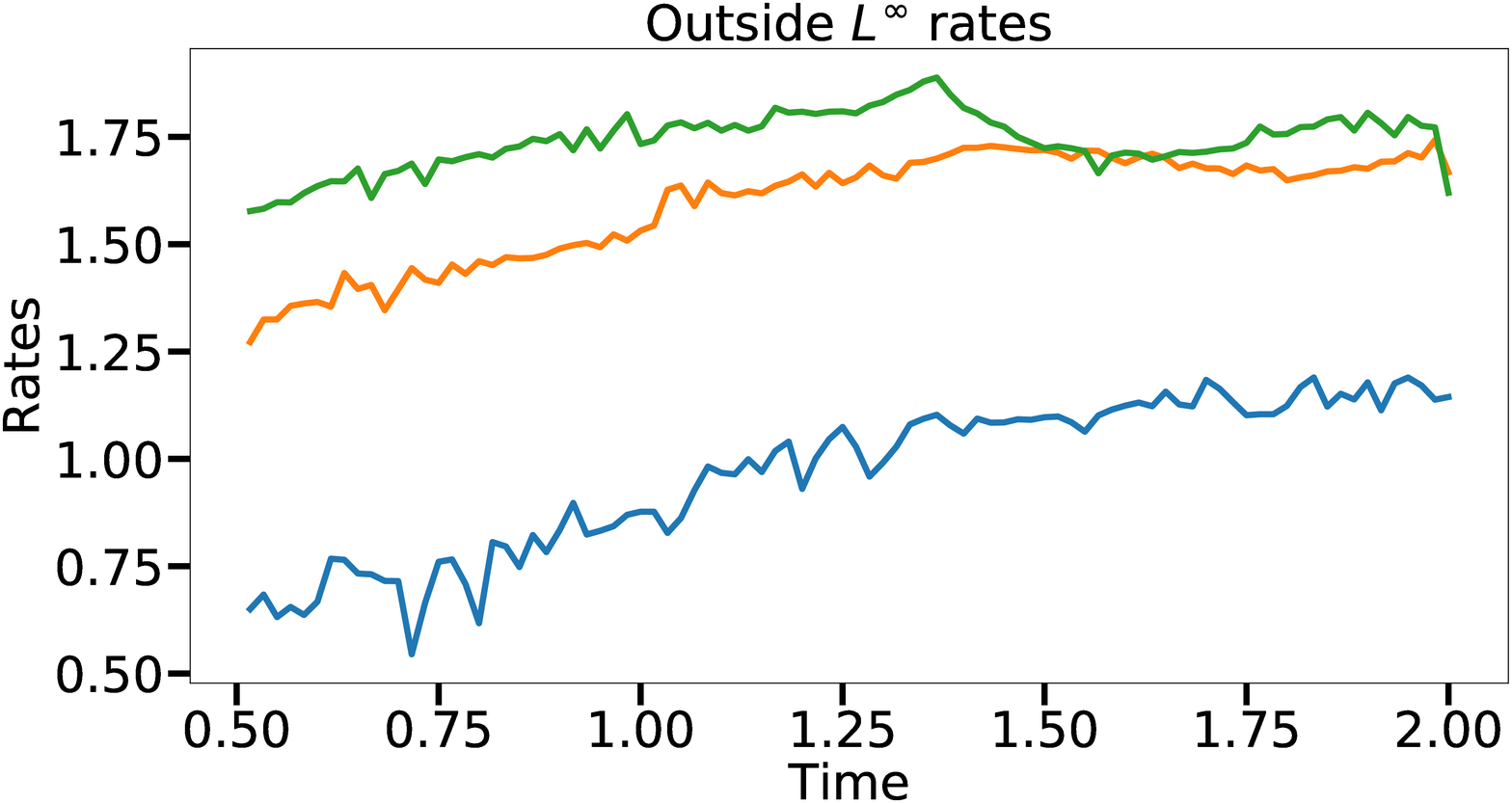}
\end{center}
\caption{Numerical convergence study for the exterior of the disk but inside the two cylinders in the $L^1$ (\subref{fig:reaction_out_converge:l1_rates}), $L^2$ (\subref{fig:reaction_out_converge:l2_rates}), and $L^\infty$ (\subref{fig:reaction_out_converge:max_rates}) norms. The convergence study is computed from simulations with $N = 128,256,\text{ and }512$ points (blue), $N = 256,512, \text{ and } 1024$ points (orange), and $N = 512,1024,\text{ and } 2048$ points (green).}
\label{fig:reaction_out_converge}
\end{figure}

\begin{figure}
\begin{center}
\phantomsubcaption\label{fig:reaction_in_norms:l1_norms}\phantomsubcaption\label{fig:reaction_in_norms:l2_norms}
\phantomsubcaption\label{fig:reaction_in_norms:max_norms}\phantomsubcaption\label{fig:reaction_in_norms:bdry_norms}
\subfigimg[width=0.45\textwidth,hsep=-1em,vsep=2em,pos=ul]{(a)}{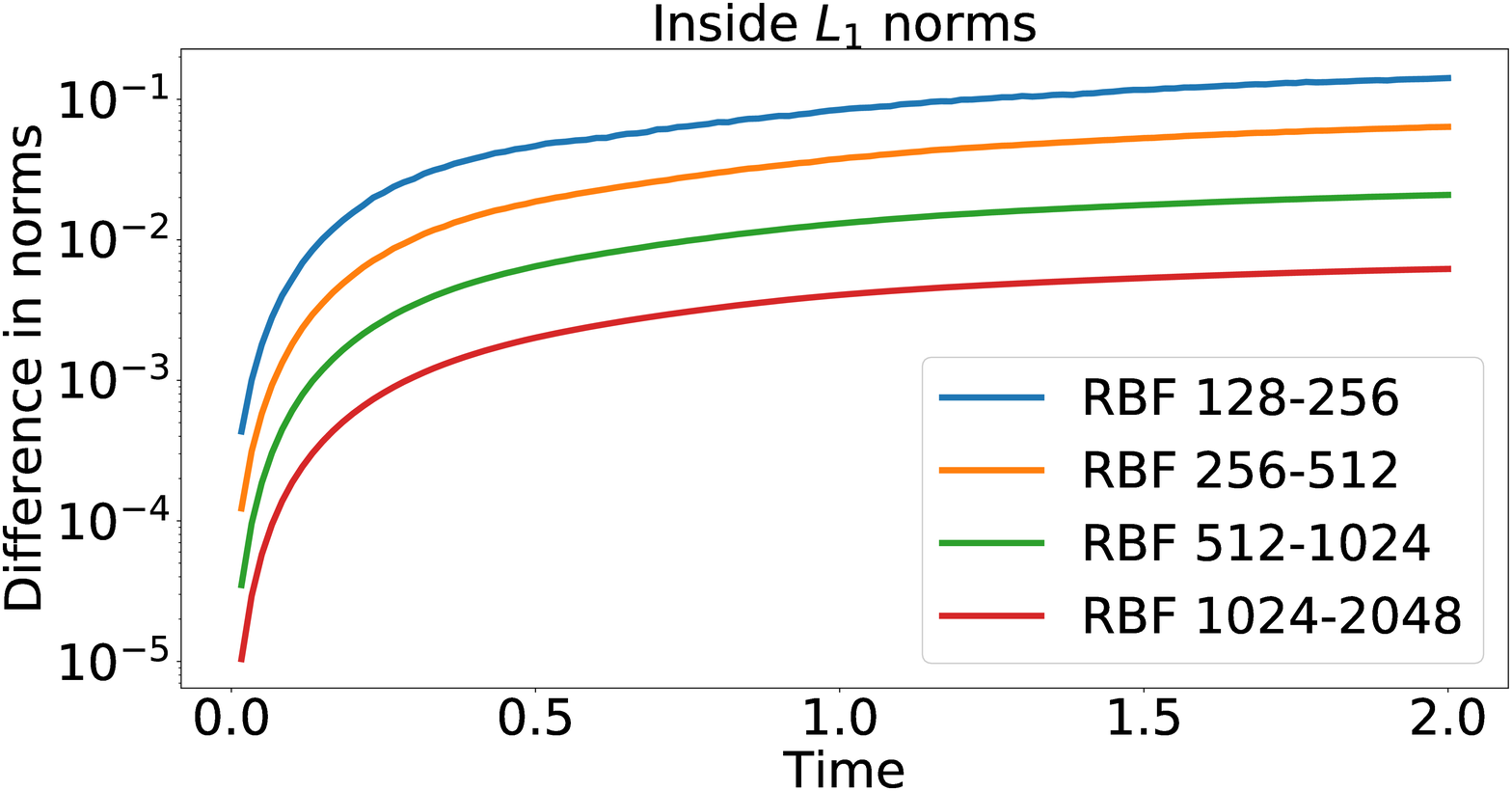}
\subfigimg[width=0.45\textwidth,hsep=-1em,vsep=2em,pos=ul]{(b)}{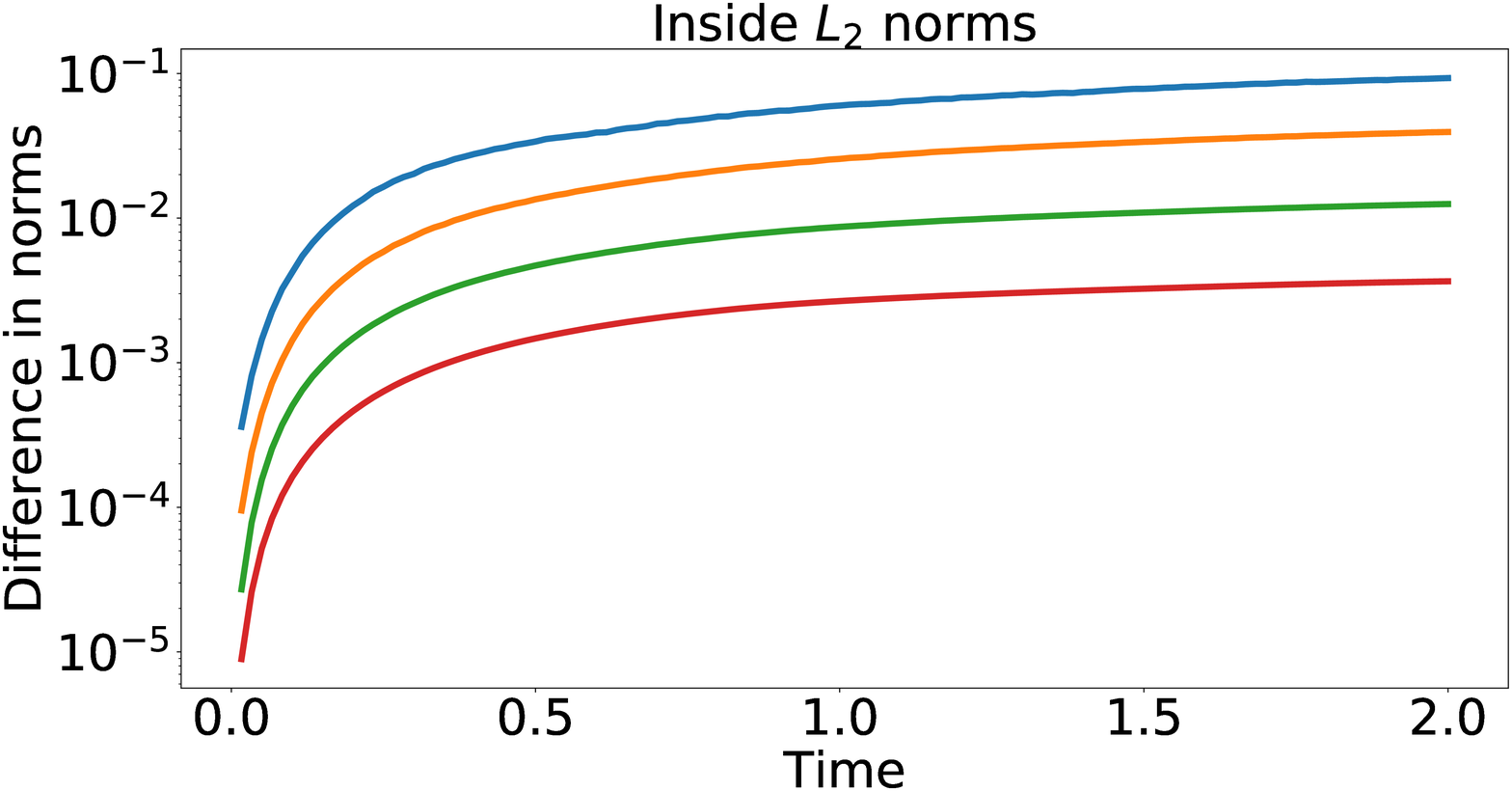}\\
\subfigimg[width=0.45\textwidth,hsep=-1em,vsep=2em,pos=ul]{(c)}{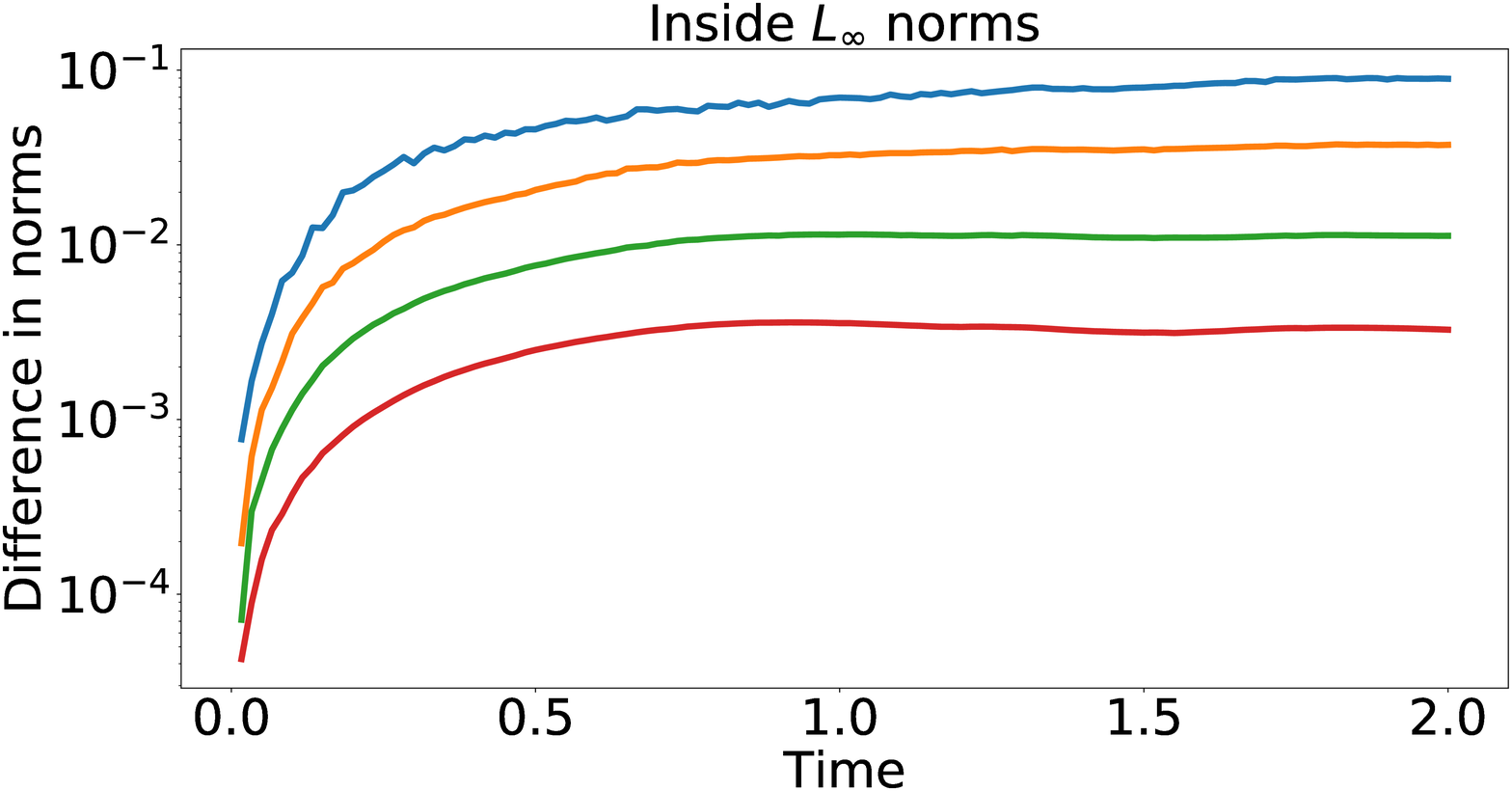}
\subfigimg[width=0.45\textwidth,hsep=-1em,vsep=2em,pos=ul]{(d)}{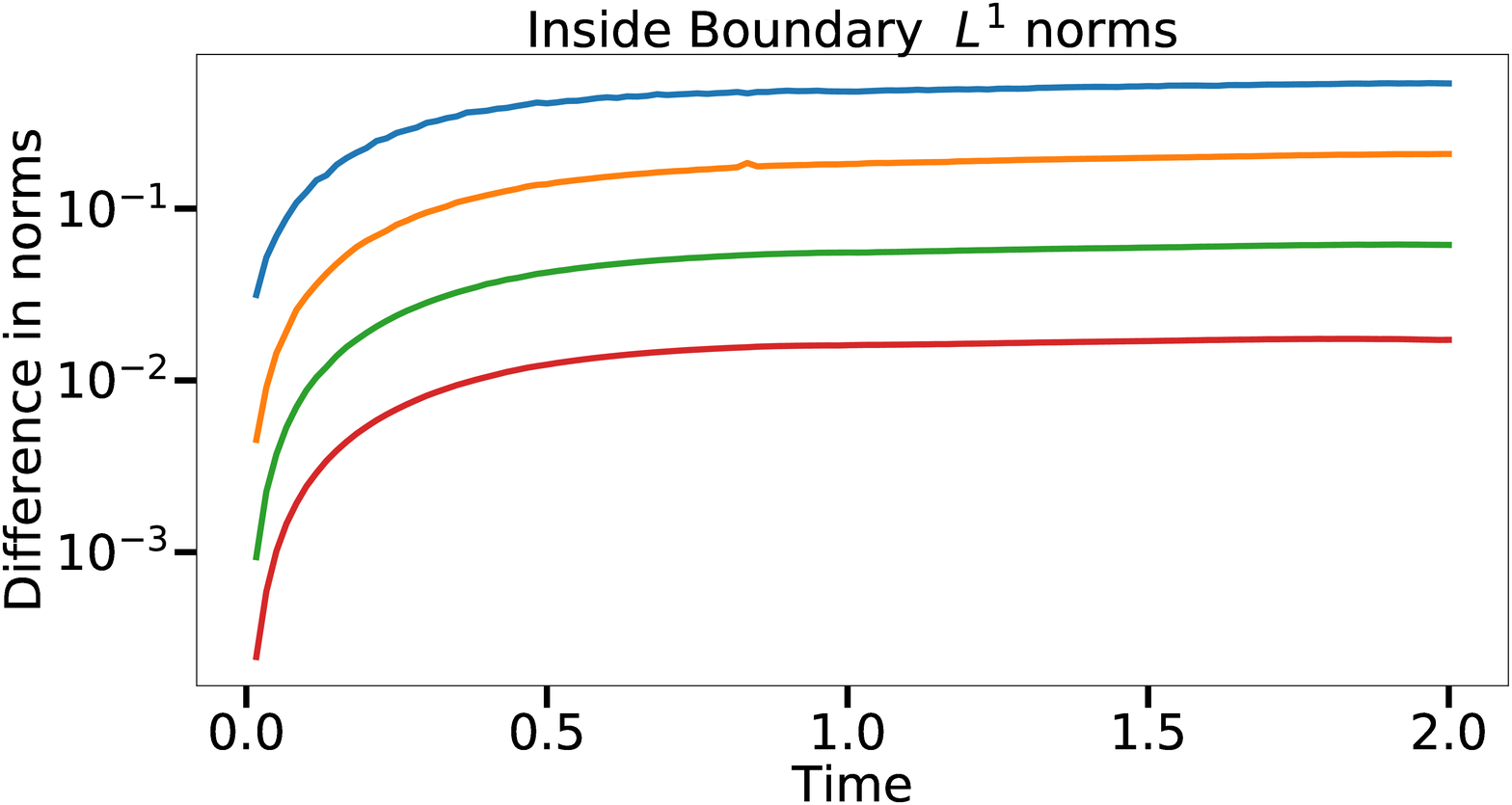}
\end{center}
\caption{Norms of solution differences for the concentration field inside the domain for the $L^1$ (\subref{fig:reaction_in_norms:l1_norms}), $L^2$ (\subref{fig:reaction_in_norms:l2_norms}), and $L^\infty$ (\subref{fig:reaction_in_norms:max_norms}) norms. Also shown are the norms of the differences between solutions extrapolated to the boundary (\subref{fig:reaction_in_norms:bdry_norms}).}
\label{fig:reaction_in_norms}
\end{figure}

\begin{figure}
\begin{center}
\phantomsubcaption\label{fig:reaction_out_norms:l1_norms}\phantomsubcaption\label{fig:reaction_out_norms:l2_norms}
\phantomsubcaption\label{fig:reaction_out_norms:max_norms}
\subfigimg[width=0.45\textwidth,hsep=-1em,vsep=2em,pos=ul]{(a)}{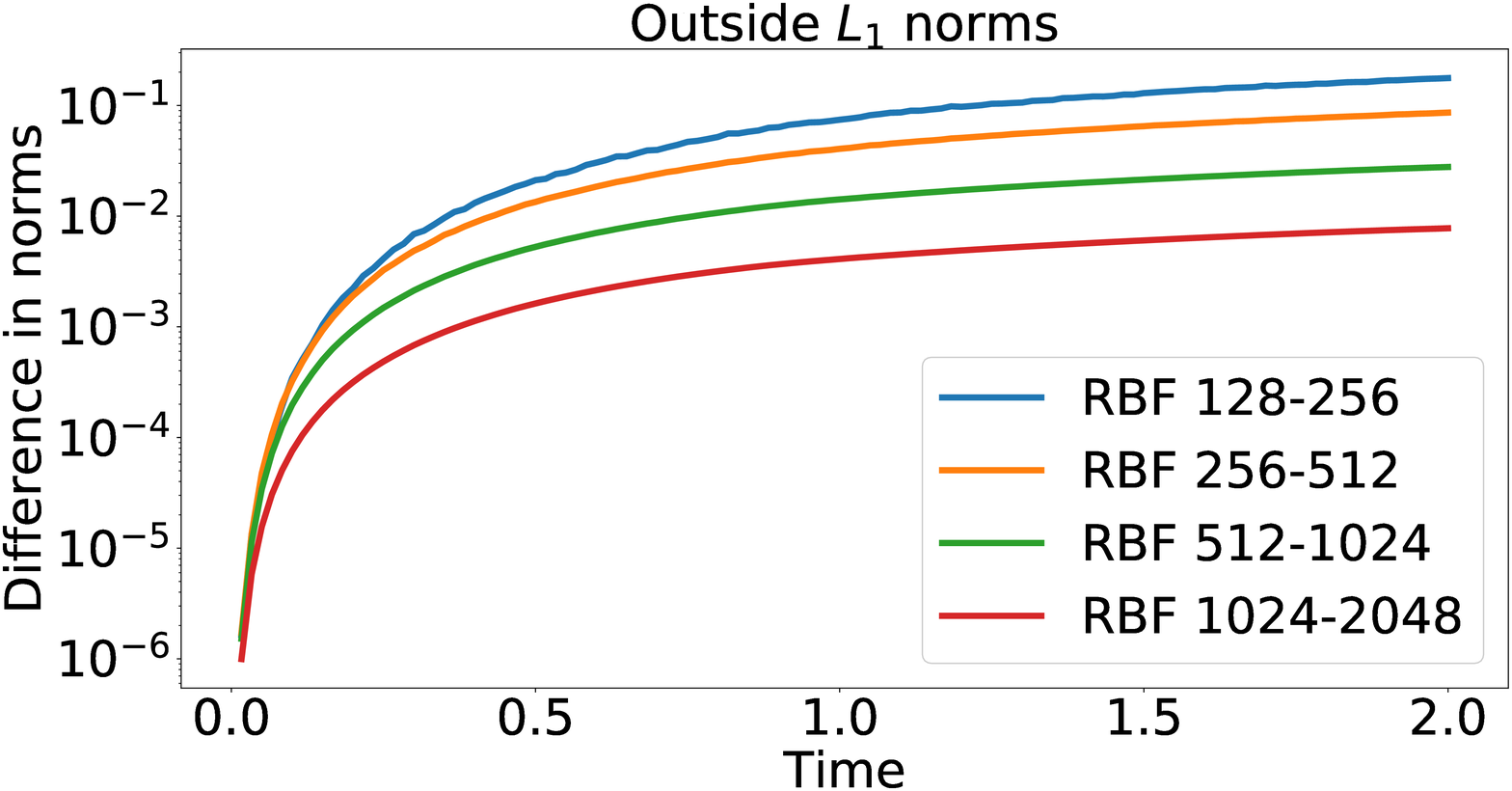}
\subfigimg[width=0.45\textwidth,hsep=-1em,vsep=2em,pos=ul]{(b)}{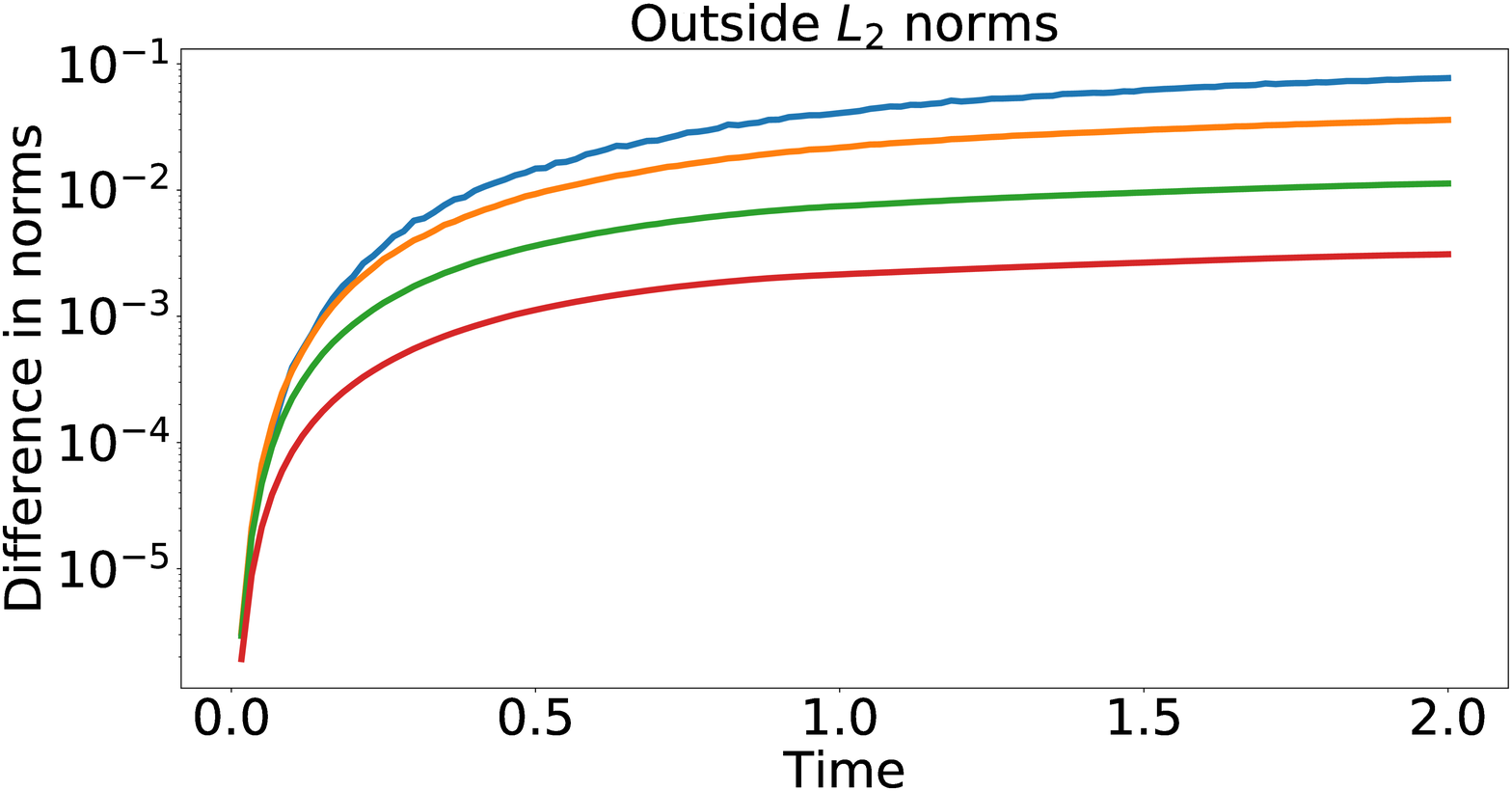}
\subfigimg[width=0.45\textwidth,hsep=-1em,vsep=2em,pos=ul]{(c)}{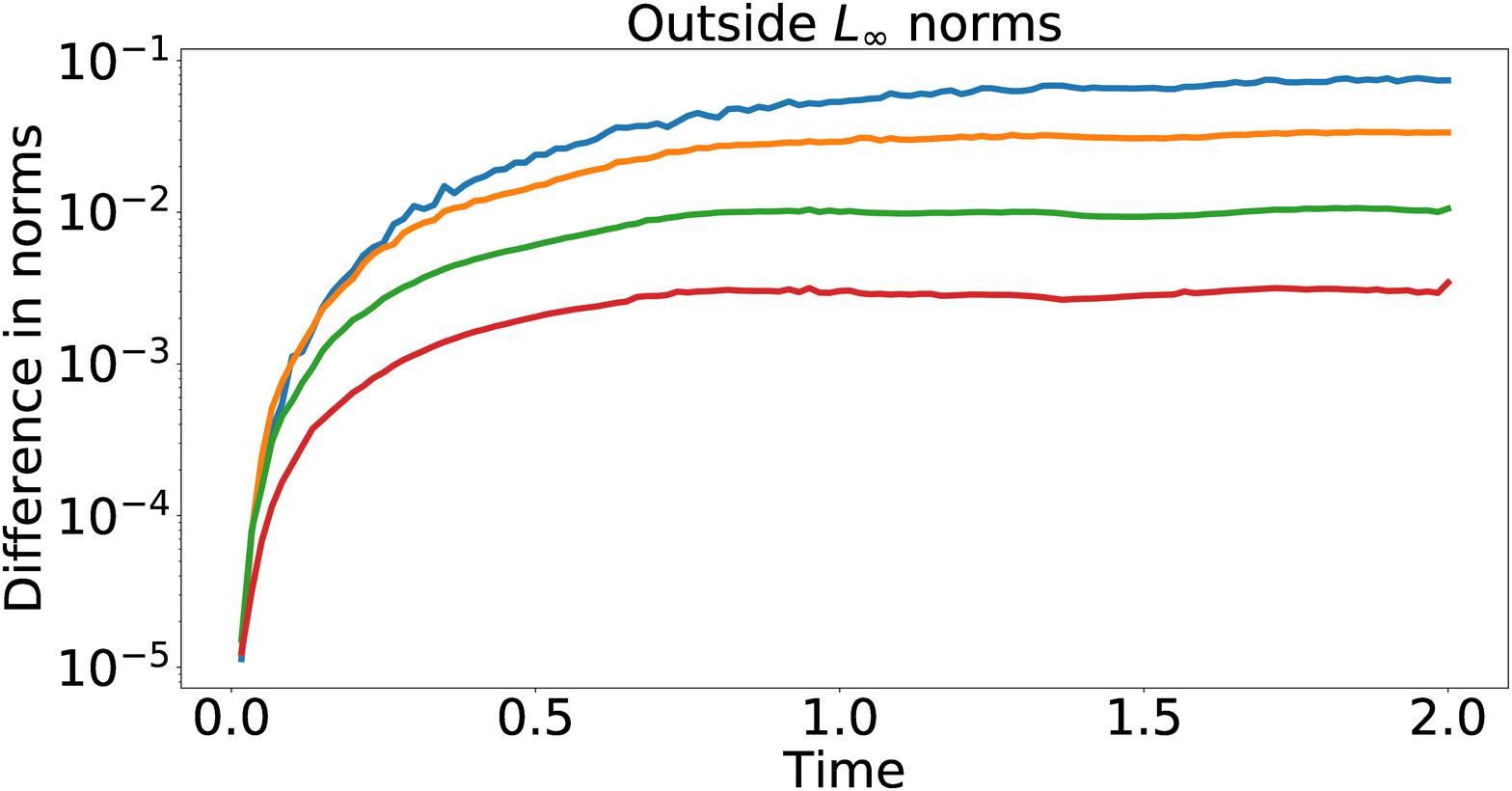}
\end{center}
\caption{Norms of solution differences for the concentration field outside the domain for the $L^1$ (\subref{fig:reaction_out_norms:l1_norms}), $L^2$ (\subref{fig:reaction_out_norms:l2_norms}), and $L^\infty$ (\subref{fig:reaction_out_norms:max_norms}) norms.}
\label{fig:reaction_out_norms}
\end{figure}

As in the previous section, we perform a numerical convergence study for both the interior and exterior fields. Figures \ref{fig:reaction_in_converge} and \ref{fig:reaction_out_converge} show the convergence rates as a function of time. The convergence rates as we refine the grid appear to be converging toward second order, although the method has not yet settled down into it's asymptotic regime. We also perform a convergence study for the solution extrapolated to the boundary, which is shown in Figure \ref{fig:reaction_in_converge:bdry_rates}. We again see approximate second order rates for this reconstruction. A plot of the norms of the solution differences is shown in Figures \ref{fig:reaction_in_norms} and \ref{fig:reaction_out_norms}. The differences in the finest grid show a difference of less than one percent at the final time.

\begin{figure}
\begin{center}
\phantomsubcaption\label{fig:total_amount:0}\phantomsubcaption\label{fig:total_amount:1}
\subfigimg[width=0.45\linewidth,hsep=-1em,vsep=2em,pos=ul]{(a)}{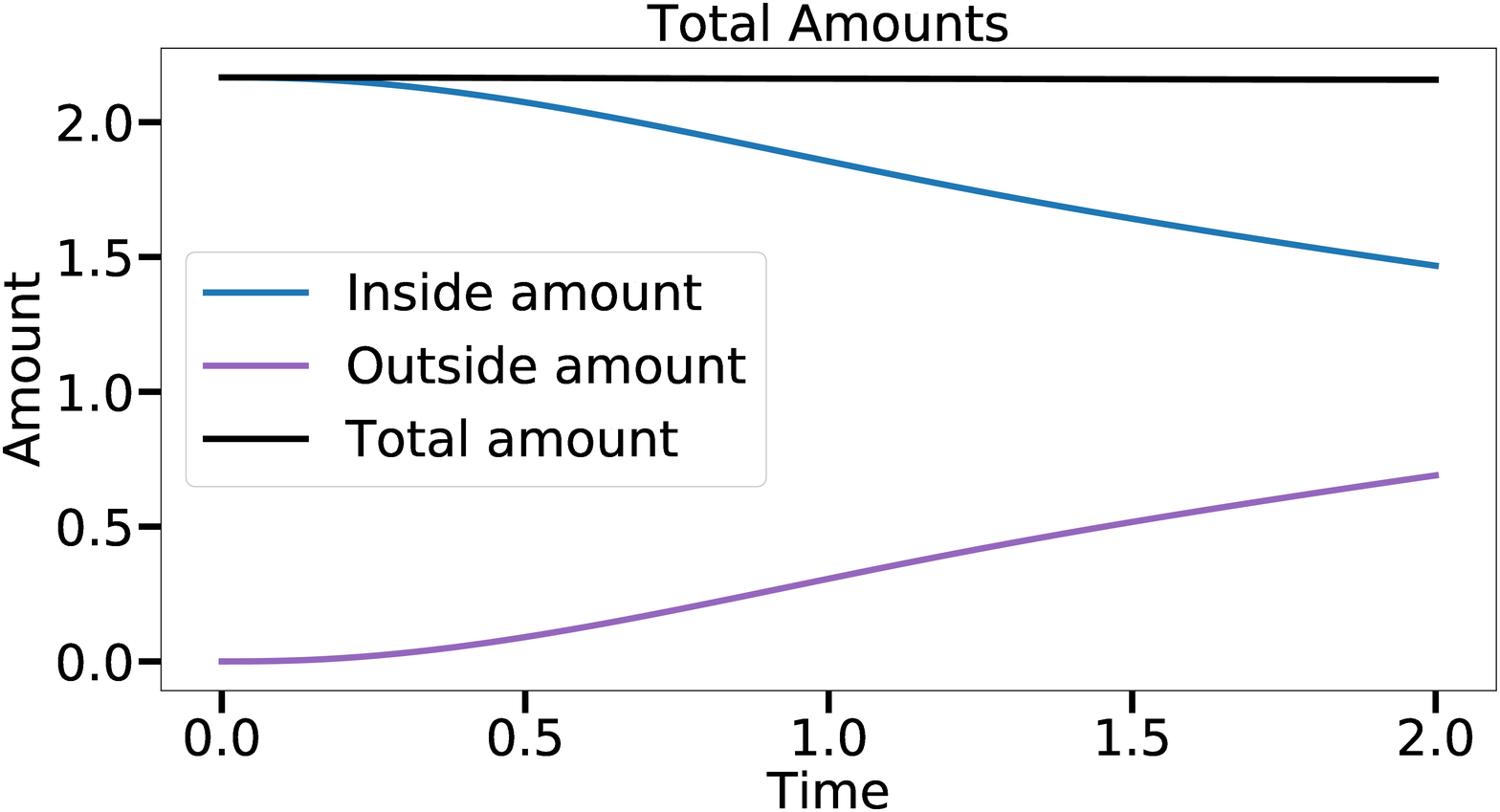}
\subfigimg[width=0.45\linewidth,hsep=-1em,vsep=2em,pos=ul]{(b)}{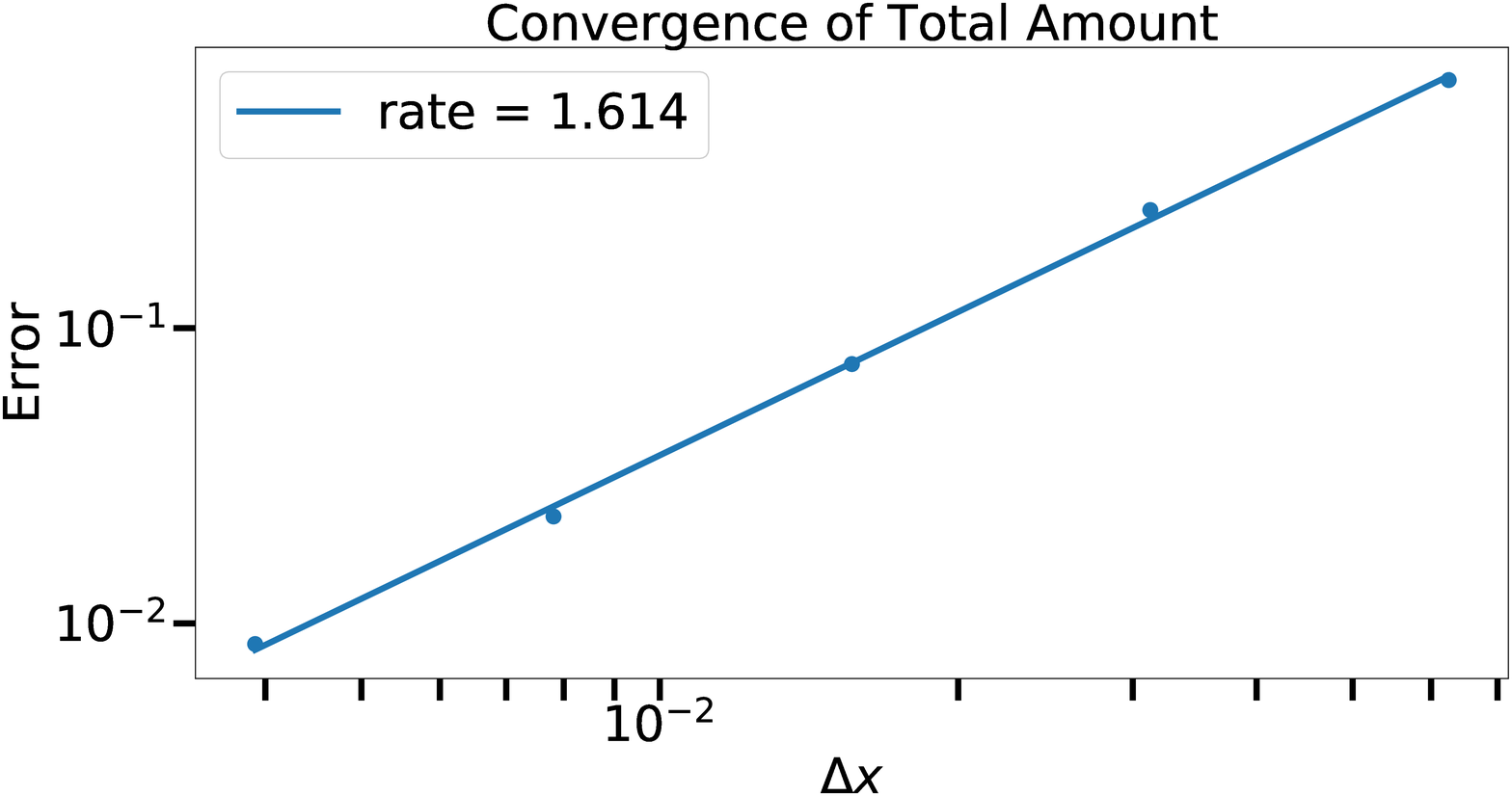}
\end{center}
\caption{Total concentration of both $q_\text{v}$ and $q_\text{o}$ for interactions with a oscillatory Couette flow (\subref{fig:total_amount:0}). The total amount converges at a rate that is between first and second order (\subref{fig:total_amount:1}).}
\label{fig:total_amount}
\end{figure}

To assess conservation, we perform a convergence study on the total amount of the two concentration fields. Figure \ref{fig:total_amount} shows the total amounts of both fields as a function of time, while the sum of both amounts remains relatively constant. Despite this method not being conservative, we see a loss of less than one percent over the course of the finest simulation. The total amount of concentration converges at a rate that is between first and second order accuracy.

\section{Conclusion}

We have presented a numerical method to simulate advection-diffusion problems with Robin boundary conditions on irregular, evolving domains. The method shows second order convergence in the $L^1$, $L^2$, and $L^\infty$ norms. We have also demonstrated the ability to accurately reconstruct solution values on the boundary, achieving second order accurate results. We have further demonstrated the method to be able to capture transformation of one concentration into another concentration with all interaction mediated through boundary conditions. Although only two dimensional tests are presented, we expect this method to be easily extended to three spatial dimensions. The use of radial basis functions when compared to moving least squares gives optimal convergence rates and more accurate solutions for the grid spacings considered. Z-splines were used in this study for computational efficiency, although any suitable reconstruction procedure should work. This method has many possible applications, including models with chemical transport in evolving bodies, such as esophageal transport, oxygen flow in the lungs, and blood flow.

\section*{Acknowledgements}
The authors thank Dr.\ Varun Shankar for helpful discussions on both radial basis functions and semi-Lagrangian schemes. This research was supported through NHLBI award 5U01HL143336 and NSF awards DMS 1664645, OAC 1450327, and OAC 1931516.
\printbibliography[heading=bibintoc]

\end{document}